\documentclass[11pt,reqno]{amsart}

\usepackage{amssymb}
\usepackage{geompsfi}

\catcode`\@=11

\long\def\@savemarbox#1#2{\global\setbox#1\vtop{\hsize\marginparwidth 
  \@parboxrestore\tiny\raggedright #2}}
\newcommand\lref[1]{\ref{#1}%
\@ifundefined{r@DisplaY #1}{}{ (#1)}}

\newcommand\fakelabel[2]{\@bsphack\if@filesw {\let\thepage\relax
   \newcommand\protect{\noexpand\noexpand\noexpand}%
\xdef\@gtempa{\write\@auxout{\string
      \newlabel{#1}{{#2}{\thepage}}}}}\@gtempa
   \if@nobreak \ifvmode\nobreak\fi\fi\fi\@esphack}

\catcode`\@=12

\def\Empty{}
\newcommand\oplabel[1]{
  \def\OpArg{#1} \ifx \OpArg\Empty {} \else
        \label{#1}
  \fi}
\newtheorem{theoremSt}{Theorem}[section]

\newtheorem{exampleSt}[theoremSt]{Example}
\newtheorem{exerciseSt}[theoremSt]{Exercise}

\newcommand\MakeStEnv[1]{
  \newenvironment{#1}[1]{
  \begin{#1St} \oplabel{##1}%
  \global\def\CrntSt{\thetheoremSt}%
}{ 
  \end{#1St} }
  \newenvironment{#1+}[1]{
  \begin{#1St} \label{##1}%
  \label{DisplaY ##1}%
  \global\def\CrntSt{\thetheoremSt}%
  \def\Labl{##1}\ifx\Labl\Empty{} \else {\em (\Labl)\,}\fi%
}{ 
  \end{#1St} }
}
\MakeStEnv{theorem}
\MakeStEnv{corollary}
\MakeStEnv{proposition}
\MakeStEnv{lemma}
\MakeStEnv{definition}
\MakeStEnv{conjecture}
\MakeStEnv{problem}
\MakeStEnv{question}

\long\def\state#1#2{
\medskip\par\noindent
{\bf #1} 
{\it #2}
\par\medskip
}

\long\def\realfig#1#2{
\begin{figure}[htbp]
\centerline{\psfig{file=#1}}
\caption[#1]{#2}
\oplabel{#1}
\end{figure}}

\newlength{\saveu}

\newenvironment{pf}{%
 \begin{proof}%
}{ 
 \end{proof}
}
\newenvironment{pf*}[1]{%
 \begin{proof}[#1]%
}{ 
 \end{proof}
}

\newcommand{\finishproof}[1]{ 
  \def\FPArg{#1}
  \ifx\FPArg\Empty
        \newcommand\FPArg{\CrntSt}  \fi
  \smallbreak\noindent\makebox[\textwidth]{\hfill\fbox{\FPArg}}
  \medbreak\noindent
}

\newcommand\AAA{{\mathcal A}}

\newcommand\CC{{\mathcal C}}
\newcommand\DD{{\mathcal D}}

\newcommand\FF{{\mathcal F}}
\newcommand\GG{{\mathcal G}}

\newcommand\II{{\mathcal I}}

\newcommand\LL{{\mathcal L}}
\newcommand\MM{{\mathcal M}}

\newcommand\PP{{\mathcal P}}

\newcommand\RR{{\mathcal R}}

\newcommand\UU{{\mathcal U}}

\newcommand\XX{{\mathcal X}}
\newcommand\YY{{\mathcal Y}}

\newcommand\PMF{{\PP\kern-2pt\MM\FF}}
\newcommand\ML{{\MM\LL}}
\newcommand\PML{{\PP\kern-2pt\MM\LL}}
\newcommand\GL{{\GG\LL}}

\newcommand\half{{\textstyle{\frac12}}}

\newcommand\ep{\epsilon}

\newcommand\hhat{\widehat}

\newcommand\union{\cup}
\newcommand\intersect{\cap}
\newcommand\bbR{{\mathord{\text{I\kern-2pt R}}}}        
\newcommand\bbH{{\mathord{\text{I\kern-2pt H}}}}        

\newcommand\C{{\mathbb C}}

\newcommand\Z{{\mathbb Z}}
\newcommand\R{{\mathbb R}}

\newcommand\Hyp{{\mathbb H}}

\newcommand\PSL[1]{\text{PSL}_{#1}}

\newcommand\bigrightarrow[1]{\hbox to #1{\rightarrowfill}}
\newcommand\bigleftarrow[1]{\hbox to #1{\leftarrowfill}}

\newcommand\homeo{\cong}
\newcommand\boundary{\partial}
\newcommand\semidir{\mathrel{\hbox{\vrule depth-.03ex height1.1ex\kern-0.15em$\times$}}}

\newcommand\del{\nabla}
\newcommand\til{\widetilde}

\newcommand{\diam}{\operatorname{diam}}

\renewcommand{\Re}{\operatorname{Re}}
\renewcommand{\Im}{\operatorname{Im}}

\numberwithin{equation}{section}

\catcode`\@=11

\def\subsection{\@startsection{subsection}{2}%
  \z@{.5\linespacing\@plus.7\linespacing}{.5em}%
  {\normalfont\bfseries\centering}}

\def\section{\@startsection{section}{1}%
  \z@{.7\linespacing\@plus\linespacing}{.5\linespacing}%
  {\normalfont\large\bfseries\centering}}

\def\subsubsection{\@startsection{subsubsection}{3}%
  \z@{.5\linespacing\@plus.7\linespacing}{-.5em}%
  {\normalfont\bfseries}}

\catcode`\@=12

\newcommand{\T}{{\mathbf T}}

\newcommand{\pleat}{\operatorname{\mathbf{pleat}}}
\newcommand{\short}{\operatorname{\mathbf{short}}}

\newcommand{\collar}{\operatorname{\mathbf{collar}}}
\newcommand{\bcollar}{\operatorname{\overline{\mathbf{collar}}}}

\newcommand{\I}{{\mathbf I}}

\newcommand{\sprec}{\prec_s}

\newcommand{\fsub}{\mathrel{\scriptstyle\searrow}}
\newcommand{\bsub}{\mathrel{\scriptstyle\swarrow}}
\newcommand{\fsubd}{\mathrel{{\scriptstyle\searrow}\kern-1ex^d\kern0.5ex}}
\newcommand{\bsubd}{\mathrel{{\scriptstyle\swarrow}\kern-1.6ex^d\kern0.8ex}}
\newcommand{\fsubeq}{\mathrel{\raise-.7ex\hbox{$\overset{\searrow}{=}$}}}
\newcommand{\bsubeq}{\mathrel{\raise-.7ex\hbox{$\overset{\swarrow}{=}$}}}
\newcommand{\tw}{\operatorname{tw}}
\newcommand{\base}{\operatorname{base}}

\newcommand{\rest}{|_}
\newcommand{\bbar}{\overline}
\newcommand{\UML}{\operatorname{\UU\MM\LL}}
\newcommand{\EL}{\mathcal{EL}}

\newcommand{\tsh}[1]{\left\{\kern-.9ex\left\{#1\right\}\kern-.9ex\right\}}
\newcommand{\Tsh}[2]{\tsh{#2}_{#1}}
\newcommand{\qeq}{\mathrel{\approx}}
\newcommand{\Qeq}[1]{\mathrel{\approx_{#1}}}
\newcommand{\qle}{\lesssim}
\newcommand{\Qle}[1]{\mathrel{\lesssim_{#1}}}
\newcommand{\simp}{\operatorname{simp}}
\newcommand{\vsucc}{\operatorname{succ}}
\newcommand{\vpred}{\operatorname{pred}}

\newcommand\sbtop{_{\text{top}}}
\newcommand\sbot{_{\text{bot}}}
\newcommand\sll{_{\mathbf l}}
\newcommand\srr{_{\mathbf r}}

\newcommand\boundvert{{\boundary_{||}}}
\newcommand\storus[1]{U(#1)}
\newcommand\Homega{\omega_H}
\newcommand\Momega{\omega_M}
\newcommand\nomega{\omega_\nu}
\newcommand\twist{\operatorname{tw}}
\newcommand\modl{M_\nu}
\newcommand\MT{{\mathbb T}}
\newcommand\Teich{{\mathcal T}}

\newcommand\bersfcn{\LL}
\newcommand\neck{\varepsilon}
\newcommand\epshrink{\ep_{T}}
\newcommand\shear{\operatorname{\mathbf{shear}}}

\newcommand\sitem{\medskip\item}

\begin{document}

\title[The classification of Kleinian surface groups, I]{The
classification of Kleinian surface groups, I: Models and bounds} 

\author{Yair N. Minsky}
\address{Yale University}
\date{November 30, 2004}
\thanks{Partially supported by NSF grant DMS-9971596. The bulk of this
work was completed at SUNY at Stony Brook.}

\begin{abstract}
We give the first part of a proof of Thurston's Ending Lamination
conjecture. In this part we show how to construct from the end
invariants of a Kleinian surface group a ``Lipschitz model'' for the
thick part of the corresponding hyperbolic manifold. This enables us
to describe the topological structure of the thick part, and to give 
a-priori geometric bounds.
\end{abstract}

\maketitle

\renewcommand\marginpar[1]{} 

\newcommand\ME{M\kern-4pt E}
\newcommand\bME{\overline{M\kern-4pt E}}

\setcounter{tocdepth}{1}
\tableofcontents

\section{Introduction}
\label{intro}

This paper is the first in a two-part series addressing the question:
to what extent is a hyperbolic 3-manifold determined by its asymptotic
geometry? This question underlies the deformation theory of
Kleinian groups, as pioneered by Ahlfors and Bers in the 60's and by
Thurston and Bonahon in the 70's and early 80's.  Their work provides us
with a theory of {\em end invariants} assigned to the ends of a hyperbolic
3-manifold, and determined by their asymptotic geometric properties.
Thurston \cite{wpt:bull} formulated this conjecture which has
been a guiding question in the field:

\state{Ending Lamination Conjecture.}{
A hyperbolic 3-manifold with
finitely generated fundamental group is uniquely determined by its
topological type and its end invariants.}

When the manifold has finite volume its ends are either empty or
cusps, the end invariants are empty,
and the conjecture reduces to the well-known rigidity
theorems of Mostow and Prasad \cite{mostow:hyperbolic,prasad}. When
the manifold has infinite volume but is ``geometrically finite'', the
end invariants are Riemann surfaces arising from the action of
$\pi_1(N)$ on the Riemann sphere, and the conjecture follows from 
the work of Ahlfors-Bers \cite{ahlfors-bers,bers:simultaneous,bers:spaces}
and Marden-Maskit \cite{marden-maskit}, Maskit \cite{maskit:self},
Kra \cite{kra:spaces} and others.

The remaining cases are those where the manifold has a ``geometrically
infinite'' end, for which the end invariant is a {\em lamination.}
Here the discussion splits into two, depending on whether the boundary
of the compact core is compressible or incompressible. If it is
incompressible, the work of Thurston \cite{wpt:notes} and Bonahon
\cite{bonahon} gives a preliminary geometric description of the end,
and allows the ending laminations to be defined 
(see also Abikoff \cite{abikoff:survey} for a survey). If the core boundary
is compressible, the deeper question of {\em geometric tameness of the end} is
yet to be resolved, although Canary showed \cite{canary:ends} that  
this is equivalent to {\em topological tameness}, namely that the
manifold is the interior of a compact manifold. Marden conjectured in
\cite{marden:geometry} that this is always the case, and this remains
one of the central open questions in the field. 

We therefore restrict ourselves to the incompressible boundary case. 
This case reduces, by restriction to boundary
subgroups, to the case of {\em (marked) Kleinian surface
groups}, with which we will be concerned for the remainder of the paper.

A marked Kleinian surface group is a discrete, faithful representation
$$
\rho:\pi_1(S) \to \PSL 2(\C)
$$
where $S$ is a compact surface, and $\rho$ sends elements representing
$\boundary S$ to parabolic elements. Each $\rho$ determines
a set of end invariants $\nu(\rho)$, which for each end
give us Ahlfors-Bers Teichm\"uller data or an ending lamination, as
appropriate. 

In broadest outline, our plan for establishing the Ending Lamination
Conjecture is to construct a ``model manifold'' $\modl$, depending
only on the invariants $\nu(\rho)$, together with a bilipschitz
homeomorphism $f:\modl \to N_\rho$. Then if $\rho_1$ and $\rho_2$ are
two Kleinian surface groups with the same end invariants $\nu$, we
would obtain a bilipschitz homeomorphism between $N_{\rho_1}$ and
$N_{\rho_2}$ (in the right homotopy class), and an application of
Sullivan's rigidity theorem \cite{sullivan:rigidity} would then imply
that the map can be deformed to an isometry.

In this paper we will construct the model manifold together with a map
saisfying some {\em Lipschitz} bounds (and some additional geometric
properties, including detailed information about the thick-thin
decomposition of $N_\rho$). In the second paper, with Brock and Canary
\cite{brock-canary-minsky:ELCII}, this map will be promoted to a
bilipschitz homeomorphism.

\subsubsection*{Structure of the model}

For simplicity, let us describe the model
manifold $\modl$ when $S$ is 
a closed surface, and when $\nu$ are invariants of a manifold $N_\rho$
without parabolics, and without geometrically finite ends. 
(This avoids discussion of parabolic cusps and boundaries of the convex
core). In this case, $\modl$ is homeomorphic to $S\times\R$, and we
fix such an identification.

Within $\modl$ there is a subset $\UU$, which consists of open solid
tori called ``tubes'' of the form $U = A\times J$, where $A$ is an
annulus in $S$ and 
$J$ is an interval in $\R$. No two components of $\UU$ are homotopic.

$\modl$ comes equipped with a piecewise-Riemannian metric, with
respect to which each tube boundary $\boundary U$ is a Euclidean torus.
The geometry of $\boundary U$ is described by a
coefficient we call $\Momega(U)$, which lies in the upper half plane
$\Hyp^2 = \{z:\Im z>0\}$, thought of as the Teichm\"uller space of the
torus.  $U$
itself is isometric to a tubular neighborhood of a hyperbolic
geodesic, whose length goes to 0 as $|\Momega|\to \infty$.

Let $\UU[k]$ denote the union of components of $\UU$ with $|\Momega|\ge
k$, and let $\modl[k] = \modl \setminus \UU[k]$.
Then $\modl[0] = \modl\setminus \UU$ is a union of ``blocks'',
which have a finite number of possible isometry types.
This describes a 
sort of ``thick-thin'' decomposition of $\modl$. 

There is a
corresponding decomposition of $N_\rho$, associating a Margulis tube
$\MT_{\ep_1}(\gamma)$ to each sufficiently short geodesic $\gamma$. 
Let $\MT[k]$ denote the set of such Margulis tubes (if any) associated 
to the homotopy classes of components of $\UU[k]$ under the
homotopy equivalence between $M_\nu$ and $N_\rho$ determined by $\rho$.

Let $\hhat C_{N_\rho}$ denote the ``augmented convex core'' of
$N_\rho$ (see \S\ref{augmented core}), which in our simplified case is
equal to $N_\rho$ itself.
Our main theorem asserts that $M_\nu$ can be mapped to $\hhat
C_{N_\rho}$ by a Lipschitz map that respects the thick-thin decompositions of
both.

\state{Lipschitz Model Theorem.}{
Fix a compact oriented surface $S$. There exist $K,k>0$ such that, if
$\rho:\pi_1(S)\to \PSL 2(\C)$ is a Kleinian surface group with end
invariants $\nu(\rho)$,
then there is a map
$$
f: \modl \to \hhat C_{N_\rho}
$$
with the following properties:
\begin{enumerate}
\item $f$ induces $\rho$ on $\pi_1$, is proper, and has degree 1.
\item $f$ is $K$-Lipschitz on $\modl[k]$, with respect to the induced
  path metric.
\item $f$ maps $\UU[k]$ to $\MT[k]$, and $\modl[k]$ to
$N_\rho\setminus \MT[k]$.
\item $f:\boundary \modl \to \boundary \hhat C_{N_\rho}$ is a
$K$-bilipschitz homeomorphism of the boundaries. 
\item For each tube $U$ in $\UU$ with $|\omega_M(U)|<\infty$, 
$f|_U$ is   $\lambda$-Lipschitz, where $\lambda$ depends only on
$K$ and $|\omega_M(U)|$.
\end{enumerate}
}

Remarks: 
The condition on the degree of $f$, after appropriate orientation
conventions, amounts to the fact that $f$ maps the ends of $\modl$ to
the ends of $\modl$ in the ``correct order''.

\subsubsection*{The extended model map}
In the general case, $N_\rho$ may have parabolic cusps and $\hhat
C_{N_\rho}$ 
may not be all of $N_\rho$. The statement of the Lipschitz Model
Theorem is unchanged, but
the structure of $\modl $ is 
complicated in several ways: Some of the tubes of $\UU$ will be
``parabolic'', meaning that their boundaries are annuli rather than
tori, and the coefficients $\Momega$ may take on the special value
$i\infty$.  $\modl$ will have a boundary, and the condition that $f$
is proper is meant to include both senses: it is proper as a map of topological
spaces, and it takes $\boundary \modl$ to $\boundary \hhat
C_{N_\rho}$. The blocks of $\modl[0]$ will still have a finite number
of topological types, but for 
a finite number of blocks adjacent to the
boundary the isometry types will be unbounded, in a controlled way.

In \S\ref{augmented core} we will describe the geometry of the exterior of the
augmented core, $E_N = \overline{N\setminus \hhat C_N}$, in terms of a 
model $E_\nu$ that depends only on the end invariants of the
geometrically finite ends of $N$, and is a variation of
Epstein-Marden's description of the exterior of the convex hull. 
$\modl$ and $E_\nu$ attach along their boundaries to yield an extended
model manifold $\ME_\nu$. 

$N$ has a natural conformal boundary at infinity $\boundary_\infty N$, 
and $\ME_\nu$ has a conformally equivalent boundary $\boundary_\infty\ME_\nu$.
The Lipschitz Model Theorem will then generalize to:

\state{Extended Model Theorem.}{
The map $f$ obtained in the Lipschitz Model Theorem extends to 
a proper degree 1 map
$$
f' : \ME_\nu \to N
$$
which restricts to a $K$-bilipschitz homeomorphism from $E_\nu$ to
$E_N$, and extends to a conformal map from $\boundary_\infty \ME_\nu$ 
to $\boundary_\infty N$.
}

\subsubsection*{Length bounds}

Note that Part (3) of the Lipschitz Model Theorem implies that for
every component of $\UU[k]$ there is in fact a corresponding Margulis
tube in $\MT[k]$, to which it maps properly. On the other hand the bounded
isometry types of blocks (ignoring the
boundary case)  and the Lipschitz bound on $f$ will imply that 
there is a lower bound $\ep>0$ on the injectivity radius of $N_\rho$
outside the image of $\UU$. In other words, the structure of $\modl$
determines the pattern of short geodesics and their Margulis tubes in
$N_\rho$. 

The following theorem makes this connection more precise. If $\gamma$
is a homotopy class of curves in $S$ let $\lambda_\rho(\gamma)$ denote
the complex translation length of the corresponding conjugacy class
$\rho(\gamma)$ in
$\rho(\pi_1(S))$. Its real part $\ell_\rho(\gamma)$, which we may
assume positive if $\rho(\gamma)$ is not parabolic, is
the length of the geodesic representative of this homotopy class in
$N_\rho$. If $\gamma$ is homotopic to the core of some tube $U$
in $\UU$, we define $\Momega(\gamma) \equiv \Momega(U)$.

\state{Short Curve Theorem.}{
There exist $\bar \ep>0$ and $c>0$ depending only on $S$, 
and for each $\ep>0$ there exists $K>0$, such that the following holds:
Let $\rho:\pi_1(S)\to \PSL 2(\C)$ be a Kleinian surface group and
$\gamma$ a simple closed curve in $S$.
\begin{enumerate}
\item If $\ell_\rho(\gamma) < \bar\ep$ then $\gamma$ is homotopic to a
core of some component $U$ in $\UU$.
\item (Upper length bounds)
If $\gamma$ is homotopic to the core of a tube in $\UU$ then
$$
|\Momega(\gamma)| \ge K \implies \ell_\rho(\gamma) \le \ep.
$$
\item (Lower length bounds)
If $\gamma$ is homotopic to the core of a tube in $\UU$ then
$$
|\lambda_\rho(\gamma)| \ge \frac c {|\Momega(\gamma)|}
$$
and 
$$
\ell_\rho(\gamma) \ge \frac c {|\Momega(\gamma)|^2}.
$$
\end{enumerate}
}

Part (2) of this theorem is actually a restatement of the main theorem
of \cite{minsky:kgcc}; part (3) is the main new ingredient.

\subsection{Outline of the proofs}

In the following summary of the argument, we will continue making
the assumptions that $S$ is closed and $N_\rho$ has no
geometrically finite ends or cusps. This greatly simplifies the logic
of the discussion, while retaining all the essential elements of the
proof. The reader is encouraged to continue making this assumption
on a first reading of the proof itself.

In this case, $N_\rho$ has two ends, which we label with $+$ and $-$
(see \S\ref{cores and ends} for the orientation conventions), and
the end invariants $\nu(\rho)$ become two filling laminations $\nu_+$ and
$\nu_-$ on $S$ (\S\ref{end invariants}).

\subsubsection*{Quasiconvexity and the complex of curves}
The central idea is to use the geometry of the 
{\em complex of curves} $\CC(S)$ to obtain a priori bounds on lengths
of curves in $N_\rho$. The vertices of $\CC(S)$ are the essential
homotopy classes of simple loops in $S$ (see \S\ref{complexes} for
details), and we will study the sublevel sets
$$
\CC(\rho,L) = \{v\in \CC_0(S): \ell_\rho(v) \le L\}
$$
where $\ell_\rho(v)$ for a vertex $v\in \CC_0(S)$ denotes the
length of the corresponding closed geodesic in $N_\rho$.

In \cite{minsky:boundgeom} we showed that $\CC(\rho,L)$ is {\em quasiconvex}
in the natural metric on $\CC(S)$. The main tool for the proof of this
is the ``short curve projection'' $\Pi_{\rho,L}$, which maps $\CC(S)$
to $\CC(\rho,L)$ by constructing for any 
vertex $v$ in $\CC(S)$ the set of pleated surfaces in $N$ with $v$ in their
pleating locus, and finding the curves of length at most $L$ in these
surfaces. This map satisfies certain contraction properties which make
it coarsely like a projection to a convex set, 
and this yields the quasiconvexity of $\CC(\rho,L)$.

As a metric space $\CC(S)$ is
$\delta$-hyperbolic, and the ending laminations $\nu_\pm(\rho)$ 
describe two points on its Gromov boundary $\boundary\CC(S)$
(see Masur-Minsky \cite{masur-minsky:complex1}, Klarreich
\cite{klarreich:boundary} and Section \ref{complexes}).
In fact
$\nu_\pm(\rho)$ are the accumulation points of $\CC(\rho,L)$ on
$\boundary\CC(S)$, and this together with quasiconvexity of
$\CC(\rho,L)$ appears to
give a coarse type of control on $\CC(\rho,L)$ -- in particular an
infinite geodesic in $\CC(S)$ joining $\nu_-$ to $\nu_+$ must lie in a
bounded neighborhood of $\CC(\rho,L)$.  However, since $\CC(S)$ is
locally infinite this estimate is not sufficient for us.

\subsubsection*{Subsurfaces and hierarchies}
In Section \ref{coarse projection} we generalize the quasiconvexity
theorem of \cite{minsky:boundgeom} to a {\em relative} result which 
incorporates the structure of subsurface complexes in $\CC(S)$. In
order to do this we recall in Sections \ref{complexes} and
\ref{hierarchies} some of the 
structure of the subsurface projections and hierarchies in $\CC(S)$
which were developed in Masur-Minsky \cite{masur-minsky:complex2}.
To  an essential subsurface $W\subset S$ we associate a 
``projection'' 
$$
\pi_W : \AAA(S) \to \AAA(W)
$$
(where $\AAA(W)$ is the arc complex of $W$, containing and
quasi-isometric to $\CC(W)$).
This, roughly speaking, is a map that associates to a curve (or arc) system
in $S$ its essential intersection with $W$. This map has properties
analogous to the orthogonal projection of $\Hyp^3$ to a horoball,
see particularly Lemma \ref{Geodesic Projection Bound}.

A {\em hierarchy} is a way of 
enlarging a geodesic in $\CC(S)$ to a system of geodesics in
subsurface complexes $\CC(W)$ that together produces families of
markings of $S$. Such a hierarchy, called $H_\nu$, is constructed
in \S\ref{hierarchies} so that its base geodesic $g$ connects $\nu_-$ to
$\nu_+$ (the construction is nearly the same as in
Masur-Minsky \cite{masur-minsky:complex2}, except for the need to treat infinite
geodesics). The vertices which appear in $H_\nu$ are all within
distance 1 of $g$ in $\CC(S)$. The structure of $H_\nu$ is strongly
controlled by the maps $\pi_W$, as in Lemma \ref{Large Link}.

\subsubsection*{Projections and length bounds}
Once this structure is in place, we revisit the map $\Pi_{\rho,L}$. We
prove, in Theorem \lref{Relative Coarse Projection}, that the
composition $\pi_Y \circ \Pi_{\rho,L}$ for a subsurface $Y$ has
contraction properties generalizing those shown in \cite{minsky:boundgeom}.
We then prove Theorem \ref{Projection
Bounds}, which states in particular that
$$
d_Y(v,\Pi_{\rho,L}(v))
$$
is uniformly bounded for any subsurface $Y$ and all vertices $v$
appearing in $H_\nu$, provided $v$ intersects $Y$ essentially. Here
$d_Y(x,y)$ denotes distance in $\AAA(Y)$ between $\pi_Y(x)$ and $\pi_Y(y)$.

This bound implies that $v$ and the
bounded-length curves $\Pi_{\rho,L}(v)$ are not too different
in some appropriate combinatorial sense, and indeed we go on to apply
this to obtain, in Lemma \ref{Upper Bounds}, an a-priori upper bound on
$\ell_\rho(v)$ for all vertices $v$ that appear in $H_\nu$. 
Another crucial result we prove along the way is Lemma \lref{Tube
Penetration}, which limits the ways in which pleated surfaces
constructed from vertices of the hierarchy can penetrate Margulis
tubes.

\subsubsection*{Model manifold construction}
At this point we are ready to build the model manifold. In Section
\ref{model} we construct $\modl$ out of the combinatorial data in
$H_\nu$. The blocks of $\modl[0]$ 
are constructed from edges of geodesics in $H_\nu$ associated to
one-holed torus and 4-holed sphere subsurfaces, and glued together
using the ``subordinacy'' relations in $H_\nu$. The structure of
$H_\nu$ is also used to embed $\modl[0]$ in $S\times\R$ (after which
we identify it with its embedded image), and
the tubes $\UU$ are the solid-torus components of $S\times\R
\setminus\modl[0]$,  
and are in one-to-one correspondence with the vertices of $H_\nu$. In
\S\ref{define metric} we 
introduce the meridian coefficients $\Momega(v)$, which encode for
each vertex $v$ the geometry of the associated tube boundary. The
metric of $\modl$ is described in this section too. 

In Section \ref{count} we define alternative meridian coefficients
$\Homega$ and $\nomega$, which are computed, respectively, directly
from the data of $H_\nu$ and directly from $\nu$ itself. It is useful
later in the proof to compare all three of these and in Theorem
\ref{Omegas close} we show that they are essentially equivalent. The
proof requires a somewhat careful analysis of the geometry of the
model, and some counting arguments using the structure of the hierarchy.

\subsubsection*{Lipschitz bounds}
In Section \ref{lipschitz} we finally build the Lipschitz map from
$\modl$ to $N_\rho$, establishing the Lipschitz Model Theorem. This is
done in several steps, starting with the ``gluing boundaries'' of
blocks, where the a-priori bound on vertex lengths from Lemma
\ref{Upper Bounds} provides the Lipschitz control. Extension to the
``middle surfaces'' of blocks  (Step 2) requires another application of
Thurston's Uniform Injectivity Theorem, via Lemma \ref{halfway
bounds}. Control of the extension to the rest of the blocks requires a
reprise of the ``figure-8 argument'' from \cite{minsky:torus} to bound
homotopies between Lipschitz maps of surfaces (Step 4). The map can be
extended to tubes, and the last subtle point comes in Step 7, where we
need Lemma \ref{Big omega short curve} to relate large meridian
coefficients to short curves (this is the point where we apply the
results of Section \ref{count}, as well as the main theorem of
\cite{minsky:kgcc}).

The proof of the Short Curve Theorem, carried out in Section
\ref{margulis}, is now a simple consequence of the Lipschitz Model
Theorem together with the properties of Margulis tubes.  Roughly
speaking, an upper bound on $|\Momega(\gamma)|$ gives an upper bound
on the meridian disk of $U(\gamma)$, and hence a lower bound on the
length of its geodesic core. 

\subsubsection*{Preliminaries}
Sections \ref{end invariants} through \ref{hierarchies} provide 
some background and notation before the proof itself
starts in Section \ref{coarse projection}. Section \ref{end
invariants} introduces compact cores, ends and laminations. Section
\ref{hypstuff} introduces pleated surfaces, Margulis tubes and collars
in surfaces, and the augmented convex core. There is only a little bit
of new material here:  the augmented convex core and particularly
the geometric structure of its boundary, via Lemma \ref{infinity to
aug core}, and a slightly technical variation (Lemma \ref{collar properties})
on the standard collar 
of a short geodesic or cusp in a hyperbolic surface.
Sections \ref{complexes} and \ref{hierarchies} introduce
the complexes of curves and arcs, and hierarchies. Most of this is
review of material from \cite{masur-minsky:complex2}, with 
certain generalizations to the case involving
infinite geodesics. In particular the existence of an infinite hierarchy
connecting the invariants $\nu_+$ and $\nu_-$ (and in fact any pair
of generalized markings) 
is shown in Lemma \ref{Existence of Hierarchies}, and the existence of
a ``resolution'' of a hierarchy, which is something like a sequence of
markings separated by elementary moves sweeping through all the data
in the hierarchy, is shown in Lemmas \ref{Resolution exists},
\ref{Resolution sweep}.

An informal but extensive summary of the argument, focusing on the
case where $S$ is a five-holed sphere, can be found in the lecture
notes \cite{minsky:warwick}.


\subsubsection*{Bilipschitz control of the model map}
In \cite{brock-canary-minsky:ELCII} we will show that the map $f:\modl
\to \hhat C_N$ can be made a bilipschitz homeomorphism, and thereby
establish the Ending Lamination Conjecture. The argument begins with
a decomposition of $\modl$ along quasi-horizontal slices into pieces
of bounded size. The surface embedding
machinery of Anderson-Canary \cite{anderson-canary:cores} is used to
show that the slices in the 
boundaries of the pieces can be deformed to embeddings in a uniform
way. The resulting map preserves 
a certain topological partial order among these slices -- this is
established using an argument by contradiction and passage to
geometric limits -- and this can be used to make the map an
orientation-preserving embedding
on each piece separately. Uniform bilipschitz bounds on these embeddings
are obtained again by contradiction and geometric limit, and a global
bilipschitz bound follows.

We note that the structure of the model manifold provides new
information even in the geometrically finite case, for which the
Ending Lamination Conjecture itself reduces to the quasiconformal deformation
theory of Ahlfors and Bers. In particular it describes the thick-thin
decomposition of the manifold, its volume and other geometric
properties explicitly in terms of the Riemann surfaces at infinity.

\subsubsection*{Acknowledgements} The author is grateful to Howard
Masur, Dick Canary and Jeff Brock, without whose collaborations and
encouragement this project would never have been completed.

\section{End invariants}
\label{end invariants}

Before we discuss the end invariants of a hyperbolic 3-manifold, we
set some notation which will be used throughout this paper.
Let $S_{g,b}$ denote an oriented, connected and compact surface of genus
$g$ with $b$ boundary components. 
Given a surface $R$ of this type, 
an {\em essential subsurface} $Y\subseteq R$ is a compact, connected
subsurface all of whose  
boundary components are homotopically nontrivial, and so that $Y$
is not homotopic into a boundary component of $R$.
(Unless otherwise mentioned we assume throughout the paper that any
subsurface is essential.)

Define the {\em complexity} of a surface to be
$$\xi(S_{g,b}) \equiv 3g+b.$$
Note that for an essential subsurface $Y$ of $R$, 
$\xi(Y)<\xi(R)$ unless $Y$ is homeomorphic to $R$.

A {\em hyperbolic structure on $int(R)$} will be a hyperbolic metric
whose completion is a hyperbolic surface with cusps and/or geodesic
boundary components, the latter of which we may identify with
components of $\boundary R$. 
If the completion is compact then its boundary is identified with all
of $\boundary R$, and we call this a
{\em hyperbolic structure on $R$}.

The Teichm\"uller space $\Teich(R)$ is the space of (marked)
hyperbolic structures on $int(R)$ which are complete -- that is, all
ends are cusps. Alternatively $\Teich(R)$ is the space of marked
conformal structures for which the ends are punctures.

\subsection{Geodesic laminations}
\label{laminations}

We will assume the reader is familiar with the basics of geodesic
laminations on hyperbolic surfaces --
see Casson-Bleiler \cite{casson-bleiler} or the recently written
Bonahon \cite{bonahon:laminations} for an
introduction to this subject. Given a surface $R$ with a complete hyperbolic
structure on $int(R)$ (in the above sense, with all ends cusps), we
will be using the following spaces:

$\GL(R)$ is the space of geodesic laminations on $R$. 
$\ML(R)$ is the space of {\em tranversely measured} laminations on $R$
with compact support (if $\lambda\in\ML(R)$ its {\em support},
$|\lambda|$, lies in 
$\GL(R)$). $\GL(R)$ is usually topologized using the Hausdorff topology
on closed subsets of $R$, whereas
$\ML(R)$ admits a topology (due to Thurston) coming from the
weak-* topology of the measures induced on transverse arcs.

$\UML(R)$ is the   quotient space of $\ML(R)$ obtained by forgetting
the measures (the ``unmeasured laminations'').
Its topology is non-Hausdorff and it is rarely used as an object on
its own.
Note that $\UML(R)$ is set-theoretically contained in $\GL(R)$, but
the topologies are different.

It is part of the basic structure theory of laminations that 
every element of $\UML(R)$ decomposes into a finite union of disjoint
connected components, each of which is a 
minimal lamination.

Let $\EL(R)$ denote the image in $\UML(R)$ of the {\em filling
  laminations} in $\ML(R)$.
  $\mu\in\ML(R)$ is called {\em filling} if it has transverse intersection
  with any $\mu'\in\ML(S)$, unless $\mu$
  and $\mu'$ have the same support. 
  An equivalent condition is that $\mu$ is both
  minimal and maximal as a measured lamination, and another is that
  the complementary components of $\mu$ are ideal polygons or
  once-punctured ideal polygons.
(see \cite{wpt:notes}).
   In the topology inherited from $\UML(S)$, 
  $\EL(S)$ is Hausdorff (see \cite{klarreich:boundary} for a proof).

{\em Remarks:} 1. All of these spaces do not really depend on
the hyperbolic structure of $int(R)$ -- that is, the spaces obtained from
two different choices of structure are canonically homeomorphic.
2. The laminations in $\EL(R)$ are exactly those that appear as
ending laminations of Kleinian surface groups $\rho\in\DD(R)$ without
accidental parabolics -- this is the reason for the notation $\EL$.

\subsection{Cores and ends}
\label{cores and ends}

Let $N$ be an oriented complete hyperbolic 3-manifold with finitely generated
fundamental group. $N$ may be expressed as the quotient
$\Hyp^3/\Gamma$ by a Kleinian group 
$\Gamma \homeo \pi_1(N)$. Let $\Lambda$ be the limit set of $\Gamma$
in the Riemann sphere $\hhat \C$, 
and let $\Omega = \hhat \C \setminus \Lambda$ be the
domain of discontinuity  of $\Gamma$. (See
e.g. \cite{maskit:book,ceg,benedetti-petronio} for background)
Then
$$
\overline N \equiv (\Hyp^3 \union \Omega)/\Gamma
$$
is a 3-manifold with boundary $\boundary \overline N = \Omega/\Gamma$
and interior $N$. 
The boundary inherits a conformal structure from $\Omega$. 

Let $C_N$ denote the {\em convex core} of $N$, which is the
quotient by $\Gamma$ of the convex hull of $\Lambda$ in $\Hyp^3$. 
$C_N$ is homeomorphic to $\overline N$ by a map homotopic to the inclusion
(except when $\Gamma$ is Fuchsian or elementary and $C_N$ has
dimension 2 or less; we will assume this is not the case). Thurston
showed that $\boundary C_N$ is a union of  
hyperbolic surfaces in the induced path metric from $N$, and by a
theorem of Sullivan, if $\boundary C_N$ is incompressible 
these metrics are within universal bilipschitz
distortion from the Poincar\'e metric on $\boundary \overline N$
(see Epstein-Marden \cite{epstein-marden}). Ahlfors' finiteness
theorem \cite{ahlfors:finitegen} states that $\boundary \overline N$ and
$\boundary C_N$ have finite hyperbolic area. 

Let $Q$ denote
the union of (open) $\ep_0$-Margulis tubes of cusps of $N$ (see
\S\ref{margulis tubes} for a discussion of Margulis tubes), and let $N_0 =
N\setminus Q$. 
By Scott's compact core theorem \cite{scott:core} there is a compact
3-manifold $K\subset N$ whose inclusion is a homotopy equivalence.
The relative core theorem of McCullough \cite{mccullough:relative}
and Kulkarni-Shalen \cite{kulkarni-shalen} tells us that $K$ can be
chosen in $N_0$ so that 
$\boundary K$ meets the boundary of each rank-1 cusp of $Q$ in
an essential annulus, and contains the entire torus boundary of each
rank-2 cusp. 
Let $P=\boundary K \intersect \boundary Q$. Note that no two
components of $P$ can be homotopic, since
no two Margulis tubes can have homotopic core curves.
The topological ends of
$N_0$ are in one to one correspondence with the components of
$N_0\setminus K$ (which are {\em neighborhoods} of the ends), and
hence with the components of $\boundary K \setminus P$ (see Bonahon
\cite{bonahon}). If $R$ is the closure of a component 
of $\boundary K\setminus P$, let $E_R$ be the component of
$N_0\setminus K$ adjacent to $R$. We say that the end {\em faces $R$}.

\subsubsection*{Geometrically finite ends}
An end of $N_0$ is {\em geometrically finite} if its associated
neighborhood $E_R$ meets $C_N$ in a bounded subset. This implies
that the boundary of $E_R$ in $\overline N$ consists of $R$, 
some annuli in $\boundary Q$, and a component $X$ of $\boundary
\overline N$ which is homotopic to $R$.
Indeed, $K$ may be chosen so that $E_R \homeo int(R)\times (0,1)$.
The conformal structure of $X$ gives rise to a point in $\Teich(R)$,
and this is the end invariant associated to $R$, which we name $\nu_R$.

\subsubsection*{Geometrically infinite ends}
An end of $N_0$ is {\em geometrically infinite} if its associated
$E_R$ intersects $C_N$ in an unbounded set. Ahlfors' finiteness
theorem \cite{ahlfors:finitegen}
implies that $\boundary C_N \intersect N_0$ is compact, and hence
cannot separate $E_R$ into two unbounded sets. Thus
it follows that in fact that there is a (possibly smaller) neighborhood of the
end which is contained in $C_N$. 

In order to describe the end invariant for a geometrically infinite
end, we must consider the following definition from Thurston
\cite{wpt:notes}:

\begin{definition}{def degenerate}
An end $E_R$ of $N_0$ is {\bf simply degenerate} if there exists a sequence of
essential simple closed curves $\alpha_i$ in $R$ whose geodesic
representatives $\alpha_i^*$ exit the end.
\end{definition}

Here ``exiting the end'' means that the geodesics are eventually
contained in $E_R$ minus 
any bounded subset. Note that a geometrically finite
end cannot be simply degenerate, since all closed geodesics are
contained in the convex hull. 

Thurston established this theorem  (see also Canary \cite{canary:ends}):

\begin{theorem}{EL exist}{\rm (Thurston \cite{wpt:notes})}
Let $e$ be an end of $N_0$ facing $R \subset\boundary K$,
and suppose that $R$ is incompressible in $K$. 
If $e$ is simply degenerate then there exists a unique
lamination $\nu_R\in\UML(R)$ such that for any sequence of simple closed
curves $\alpha_i$ in $R$, 
$$
\alpha_i \to \nu_e \iff \text{$\alpha_i^*$ exit the end $e$.}
$$
A sequence $\alpha_i\to \nu_e$ can be chosen so that the lengths
$\ell_N(\alpha_i^*) \le L_0$, where $L_0$ depends only on $S$.

Furthermore, $\nu_e\in\EL(R)$ -- that is, it {\em fills} $R$.
\end{theorem}

Thurston also proved that a simply degenerate end is
{\em topologically tame}, meaning that it has a neighborhood homeomorphic to
$R\times(0,\infty)$, 
and that manifolds obtained as limits of quasifuchsian manifolds
have ends that are geometrically finite or simply degenerate. Bonahon
completed the picture, in the incompressible boundary case,  with his
``tameness theorem'',  

\begin{theorem}{Bonahon tameness}{\rm \cite{bonahon}}
Suppose that each component of $\boundary K \setminus P$ is
incompressible in $K$. Then 
the ends of $N_0$  are either geometrically finite or simply degenerate.
\end{theorem}
An end that is geometrically finite or simply degenerate is known as
{\em geometrically tame}.

The case of an end facing a compressible boundary component is
considerably harder to understand. Canary \cite{canary:ends} showed
that the analogue of Bonahon's theorem holds for such ends (with a
suitable strengthened notion of simple degeneracy)
if the end is known to be topologically tame. 
Marden had conjectured in \cite{marden:geometry} that all hyperbolic
3-manifolds with finitely generated fundamental groups have
topologically tame ends, and this conjecture remains open in the
compressible-boundary case (see \cite{brock-bromberg-evans-souto:ahlfors}
for recent progress on this question).

\subsubsection*{Ends for Kleinian surface groups}
From now on restrict to the case of a Kleinian surface group
$\rho:\pi_1(S) \to \PSL 2(\C)$, and let $N=N_\rho$ be the quotient
manifold $\Hyp^3/\rho(\pi_1(S))$.  Let $Q$ and $K$ be defined as above.

Let  $Q_0\subseteq Q$ be the set of cusp tubes associated to the
boundary of $S$ (if any), and let $P_0$ be the union of annuli in
$\boundary Q_0\intersect K$. Then $K$ is homeomorphic to 
$S\times[-1,1]$ with $P_0$ identified with $\boundary S \times
[-1,1]$.

Divide the remaining components of $Q$ into $Q_+$ and $Q_-$ according
to whether they meet $K$ on
$S\times\{1\}$ or $S\times\{-1\}$.  Let $P$ denote the union of annuli
$\boundary Q\intersect K$, divided similarly into $P_0, P_+$ and $P_-$.

The closure $R$ of a component of $\boundary K \setminus P$ then has
an associated end invariant $\nu_R$ as above, in either $\Teich(R)$ or
$\EL(R)$. We group these according to whether they come from the
``$+$'' or ``$-$'' side. More specifically, let
$p_+$ denote the set of core curves of $P_+$. 
Let $R^L_+$ denote the union of (closures of) components of
$S\times\{1\}\setminus P_+$, for which the invariant $\nu_R$ is a
lamination in $\EL(R)$. Let $R^T_+$ denote the remaining components. 
Define $\nu_+$ as a pair $(\nu^L_+,\nu^T_+)$, where
$\nu^L_+\in\UML(S)$ is the union of $p_+$ and the laminations 
of components of $R^L_+$, and $\nu^T_+$ is the set of Teichm\"uller
end invariants of components of $R^T_+$, which can be seen as an element of
$\Teich(R^T_+)$. 
We also have to allow either $\nu^T_+$ or $\nu^L_+$
to be empty, in the case that there are no lamination or Teichm\"uller
invariants, respectively. Define $\nu_-$ in the analogous way.

We remark that, because distinct cusps in a hyperbolic manifold cannot
have homotopic curves, the parabolics $p_+$ and $p_-$ have no elements
in common. This implies, when $p_+$ or $p_-$ are nonempty, that 
$\nu^L_+$ and $\nu^L_-$ have no infinite-leaf components in common
either. In fact, this is true if $p_+=p_-=\emptyset$ as well, as shown by
Thurston \cite{wpt:notes}.

\realfig{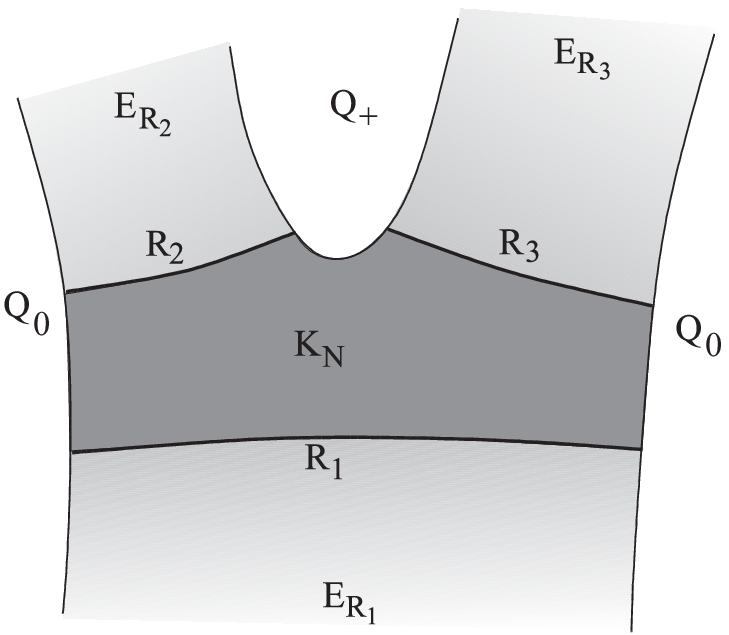}{A schematic of the division of $N$
into cusps, compact core and ends.}

It may be helpful to discuss special cases for clarity: 
If $\rho$ is
quasifuchsian, then $p_\pm$ are both empty (there are no parabolics
aside from those associated to $\boundary S$) and $R^T_+$ and $R^T_-$
are both copies of $S$. Thus $\nu_\pm = (\emptyset,\nu^T_\pm)$, where
$\nu^T_\pm\in\Teich(S)$ are the classical Ahlfors-Bers parameters for
$\rho$. In this case of course the Ending Lamination Conjecture is
well established, but our construction will still provide new
information about the geometry of $N_\rho$.

If $\rho $ has no parabolics aside from $\boundary S$ and no
geometrically finite ends, then $\nu_\pm = (\nu^L_\pm,\emptyset)$
with $\nu^L_\pm \in \EL(S)$. This is called the {\em doubly degenerate
case}.

In Section \ref{define Hnu}, we will replace $\nu_\pm$ with a pair of {\em
generalized markings} which will be used in the rest of the construction.

\subsubsection*{Orientation convention}
In order to keep track of which end of $N_\rho$ is which, let us
fix an orientation convention throughout the paper.
An orientation on a manifold $X$ will induce an orientation on
$\boundary X$ by the convention that, if $e$ is a baseframe for
$T_p(\boundary X)$ and $e'$ is an inward-pointing vector in $T_pX$,
then $e$ is positively oriented if and only if $(e,e')$ is.
We orient hyperbolic space $\Hyp^3$ so that it induces the standard
orientation on its boundary $\hat \C$, and note that the induced 
orientation on the upper half plane 
$\Hyp^2$ in turn induces the standard orientation on its boundary
$\hat\R$. 

The manifold $N_\rho$ inherits an orientation from $\Hyp^3$, and the
compact core $K$ inherits one as well. 
We are given a fixed orientation on $S$, and this determines
(up to proper isotopy) an identification of $K$ with $S\times [-1,1]$
by the condition that the induced boundary orientation on
$S\times\{-1\}$ agrees with the given orientation on $S$. Thus
we know which end is up.

\section{Hyperbolic constructions}
\label{hypstuff}

\subsection{Pleated surfaces}
\label{pleated}
\label{pleated facts}
A {\em pleated surface} is a map $f:int(S)\to N$ together with a hyperbolic
 structure on $int(S)$, written $\sigma_f$ and called the {\em induced
 metric}, and a $\sigma_f$-geodesic lamination $\lambda$ on $S$,
 so that the following holds: 
 $f$ is length-preserving on paths, maps leaves of $\lambda$ to geodesics, and
 is totally geodesic on the complement  of $\lambda$. 
Pleated surfaces were
 introduced by Thurston \cite{wpt:notes}. See Canary-Epstein-Green
 \cite{ceg} for more details. 
We include the case that 
the hyperbolic structure on $int(S)$ is incomplete and the 
 boundary of the completed surface is mapped geodesically. 
In almost every case in this paper, however, 
 boundary components of $S$ are mapped to cusps, so that
 $\sigma_f$ will be a complete metric with all ends cusps, 
 and leaves of $\lambda$ will go straight out this cusp in
 the $\sigma_f$ metric.
 
As in \cite{minsky:kgcc} we extend the definition slightly to include
{\em noded pleated surfaces}:
Suppose  $S'$ is an essential subsurface
of $S$  whose complement is a disjoint union of 
open collar neighborhoods of simple curves $\Delta$.
Let $[f]$ be a homotopy class of  maps from $S$ to $N$, such that
$f$ takes $\Delta$ to cusps of $N$.
We say that $g:S'\to N$ is a {\em noded
  pleated surface} in the class $[f]$ if $g$ is pleated with respect
to a hyperbolic metric on $S'$ (in which the ends are cusps), and $g$ is
homotopic to the restriction to $S'$ of an element of $[f]$. We say
that $g$ is ``noded on $\Delta$''. By convention, we represent the ``pleating
locus'' $\lambda$ of $f$ as a lamination in $\GL(S)$ that contains $\Delta$ as
components, and 
leaves of $\lambda$ that spiral onto $\Delta$ will be taken to leaves
that go out the corresponding cusp.

Now fixing a Kleinian surface group $\rho:\pi_1(S)\to \PSL 2(\C)$, 
we have a natural homotopy class of maps $S\to N_\rho$ inducing $\rho$
on fundamental groups. If $\lambda\in\GL(S)$, define
$$\pleat_\rho(\lambda)$$ 
to be the set of pleated surfaces $g:int(S)\to N_\rho$ in
the homotopy class of $\rho$, which map the leaves of $\lambda$
to geodesics (if some closed leaves of $\lambda$ correspond to 
parabolic elements in $\pi_1(N_\rho)$ then we allow the surface to be
noded on these curves). 

When $v$ is a pants decomposition there is only a finite number of
laminations  containing the curves of $v$, so $\pleat_\rho(v)$ is 
finite up to isotopy equivalence (the natural equivalence of precomposition by
homeomorphisms homotopic to the identity). Each of these laminations
consists of the curves of $v$, together with finitely many arcs that
spiral around them. We call $g\in\pleat_\rho(v)$ {\em leftward} if
this spiraling is to the left for each leaf. This determines a unique
isotopy equivalence class. Figure \ref{spin} shows a 
leftward-spinning lamination, restricted to a pair of pants, and also
indicates how it can be obtained from a triangulation by ``spinning''.

\realfig{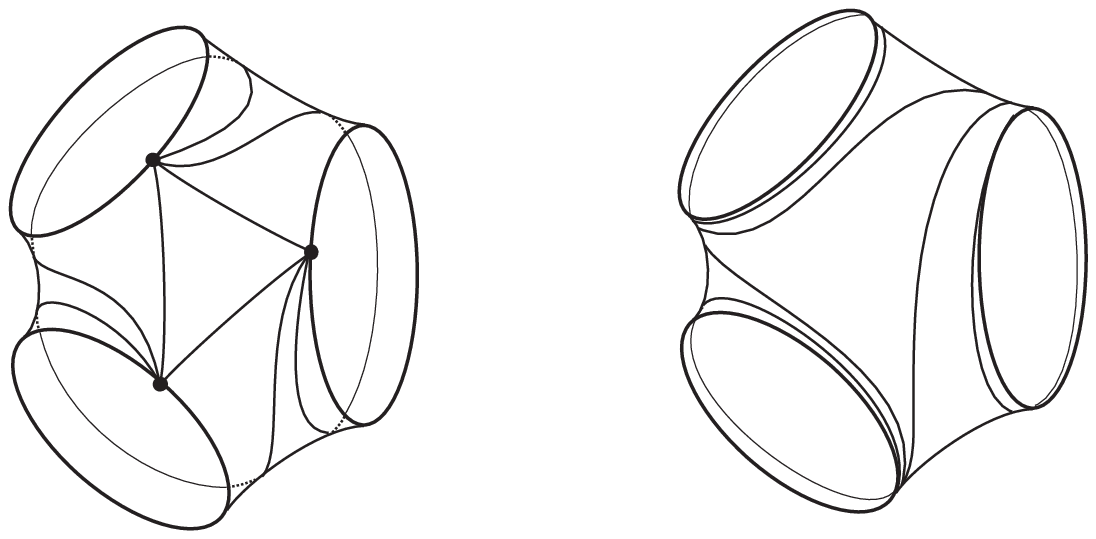}{A triangulation of a pair of pants and
  the lamination obtained by spinning it -- moving the vertices
  leftward along the boundaries and taking a limit}

In particular we find that $\pleat_\rho(v)$ is non-empty for any
lamination $v$ consisting of simple closed curves. The same is true
for any {\em finite-leaved} laminations, i.e. one in which we also
allow infinite arcs that either spiral on closed curves, or go out the
cusps of $S$. See \cite{ceg}.

\subsection*{Halfway surfaces}
\label{halfway}
Let $v$ and $v'$ be pants decompositions that agree except on one
curve: Thus $v = v_0 \union w$ and $v'=v_0 \union w'$ where $v_0$
divides $S$ into three-holed spheres and one component $W$ of type
$\xi(W)=4$. Assume also that $v$ and $v'$ intersect a minimal number
of times in $W$ (once if $W=S_{1,1}$ and twice if $W=S_{0,4}$). We say
that $v$ and $v'$ differ by an {\em elementary move on pants decompositions}.

As in \cite{minsky:boundgeom},  define a lamination $\lambda_{v,v'}$ as
follows: $\lambda_{v,v'}$ contains $v_0$ as a sublamination;
in each 3-holed sphere complementary component of $v_0$ it is the same
as the lamination in Figure \ref{spin};
and in $W$,  $\lambda_{v,v'}$
is given by the following diagram in the $\R^2$-cover of $W$: The
lift of $\boundary W$ is the lattice of small circles at integer
points. $W$ is the quotient under the action of $\Z^2$ if $W=S_{1,1}$,
and under the group generated by $(2\Z)^2$ and $-I$ if $W=S_{0,4}$.
The standard generators of $\Z^2$ in the figure correspond to the
curves $w$ and $w'$. 

\realfig{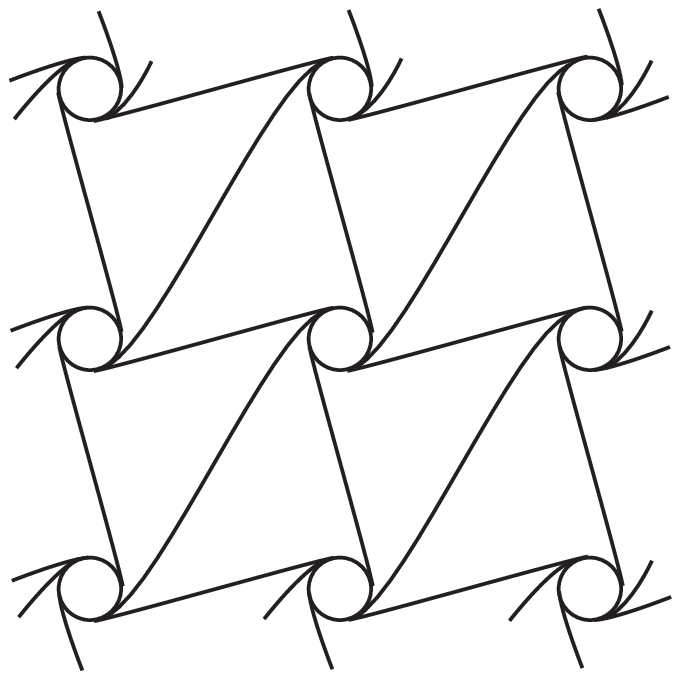}{The lamination
  $\til\lambda_{v,v'}$
in the planar cover of $W$.}

There is then a unique isotopy equivalence class 
$g_{v,v'}\in \pleat_\rho(\lambda_{v,v'})$, which we call the {\em halfway
  surface} associated to the pair $(v,v')$.

Let $g_v\in\pleat_\rho(v)$ and $g_{v'}\in\pleat_\rho(v')$ be leftward
pleated surfaces. Figure \ref{twopleats} indicates the laminations
$\lambda_v$ and $\lambda_{v'}$
associated to these surfaces, restricted to $W$ and lifted to the
planar cover as above. In the complement of $W$, all three laminations
agree. 

Consider the set of leaves labeled $\til l$ in $\til\lambda_{v}$ in
Figure \ref{twopleats}. These are disjoint from the lifts of $w$
and project to either one or two leaves $l$ in $W$ (depending on whether
$W$ is a 1-holed torus or 4-holed sphere) that spiral on its
boundary. The leaves $l$ are common to both $\lambda_v$ and
$\lambda_{v,v'}$. Furthermore it is easy to see that any closed curve $\alpha$
that has an essential intersection with $\gamma_w$ must also have an
essential intersection with $l$. If $\alpha$ is contained in $W$ this
is evident from the fact that $l$ cuts $int(W)$ into an
annulus with core $\gamma_w$; if not then $\alpha$ must cross
$\boundary W$, and we use the fact that $l$ spirals around every
component of $W$.
The corresponding facts hold for the leaves
$\til l'$ which project to $l'$ in $\lambda_{v'}$. 

\realfig{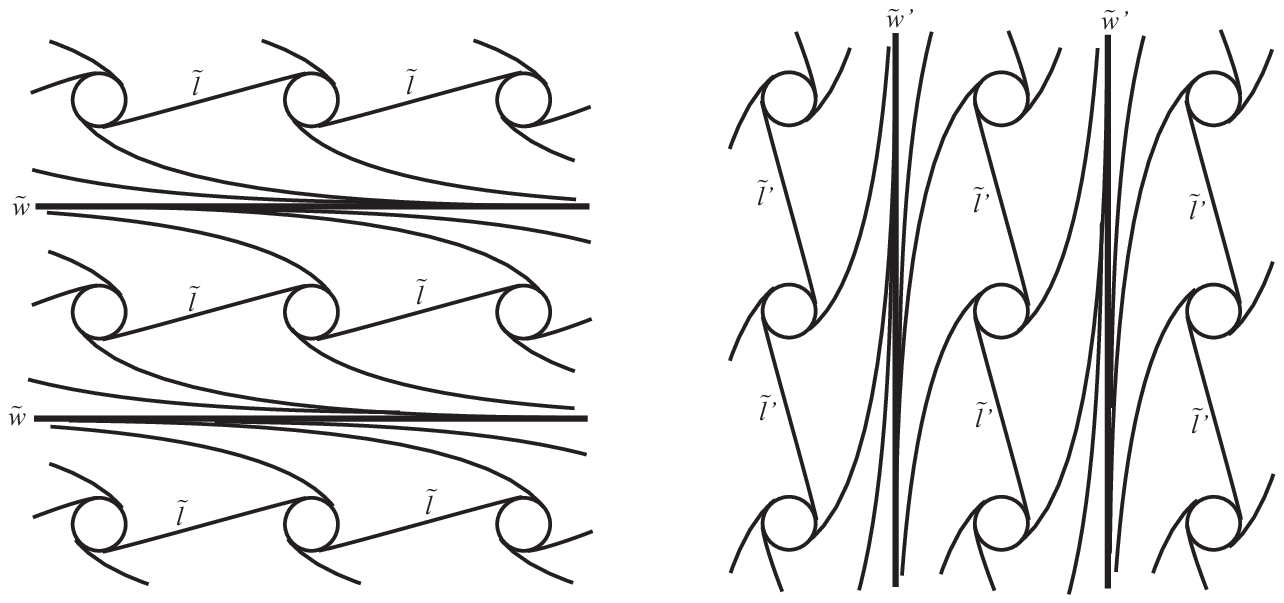}{The laminations associated to
  $g_v$ and $g_{v'}$. The representative of $w$ lifts to a horizontal
  line in the left figure, and  the representative of $w'$ lifts to a
  vertical line in the right figure.}

The following lemma, a restatement of Lemma 4.2 in \cite{minsky:boundgeom}, serves to
control the geometry of a halfway surface. 

\begin{lemma}{halfway bounds}
  Let $v=v_0\union w$ and $v'=v_0\union w'$ be pants decompositions
  that differ by an elementary move. Let $\sigma_v$, $\sigma_{v'}$ 
  and $\sigma_{v,v'}$ be the
  metrics induced on $S$ by the pleated surfaces $g_v$, $g_{v'}$ and  $g_{v,v'}$,
  respectively. Then for a constant $C$ depending only on the topology
  of $S$, 
  $$
  \ell_{\sigma_v}(w)  \le \ell_{\sigma_{v,v'}}(w) \le
  \ell_{\sigma_v}(w) + C,
$$
and
$$
  \ell_{\sigma_{v'}}(w') \le \ell_{\sigma_{v,v'}}(w') \le \ell_{\sigma_{v'}}(w') + C,
  $$
\end{lemma}
Note that $\ell_{\sigma_v}(w) = \ell_\rho(w)$ since the curve
representing $w$ is mapped to its geodesic representative by $g_v$; and
similarly for $\ell_{\sigma_{v'}}(w')$. Hence the left-hand inequality
on each line is immediate. The right-hand inequalities follow from an
application of Thurston's Efficiency of Pleated Surfaces \cite{wpt:I}.

In \cite{minsky:boundgeom} this result is proved and used without any
assumption about the lengths $\ell_\rho(v)$ or $\ell_\rho(v')$. In this
paper we will only use it in the case where these lengths are already
bounded both above and below.

\subsection{Tubes and constants}
\label{thick thin}

\subsubsection{Hyperbolic tubes}
\label{hyperbolic tubes}
A {\em hyperbolic tube} is the quotient of an $r$-neighborhood of a geodesic
in $\Hyp^3$ by a translation.

Given $\lambda\in\C$ with $\Re \lambda>0$, and $r>0$, we define
$\MT(\lambda,r)$ to be the quotient of the open $r$-neighborhood of the
vertical line above $0\in\C$ in the upper half-space model of $\Hyp^3$
by the loxodromic $\gamma:z\mapsto e^\lambda z$. Let $\MT(\lambda,\infty)$
denote the quotient of $\Hyp^3 \union \C\setminus\{0\}$ by $\gamma$.
Any hyperbolic tube is isometric to some $\MT(\lambda,r)$, but we note
that the imaginary part of $\lambda$ is, so far, only determined
modulo $2\pi$. 

\subsubsection*{Marked boundaries}
We discuss now
how to describe the geometry of a hyperbolic tube $\MT$ in terms of
the structure of its boundary torus.

If $T$ is an oriented Euclidean torus, a {\em marking} of it is an
ordered pair $(\alpha,\beta)$ of homotopy classes of unoriented simple
closed curves with 
intersection number 1. There is a unique $t>0$ and $\omega\in\Hyp^2$
such
that $T$ can be identified with $\C/t(\Z+\omega\Z)$ by an
orientation-preserving isometry, so that the images of $\R$ and
$\omega\R$ are in the classes $\alpha$ and $\beta$, respectively.
The parameter $\omega$ describes the conformal structure of $T$ as a
point in the Teichm\"uller space $\Teich(T)\equiv \Hyp^2$.

The boundary torus of a hyperbolic tube $\MT$ inherits a Euclidean
metric and an orientation from $\MT$ (see \S\ref{cores and ends} for
orientation conventions), and it admits an
almost uniquely defined marking: Let $\mu$ denote the homotopy class
of a meridian of the torus (i.e. the boundary of an essential disk in
$\MT$) and let $\alpha$ denote a homotopy class in $\boundary\MT$
of simple curves homotopic to the core curve of
$\MT$. While $\mu$ is unique, $\alpha$ is only defined up to multiples
of $\mu$. Fixing such a choice $(\alpha,\mu)$, we obtain
{\em boundary parameters} $(\omega,t)$ as above. 

Given $\lambda\in\C$ we can determine a marking $(\alpha,\mu)$ of
$\boundary\MT(\lambda,\infty)$, in a way that relates
ambiguities in the choice of $\alpha$ to the freedom of adding $2\pi i$
to $\lambda$.
The boundary at infinity $\boundary\MT(\lambda,\infty)$ is the
quotient of $\C\setminus \{0\}$ by $z\mapsto e^\lambda z$. Using $e^z$
as the universal covering $\C\to \C\setminus\{0\}$, we obtain 
$\boundary\MT(\lambda,\infty)$ as the quotient $\C/(\lambda \Z + 2\pi i\Z)$. 
The line $i\R$ maps to the meridian, and we let $\alpha$ be the image
of $\lambda\R$. Thus adding $2\pi i$ to $\lambda$ corresponds to
twisting $\alpha$ once around $\mu$. We note that the corresponding
boundary parameters for $\boundary\MT(\lambda,\infty)$ (with the
standard metric and orientation inherited from $\C$)  are
$\omega = 2\pi i/\lambda$ and $t=|\lambda|$. 

The marking $(\alpha,\mu)$ at $\boundary\MT(\lambda,\infty)$
determines a unique 
marking via orthogonal projection, denoted also $(\alpha,\mu)$,  on
$\boundary\MT(\lambda,r)$,
and for this marked torus we also have
boundary data, $(\omega_r,t_r)$.

The next lemma allows us to recover
the length-radius parameters $(\lambda,r)$ from $(\omega_r,t_r)$, and
indeed to construct a tube realizing any desired boundary data. 

\begin{lemma}{torus boundary data}
Given $\omega\in\Hyp^2$ and $t>0$, 
there is a unique pair $(\lambda,r)$ with $\Re\lambda>0,r>0$
such that $\MT(\lambda,r)$ 
has boundary data $(\omega,t)$ in the marking
determined by $\lambda$. 
\end{lemma}

\begin{proof}
Let $X_r$ be the boundary of the $r$-neighborhood of the vertical
geodesic over $0\in\C$. This is a cone in the upper half-space, and 
in the induced metric is a Euclidean cylinder with circumference
$2\pi\sinh r$. We may identify the universal cover $\til X_r$
isometrically with $\C$, with deck translation $z\mapsto z+ 2\pi i
\sinh r$. 
Let $X_\infty = \C\setminus \{0\}$ with the Euclidean
metric of circumference $2\pi$, and identify 
the universal cover $\til X_\infty$
isometrically with $\C$, with deck translation $z\mapsto z+ 2\pi i$.

Let $\Pi_r : X_\infty \to X_r$ be the orthogonal projection map. 
It lifts to a map $\Phi_r:\til X_\infty \to \til X_r$ which, in our
coordinates, can be written
\begin{equation}\label{Phi r def}
\Phi_r(x+iy) = x \cosh r +  iy \sinh r.
\end{equation}
(up to translation).

For our torus with parameters $(\omega,t)$ the meridian length is
$t|\omega|$. Hence in order to build the right hyperbolic tube 
we must choose $r$ so that 
\begin{equation}\label{r from omega}
2\pi \sinh r = t|\omega|.
\end{equation}
It remains to determine $\lambda$. In $\til X_\infty$, we see
$\lambda$ as the translation that, together with the meridian
$2\pi i$, produces the quotient torus
$\boundary\MT(\lambda,\infty)$, with marking. The map $\Phi_r$ takes
$\lambda$ to a translation $\lambda'$ that, together with the meridian
$2\pi i \sinh r$, yields the marking for $\boundary \MT(\lambda,r)$.
Our boundary parameter $\omega$ is given by the ratio
$$\omega = \frac {2\pi i \sinh r }{\lambda'}$$
so that $\lambda' = 2\pi i \sinh r /\omega$, and hence we can define
\begin{align}\label{lambda from omega}
\lambda &= \Phi_r^{-1} \left(\frac{2\pi i \sinh r}{\omega}\right)\\
&= h_r\left( \frac{2\pi i}{\omega} \right)
\end{align}
where we define 
$$
h_r(z) \equiv \Phi_r^{-1}(z\sinh r) = \Re z\tanh r  + i\Im z.
$$
These parameters therefore yield the desired torus. Notice that  for
large $r$,  $\lambda$ and $2\pi i/\omega$ are nearly the same. 
\end{proof}

\subsubsection*{Parabolic tubes}
The parabolic transformation $z\mapsto z+t$, acting on the region of
height $> 1$ in the upper half space model of $\Hyp^3$, gives a
quotient that we call a (rank 1) parabolic tube, with boundary
parameters $(i\infty,t)$. This can be obtained as the geometric limit of
hyperbolic tubes with parameters $(\omega_n,t)$ with $\Im\omega_n\to
\infty$ and $\Re\omega_n$ bounded.

\subsubsection{Margulis tubes}
\label{margulis tubes}
Let $N_J$ for $J\in\R$ denote the region in $N$ where the injectivity
radius times 2 is in $J$. Thus
$N_{(0,\ep)}$ denotes the (open) $\ep$-thin part of a hyperbolic manifold
$N$, and  $N_{[\ep,\infty)}$ denotes the $\ep$-thick part. 

By the Margulis Lemma (see
Kazhdan-Margulis \cite{kazhdan-margulis}) or J\o rgensen's inequality
(see J\o rgensen \cite{jorgensen:inequality}, and Hersonsky
\cite{hersonsky:inequality} or Waterman \cite{waterman:inequality}
for the higher-dimensional case), there is
a constant 
$\ep_M(n)$ known as the Margulis constant for $\Hyp^n$, such that every
component of an $\ep$-thin part when $\ep\le \ep_M(n)$ is of a standard
shape. In particular, in dimensions 2 and 3 in the orientable case,
such a component is either an open tubular neighborhood of a simple closed
geodesic, or the quotient
of an  open horoball by a parabolic group of 
isometries of rank 1 or 2. We call these components ``$\ep$-Margulis
tubes''. In dimension 3, the former type are hyperbolic tubes as in
the previous 
subsection, and the rank 1 parabolic (or ``cusp'') tubes were
mentioned above as well. 
See
Thurston \cite{wpt:textbook} for a discussion of the thick-thin
decomposition.

If $\gamma$ is a non-trivial homotopy class of closed curves in a hyperbolic
manifold with length less than $\ep$, then it corresponds
to an $\ep$-Margulis tube, which we shall denote
by $\MT_\ep(\gamma)$.
If $\Gamma$ is a union of several curves we let $\MT_\ep(\Gamma)$
denote the union of their Margulis tubes.

The radius $r(\gamma)$ of a hyperbolic $\ep$-Margulis tube depends on its
complex translation length, $\lambda(\gamma)$. Brooks-Matelski
\cite{brooks-matelski} apply J\o rgensen's inequality to show that 
\begin{equation}\label{r lambda bound}
r(\gamma) \ge \log\frac{1}{|\lambda(\gamma)|} - c_1
\end{equation} 
where the constant $c_1$ depends only on $\ep$. Writing $\lambda =
\ell + i\theta$,  
we note that when $\ell$ is small $|\lambda|$ can still be
large. However, an additional pigeonhole principle argument in
Meyerhoff \cite{meyerhoff:volumes} implies
\begin{equation}\label{r ell bound}
r(\gamma) \ge \half \log\frac{1}{\ell(\gamma)} - c_2
\end{equation}
where again $c_2$ depends only on $\ep$. 

We also have  some definite separation
between Margulis tubes for different values of $\ep$. Given 
$\ep \le \ep_M(3)$ we have, for any $\ep'<\ep$, 
\begin{equation}\label{definite tube nesting}
dist(\boundary \MT_{\ep}(\gamma), \boundary \MT_{\ep'}(\gamma)) \ge 
\half \log \frac{\ep}{\ep'} - c_3
\end{equation}
(see \cite{minsky:torus} for a discussion).

\subsubsection*{Margulis constants}
For the remainder of the paper we
fix $\ep_0 \le \ep_M(3)$.

Thurston \cite{wpt:II} pointed out that for a $\pi_1$-injective
pleated surface $f:S\to N$, the thick part of $S$ maps to the thick
part of $N$. More precisely, there is a function
$\epshrink:\R_+\to\R_+$, 
depending on $S$, such that
$$
f(S_{thick(\ep)}) \subset N_{thick(\epshrink(\ep))}.
$$
(See also Minsky \cite{minsky:torus} for a brief discussion)
Define $\ep_1 = \epshrink(\ep_0)$.  We may also assume that
$\ep_1<1$. 

\subsubsection*{Bers constant} There is a constant $L_0$, depending
only on $S$ (see Bers \cite{bers:degenerating,bers:inequality} and
Buser \cite{buser:surfaces}),
such that any hyperbolic metric on $S$ admits a pants
decomposition of total length at most $L_0$. 
In fact we can and will choose $L_0$ 
so that such a pants decomposition includes all simple geodesics of
length bounded by $\ep_0$. This $L_0$ will also do as the constant in
Theorem \ref{EL exist}, since the sequence of bounded curves there is
obtained by taking shortest curves in a sequence of pleated surfaces.

We can refine this
slightly as follows: There is a function $\bersfcn(\ep)$ such that, if
$\gamma$ is a curve of $\sigma$-length at least $\ep$, then there is a
pants decomposition of total length at most $\bersfcn(\ep)$ which
intersects $\gamma$ (equivalently, does not contain $\gamma$ as a
component). To see this, start with a pants decomposition $w$ of length
bounded by $L_0$. If $\gamma$ is a component of $w$ consider the shortest
replacement for $\gamma$ that produces an elementary move on $w$. The
upper bound on $w$ and lower bound on $\gamma$ gives an upper bound on
this replacement curve.

\subsubsection{Collars in surfaces}
\label{collar defs}
Components of the thin part in a hyperbolic surface are annuli, and the
radius estimate (\ref{r lambda bound}) can be made much more
explicit. We note also that as a result of this estimate
there is a positive function $\neck(L)$ such that, if
$\gamma$ is a curve in a hyperbolic surface  which essentially
intersects a curve $t$ of $\sigma$-length at most $L$, then
$\ell_\sigma(\gamma) > \neck(L)$. This will be used in \S\ref{coarse
proof}. 

It will also be convenient to fix a standard construction of a {\em collar} 
for each nontrivial homotopy class of simple curves in a hyperbolic surface:

\begin{lemma}{collar properties}
Fix a constant $0<c\le 1$.
For every hyperbolic structure $\sigma$ on $int(R)$ and every
homotopically non-trivial simple closed curve $\gamma$ in $R$ there is
an open annulus $\collar(\gamma,\sigma)\subset R$ whose core is homotopic to
$\gamma$, so that the following holds:

\begin{enumerate}
\item If $\ell_\sigma(\gamma) = 0 $ then $\collar(\gamma,\sigma)$ is a
horospherical neighborhood of the cusp associated to
$\gamma$. Otherwise it is an open embedded tubular neighborhood of radius
$w=w(\ell_\sigma(\gamma))$
of the geodesic representative of $\gamma$. 
\item 
If $\ell_\sigma(\gamma)<\frac{c}{2}$ (including the cusp
case) then the length of each boundary component of
$\collar(\gamma,\sigma)$ is exactly $c$.
\item If $\beta$ and $\gamma$ are homotopically distinct, disjoint
curves,
then $\collar(\gamma,\sigma)$ and $\collar(\beta,\sigma)$ have
disjoint closures. Indeed the collars are at least $d$ apart, for a
constant $d>0$, if $\ell_\sigma(\gamma) \le c$.
\end{enumerate}
\end{lemma}

\begin{proof}
This is done by a small variation of a construction of Buser
\cite{buser:surfaces} (in \cite{minsky:boundgeom} we used a slightly
different  variation).

Let 
$$
w_0(t) = \sinh^{-1}\left(\frac{1}{\sinh(t/2)}\right)
$$
and 
$$
w_c(t) = \cosh^{-1}\left(\frac{c}{t}\right)
$$
(the latter is defined only for $t\in(0,c]$).
Assuming for the moment that $\ell_\sigma(\gamma)>0$, 
define $\collar_i(\gamma,\sigma)$ to be the $w_i(\ell_\sigma(\gamma))$
neighborhood of the  geodesic representative of $\gamma$, where $i=0$
or $c$.  (If $\gamma$ is isotopic to a boundary component of $R$, we
take this neighborhood in the metric completion of $int (R)$, and
then intersect with $int(R)$.)

Buser shows that $\collar_0(\gamma,\sigma)$ is always an embedded
open annulus, and such collars are disjoint if non-homotopic.
On the other hand when $\ell_\sigma(\gamma) < c$, 
$\collar_c(\gamma,\sigma)$ has boundaries in $int(R)$ of length
exactly $c$, since the boundary length of a collar of width $w$
and core length $\ell$ is $\ell\cosh(w)$. 
We will obtain our desired collars by interpolating between these
two. 

Let $\delta(t) = w_0(t)-w_c(t)$. One can check that $\delta(t)$ is a
positive function for $t\le c/2$, and bounded away from $0$ and $\infty$.
Thus, for $\ell_\sigma(\gamma)<c/2$, $\collar_c(\gamma,\sigma)$
is a subannulus of
$\collar_0(\gamma,\sigma)$, and their boundaries are separated by a
definite but bounded distance. 

Thus we define 
$$
w(t) = \begin{cases}
        w_c(t) & t\le c/2 \\
        \alpha w_0(t) & t \ge c/2
       \end{cases}
$$
where
$\alpha=w_c(c/2)/w_0(c/2) < 1$. 
Define $\collar(\gamma,\sigma)$ to be the
$w(\ell_\sigma(\gamma))$-neighborhood 
of the geodesic representative. This collar has all the desired
properties. In particular the definite separation in part (3) follows
from the lower bound on $\delta$ when $\ell_\sigma(\gamma) < c/2$, and
from the fact that $\alpha<1$ when $\ell_\sigma(\gamma) \in [c/2,c]$.

Collars in the cusp case are easily seen to be obtained 
from these in the limit as $\ell(\gamma)\to 0$. In particular one may
check that $\delta(t)$ converges to $\log(2/c)$ as $t\to 0$,
and $\collar_0$ has boundary length 2 in the limit.
\end{proof}

\subsubsection*{Collar normalization}
In the remainder of the paper we will assume that the constant $c$ is
equal to the constant $\ep_1<1$. Note that this means that
$\collar(\gamma,\sigma)$ is contained in the component of the
$\ep_1$-thin part associated with $\gamma$. 
If $\sigma$ is implicitly understood then we may write
$\collar(\gamma)$. For a system $\Gamma$ of disjoint, homotopically distinct
simple closed curves we let $\collar(\Gamma,\sigma)$ or
$\collar(\Gamma)$ be the union of the collars of the components.

\subsection{Convention on isotopy representatives}
\label{isotopy convention}
Although we usually think of curves and subsurfaces in terms of their
isotopy classes, it will be useful for our constructions and arguments
to fix explicit representatives. Thus
we will adopt the following convention for the remainder
of the paper. We fix once and for all an oriented surface $S=S_{g,b}$.
Let $\hhat S$ denote a separate copy of $int(S)$ and fix a complete,
finite-area hyperbolic metric $\sigma_0$ on $\hhat S$. Embed $S$
inside $\hhat S$ 
as the complement of $\collar(\boundary S,\sigma_0)$. 

Now if $v$ is an essential homotopy class of simple closed curves or simple
properly embedded arcs in $S$, 
we let $\gamma_v$ denote its geodesic representative with respect to
$\sigma_0$ (occassionally we conflate $v$ and $\gamma_v$ when this is
convenient). If $v$ is a class of arcs then $\gamma_v$ is an infinite
geodesic whose ends exit the cusps of $\hhat S$.

We let $\collar(v)$ denote $\collar(v,\sigma_0)$ in $\hhat S$, and
we assume from now on that 
every open annulus is of the form
$\collar(v)$
(and every closed annulus is the closure of such a collar). 
Similarly every other subsurface (including $S$ itself) is a component of 
$$S\setminus \collar(\Gamma)$$ for a system of curves $\Gamma$.

This has the property that if two isotopy classes of subsurfaces have
disjoint representatives then these chosen representatives are already
disjoint (except for a closed collar that may share a boundary
component with an adjacent non-annular surface). 

For a nonannular surface $Y$ it is also true that all its
intersections with any closed geodesic are essential, since its
boundary is concave.

\subsection{The augmented convex core}
\label{augmented core}

Let $C^r_N$  be the closed $r$-neighborhood of the convex
core of any hyperbolic 3-manifold $N$, where $r\ge 0$. Define the {\em
augmented core of $N$} to be
$$
\hhat C_N = C^1_N \union N_{(0,\ep_0]}.
$$

Assume now that $\boundary \bbar N$ is incompressible. 
Let 
$$
\Pi_r : \overline N \to C^r_N
$$
denote the orthogonal projection (as in Epstein-Marden
\cite{epstein-marden}). Setting $r=1$ and restricting $\Pi_1$ to 
$\boundary\overline N = \boundary_\infty N$ we obtain a
map 
$$p_\infty : \boundary_\infty N \to \boundary C^1_N$$
which is shown in \cite{epstein-marden} to be $K_0$-bilipschitz
$C^1$-diffeomorphism  for a universal
$K_0$, where the metric on $\boundary_\infty N$ is the Poincar\'e
metric, denoted $\sigma_\infty$. 
Let
$$
\hat p: \boundary \hhat C_N \to \boundary C^1_N
$$
denote the restriction of $\Pi_1$ to $\boundary\hhat C_N$.
We will use these maps to compare the geometry of $\boundary_\infty N$
to that of $\boundary \hhat C_N$.  First let us define a modified
metric $\sigma^m$ on $\boundary_\infty N$ as follows: 

Let $\Gamma$ denote the
homotopy classes of curves in $\boundary_\infty N$ whose
$\sigma_\infty$-length is at most $\ep_1$
(note that this includes curves homotopic to cusps), and let
$\collar(\Gamma,\sigma_\infty)$ be their 
standard collars in $\boundary_\infty N$. 
We will let $\sigma^m$ be conformally
equivalent to $\sigma_\infty$, with $\sigma^m/\sigma_\infty$ a
continuous function 
which is equal to 1 on $\boundary_\infty N \setminus
\collar(\Gamma,\sigma_\infty)$, 
and such that  each annulus of $\collar(\Gamma,\sigma_\infty)$ is flat
in the $\sigma^m$ metric (a {\em flat annulus} for us is an annulus
isometric to the product of a circle with an interval).
This defines $\sigma^m$ uniquely. 
The flat annuli have circumference exactly $\ep_1$, by the definition
of the collars. 

\begin{lemma}{infinity to aug core}
The map
$$
\hhat p^{-1} \circ p_\infty : \boundary_\infty N \to \boundary \hhat C_N
$$
is a locally $K_1$-bilipschitz homeomorphism, with respect to the metric
$\sigma^m$ on the domain and the induced path metric from $N$ on the
range. 
\end{lemma}

\begin{proof}
Our result will follow from the following assertions.
\begin{enumerate}
\item For any $x\in \boundary_\infty N$ and $u\in T_x\boundary_\infty N$, 
$$
\frac{|u|_{\sigma_\infty}}{|u|_{\sigma^m}} \asymp
\min(1,inj_{\sigma_\infty}(x) )
$$
\item $\hat p$ is bijective.
\item If $z$ is a smooth point of $\boundary\hhat C_N$ and
$v\in T_z\boundary \hhat C_N$,
$$
\frac{|d\hat p(v)|}{|v|} \asymp \min(1,inj_{\boundary C^1_N}(\hhat p(z))).
$$
\end{enumerate}
Here by $A\asymp B$ we mean that $C^{-1}<A/B < C$ for a universal
constant $C$. 

Given $x\in \boundary_\infty N$ let $y=p_\infty(x)$ and let $z=\hhat
p^{-1}(y)$, uniquely defined by virtue of (2),
and assuming $z$ is a smooth point let $v=(d\hhat p)^{-1}dp_\infty(u)$
($d\hhat p$ is invertible by (3)).
We have already noted the fact from \cite{epstein-marden} that $p_\infty$ is
a $K_0$-bilipschitz diffeomorphism from $(\boundary_\infty N,\sigma_\infty)$ to
$\boundary C^1_N$.  Hence
the right-hand sides of (1) and (3)
are within uniform ratio, and $|u|_{\sigma_\infty}\asymp |dp_\infty(u)|
= |d\hhat p(v)|$. It follows immediately that $|u|_{\sigma^m} \asymp |v|$, 
so we have upper and lower bounds on 
$d(\hat p^{-1}\circ p_\infty)$ whenever $z$
is a smooth point of $\hhat C_N$. Since the failure of smoothness
occurs only on the transverse intersection points of
$\boundary N_{thin(\ep_0)} \intersect\boundary C^1_N$, which is
1-dimensional, we obtain the desired bilipschitz bounds.

Statement (1) follows from the definition of $\sigma^m$: outside of
$\collar(\Gamma,\sigma_\infty)$ the injectivity radius is bounded
uniformly below and $\sigma_\infty/\sigma^m = 1$. 
Inside the collar the injectivity radius is approximated by the length
of an equidistant curve from the core (or a horocycle for a cusp
collar), and the  rescaling of $\sigma^m$ is chosen to give this curve
the constant length $\ep_1$. Thus again the estimate follows. 

Before we establish statements (2) and (3), 
we consider the geometry of the map $\Pi_r$.
Let $X^r$ denote the preimage of $C^r_N$ in $\Hyp^3$. 
Then $X^0$ is the convex hull of the limit set. Letting
$\delta:\Hyp^3\to \R$ be the distance function $\delta(x) = d(x,X^0)$, 
we have $X^r = \delta^{-1}([0,r])$. In \cite{epstein-marden} it is
shown that, outside $X^0$, $\delta$ is  $C^1$ and  1-Lipschitz. The
gradient lines  of $\delta$ are geodesics that meet each
$\boundary X^r$ ($r>0$) orthogonally, and  $\Pi_r$ is defined in
the exterior of  
$X^r$ by following the gradient line back to its intersection with $\boundary
X^r$. (Here we are allowing $\Pi_r$ to denote the
map in the universal cover as well.)
These lines meet infinity in the
domain of discontinuity of the group 
and this allows the definition of $\Pi_r$ to be extended to the
boundary at infinity. Note that the gradient ray based at
any point in $\boundary X^1$ 
must intersect the lift of $\boundary \hhat C_N$, and hence $\hhat
p$ is surjective (half of (2)).

let $\til\MT$ be a component of the preimage $\til N_{thin(\ep_0)}$
of $N_{thin(\ep_0)}$ in
$\Hyp^3$. 
If the stabilizer of $\til\MT$ is hyperbolic then $\til\MT$ is an
$R$-neighborhood (for some $R$) of its geodesic axis $L$. 
Let $\beta(x) = d(x,L)$ be the distance function to $L$, so that
$\til\MT= \beta^{-1}([0,R))$. If $\til\MT$ has a parabolic stabilizer then 
it has a Busemann function, namely a function $\beta$ constant on
concentric horospheres
in $\til\MT$ and measuring signed distance between them, so that
$\til\MT=\beta^{-1}((-\infty,R))$. 

We claim now that, for any $x\in \boundary\til\MT$ outside $X^0$,
\begin{equation}\label{angle bound}
\del \delta(x) \cdot \del\beta(x) \ge \tanh \delta(x)/2. 
\end{equation}
To see this, let $y=\Pi_0(x)$. Then $\del\delta(x)$ is the outward
tangent vector to the geodesic arc $[y,x]$ at $x$.
$\del\beta(x)$ is the outward tangent vector to a geodesic $[z,x]$ at
$x$, where $z$ is either a point on the axis $L$ of $\til\MT$ or the
parabolic fixed point of $\til\MT$. In either case $z$ is in the closure
of $X^0$ in $\overline \Hyp^3$. 
Let $P$ be the 
plane through $y$ orthogonal to $[y,x]$. By convexity,
all of $X^0$, and in particular $z$, lies on the closure of the 
side of $P$ opposite from $x$. The lowest possible value for 
$\del\delta \cdot \del\beta$ is therefore obtained when $z\in
\boundary P$, and for this configuration the value $\tanh \delta(x)/2$
is obtained easily from hyperbolic trigonometry. This proves
(\ref{angle bound}).

In particular  we may conclude
that $\hhat p$ is injective. For if we follow a gradient line of
$\delta$ upward from $X^1$, the inequality 
$\del\delta\cdot\del\beta > 0$  implies that we can only exit 
$\til N_{thin(\ep_0)}$,
not enter it. Thus the gradient line meets $\boundary\til N_{thin(\ep_0)}$
at most once. This gives statement (2).

The derivative map $d\Pi_1$ is zero on $\del\delta$, and for
$v\in(\del\delta)^\perp$ we have
\begin{equation}\label{exponential contraction}
\frac{|(d\Pi_1)_x(v)|}{|v|} \asymp e^{-\delta(x)}
\end{equation}
for any $x\notin X^1$
(see \cite{epstein-marden}, especially the matrices on p. 141 and the
proof of Theorem 2.3.1). 

Since for $x\notin X^1$ we have $\delta(x)>1$, 
the lower bound (\ref{angle bound}) on $\del\delta\cdot\del \beta$ at
$x$ implies that (\ref{exponential contraction}) holds  
for $v\in (\del\beta)^\perp$ as well (with a different constant).

Now we are ready to 
establish statement (3). Consider first when $z$ is at least
$\ep_0$ away from $\boundary C^1_N$. Then $z$ is in the boundary of an
$\ep_0$ Margulis tube $\MT$, so there is an essential loop of length less than
$2\ep_0$ in $\boundary \MT$ based at $z$, which must therefore also
be in $\boundary\hhat C_N$; hence
$inj_{\boundary\hhat C_N}(z) \asymp \ep_0$. 

By
(\ref{exponential contraction}) applied to $(\del\beta)^\perp =
T\boundary\hhat C_N$ in an $\ep_0$-neighborhood of $z$, we conclude
that a shortest loop at $z$ has $\hat p$-image bounded by a uniform multiple
of $e^{-\delta(z)}$. Arguing similarly with a shortest loop at $y =
\hat p(z)$ (and lifting back using injectivity of $\hat p$), we
conclude that the ratio of injectivity radii
$inj_{\boundary C^1_N}(y)/inj_z(\boundary\hhat C_N) \asymp
inj_{\boundary C^1_N}(y)/\ep_0$
is uniformly estimated by $|d\hat p|$ in any direction at $z$. This
gives statement (3).

If $z$ is within $\ep_0$ of $\boundary
C^1_N$ then, again by (\ref{exponential contraction}) we have 
uniform upper and lower bounds on $d\hat p$ at $z$. In this case, 
since $d(z,y)\le \ep_0$, we conclude that $y$ is in the $\ep'$-thick
part of $N$, where $\ep'$ depends only on $\ep_0$. It follows that
$inj_{\boundary C^1_N}(y) \ge \ep'$, and again (3) follows.
\end{proof}



\subsubsection*{The exterior of the augmented core}
Let $E_N$ denote the closure of $N\setminus \hhat C_N$ in
$N$, and 
let $\bar E_N = E_N\union \boundary_\infty N$ be the closure of $E$
in $\bar N$. 

The previous discussion 
can also give us an explicit bilipschitz model for the geometry of
$E_N$, analogous to that of Epstein-Marden \cite{epstein-marden} for the
exterior of the convex core $C_N$.

Let $E_\nu$ denote a copy of $\boundary_\infty N \times [0,\infty)$,
  endowed with the metric 
$$
e^{2r}\sigma^m + dr^2
$$
where $r$ is a coordinate for the second factor.
We can also let $\bar E_\nu = \boundary_\infty N\times[0,\infty]$
where the ``boundary at infinity''
$\boundary_\infty E_\nu \equiv \boundary_\infty N \times\{\infty\}$
is endows with the conformal structure of $\boundary_\infty N$.
Note that the metric $\sigma^m$ is determined completely by the end
invariants $\nu$, justifying the notation
$E_\nu$ and $\bar E_\nu$.

The following lemma, taken together with the Lipschitz Model Theorem,
will give us the proof of the Extended Model Theorem (see
\S\ref{lipschitz}). 

\begin{lemma}{aug core exterior}
When $\boundary \bar N$ is incompressible, there is a homeomorphism
$$
\varphi : E_\nu \to E_N
$$
which is locally bilipschitz with uniform constant.
Furthermore $\varphi$ extends to a homeomorphism $\bar E_\nu
\to \bar E_N$ where the map on $\boundary_\infty E_\nu$
is given by $id_{\boundary_\infty N}$.
\end{lemma}

\begin{proof}
Using the notation of the proof of  Lemma \ref{infinity to aug core},
let $\varphi_0 = \hhat p^{-1}\circ p_\infty$. 
We define $\varphi$ to map $x\times[0,\infty)$ to the 
$\delta$ gradient line starting at $\varphi_0(x)$, by an
arclength-preserving map. 

First consider the map near $r=0$. By Lemma \ref{infinity to aug core} 
the map is locally $K_1$ bilipschitz on $\boundary_\infty N 
\times\{0\}$, where it 
restricts to $\varphi_0$ (in fact it is a diffeomorphism on the smooth
parts). The tangent plane at a smooth point $z$ of $\boundary \hhat C_N$
is $(\del \delta)^\perp$ if $z$ is in $\boundary C^1_N$, and 
otherwise $z\in \boundary \MT$ for a Margulis tube $\MT$, so
the tangent plane is $(\del \beta)^\perp$ for the associated distance
or Busemann function $\beta$. Inequality (\ref{angle bound}) then
tells us that in this case the tangent plane makes an angle with $\del
\delta$ that is bounded away from 0.  
It follows that we get uniform upper and lower bounds on $d\varphi$
at $r=0$.

For $r>0$, consider a smooth point $z\in \boundary \hhat C_N$ and let $q_r$ be the
map that slides points in a neighborhood of $z$ a distance $r$ along
the gradient lines of $\delta$. Thus $\varphi(x,r) =
q_r(\varphi_0(x))$ for $\varphi_0(x)$ near $z$. The derivative $dq_r$
is 1 in the $\del\delta$ direction, and expands by $e^r$ (up to a
bounded factor) in the $(\del\delta)^\perp$ direction (this follows
from the negative curvature of hyperbolic space, as in \cite{epstein-marden}).
It follows that the angle of $dq_r(T_z\boundary \hhat C_N)$ with $\del
\delta$ is even closer to $\pi/2$ than the angle between
$T_z\boundary\hhat C_N$ and $\del \delta$. Thus $d\varphi_{(x,r)}$
nearly preserves the orthogonality of the $\boundary_\infty N$ and
$\R$ directions, and is uniformly bilipschitz with respect to 
the metric $e^{2r}\sigma^m + dr^2$.

As in the proof of lemma \ref{infinity to aug core}, 
since the non-smooth points are a codimension 1 set, this gives
the desired global bilipchitz bounds.

The last statement on the extension of $\varphi$ to $\boundary_\infty E_\nu$
is immediate from the definition, since we have
mapped $y\times[0,\infty)$ for $y\in\boundary_\infty N$ to the
  $\delta$-gradient ray that terminates in $y$. 
\end{proof}

\section{Complexes of curves and arcs}
\label{complexes}

The definitions below are originally due to Harvey
\cite{harvey:boundary}, with some modifications in
\cite{masur-minsky:complex1,masur-minsky:complex2} and
\cite{minsky:boundgeom}.

The case of the one-holed torus $S_{1,1}$ and four-holed
sphere, $S_{0,4}$,  are special and will be
treated below, as will the case of the annulus $S_{0,2}$, which will 
only occur as a subsurface of larger surfaces.
We call all other cases ``generic''. 

If $S$ is generic, we define $\CC(S)$ to be the
simplicial complex whose vertices are
non-trivial, non-peripheral homotopy classes of simple curves, and whose
$k$-simplices are sets $\{v_0,\ldots,v_k\}$ of distinct vertices with
disjoint representatives. For $k\ge 0$ let $\CC_k(S)$ denote the
$k$-skeleton of $\CC(S)$. 

We define a metric on $\CC(S)$ by making each simplex regular Euclidean
with sidelength 1, and taking the shortest-path  distance.
We will more often use the shortest-path distance in the 1-skeleton,
which we denote $d_{\CC_1(S)}$ -- but 
note that the path metrics on $\CC(S)$ and $\CC_1(S)$ are quasi isometric.
These conventions also apply to the non-generic cases below.
\subsubsection*{One-holed tori and 4-holed spheres}
If $S$ is $S_{0,4}$ or $S_{1,1}$, we define the
vertices $\CC_0(S)$ as before, but let edges denote pairs
$\{v_0,v_1\}$ which have the minimal possible geometric intersection
number (2 for $S_{0,4}$ and $1$ for $S_{1,1}$.)

In both these cases, $\CC(S)$ is isomorphic to the
classical Farey graph in the plane (see e.g. Hatcher-Thurston
\cite{hatcher-thurston} or  \cite{minsky:taniguchi}).

\subsubsection*{Arc complexes}
If $Y$ is a non-annular surface with non-empty boundary, let us
also define the larger {\em arc complex} 
$\AAA(Y)$ whose vertices are 
either properly embedded essential arcs in $Y$, up to homotopy keeping
the endpoints in 
$\boundary Y$, or essential closed curves up to homotopy.
Simplices again correspond to sets of vertices with  disjoint representatives.

\subsubsection*{Subsurface projections}
Note that the vertices 
$\AAA_0(S)$ can identified with a subset of the geodesic lamination
space $\GL(S)$ -- a geodesic leaf with ends in the cusps of a complete finite
area hyperbolic structure on $int(S)$ determines a homotopy class of
properly embedded arcs in $S$.
Let $Y$ be a non-annular essential subsurface of $S$.
We can define a ``projection'' 
$$\pi_Y:\GL(S) \to \AAA(Y)\union\{\emptyset\}$$
as follows:

If $\alpha\in \GL(S)$ has no essential intersections with $Y$
(including the case that $\alpha$ is homotopic to $\boundary Y$) then
define $\pi_Y(\alpha) = \emptyset$.
Otherwise, $\alpha \intersect Y$ is a collection of disjoint essential curves
and/or properly embedded arcs (this follows from our use of geodesic
representatives of isotopy classes), which therefore span a simplex
in $\AAA(Y)$. Let $\pi_Y(\alpha)$ be an arbitrary choice of vertex of
this simplex. 

For convenience we also extend the definition of $\pi_Y$
to $\AAA_0(Y)$, where it is the identity map.

Let us denote 
$$d_Y(\alpha,\beta) \equiv
d_{\AAA_1(Y)}(\pi_Y(\alpha),\pi_Y(\beta)),
$$
provided $\pi_Y(\alpha)$ and $\pi_Y(\beta)$ are nonempty.
Similarly $\diam_Y(A)$ 
denotes $\diam_{\AAA_1(Y)}(\union_{a\in A} \pi_Y(a))$, where
 $A\subset \GL(S)$.

\subsubsection*{Annuli} 
If $Y$ is a closed annulus, 
let $\AAA(Y)$ be the complex whose
vertices are essential homotopy classes, rel endpoints, 
of properly embedded arcs, and whose simplices are sets of
vertices with representatives with disjoint interiors. Here it is
important that endpoints are {\em not} allowed to move in the boundary.

It is again easy to see that $\AAA(Y)$ is 
quasi-isometric to its 1-skeleton $\AAA_1(Y)$, and we will
mostly consider this. 

\subsubsection*{Twist numbers}
Fix an orientation of $Y$. 
Given vertices $a,b$  in $\AAA(Y)$, we will define their 
{\em twist number} $\twist_{Y}(a,b)$, which is a sort of rough signed
distance function satisfying
\begin{equation}\label{twist n dist}
|\twist_Y(a,b)|  \le d_{\AAA_1(Y)}(a,b) \le
|\twist_Y(a,b)| + 1.
\end{equation}

First, let
$\til\RR\subset\C$ be the strip $\{\Im z \in [0,1]\}$, 
$\RR$ its quotient by $z\mapsto z+1$, 
and identify 
$Y$ with $\RR$ by an orientation-preserving
homeomorphism. Choosing a lift $\til a$ in $\til \RR$ of a 
representative of $a$,
denote its endpoints as $x_a$ and $y_a+i$. 

Next define a function $\xi\mapsto\xi'$ on $\R$ by
letting $\xi'=\xi$ if $\xi\in\Z$ and $\xi' = n+1/2$ if $\xi\in(n,n+1)$
for $n\in\Z$. We can now define
\begin{equation}\label{twist definition}
\twist_Y(a,b) \equiv   (y_b-y_a)' - (x_b-x_a)'.
\end{equation}
This definition does not depend on any of the choices
made. An isotopy of $\boundary Y$ along $\boundary \RR$ corresponds
to changing the endpoints by a homeomorphism $h:\R\to\R$ satisfying $h(y+1) =
h(y)+1$, and one can check that $(h(y_2)-h(y_1))' = (y_2-y_1)'.$
Interchanging boundary components (by an orientation-preserving map of
$Y$) is taken care of by the identity
$(-\xi)'=-\xi'$, and choosing different lifts of $a$ and $b$
corresponds to the identity
$(\xi+1)'=\xi'+1$. It is also evident that
$\twist_Y(b,a) = -\twist_Y(a,b)$. 

The twist numbers are additive under concatenation of annuli:
Suppose that $Y$ is decomposed
as a union of 
annuli $Y_1$ and $Y_2$ along their common boundary and
$\alpha_1$ and $\alpha_2$ are two arcs connecting the boundaries of
$Y$ so that $\alpha_{ij} \equiv \alpha_i \intersect Y_j$
is an arc joining the boundaries of $Y_j$ for each $i=1,2$ and
$j=1,2$. The definition of twist numbers easily yields
\begin{equation}\label{additivity annuli}
\twist_Y(\alpha_1,\alpha_2) = \twist_{Y_1}(\alpha_{11},\alpha_{21})
+ \twist_{Y_2}(\alpha_{12},\alpha_{22}).
\end{equation}

One can furthermore check that 
$$
|(\xi+\eta)' - \xi' - \eta'| \le 1/2
$$
and this yields the approximate additivity property:
\begin{equation}\label{additivity arcs}
|\twist_Y(a,c) - \twist_Y(a,b) - \twist_Y(b,c)| \le 1
\end{equation}
for any three $a,b,c\in \AAA(Y)$. 

The inequality (\ref{twist n dist}) relating $\twist_Y$ to 
$d_{\AAA_1(Y)}$ can be easily verified from the definitions.

From (\ref{additivity arcs}) and (\ref{twist n dist})
we conclude that, 
fixing any $a\in\AAA_0(Y)$, the map $b \mapsto \tw_Y(a,b)$
induces a quasi-isometry from $\AAA(Y)$ to $\Z$. 
We can also use twisting to
define a notion of {\em signed length} for a geodesic: If $h$ is a
directed geodesic in $\AAA_1(Y)$ beginning at $a$ and terminating in
$b$,  we write
\begin{equation}\label{signed length}
[h] = \begin{cases}
        |h| &   \twist_Y(a,b) \ge 0, \\
        -|h| &  \twist_Y(a,b) < 0. 
      \end{cases}
\end{equation}
Note that 
$$[\overrightarrow{ab}] = -[\overrightarrow{ba}]$$
except when 
$\twist_Y(a,b) = 0$ and
$d_Y(a,b)=1$, in which case
$[\overrightarrow{ab}] = [\overrightarrow{ba}] = 1$.
This is an effect of coarseness.

\medskip

Now consider an essential annulus $Y\subset S$, and let us define
the subsurface projection $\pi_Y$ in this case. 

There is a unique
cover of $S$ corresponding to the inclusion
$\pi_1(Y)\subset\pi_1(S)$, to which we can append a boundary
using the circle at infinity
of the universal cover of $S$ to yield a closed annulus $\hhat Y$
(take the quotient of the compactified hyperbolic plane minus 
the limit set of $\pi_1(Y)$). 

We define $\pi_Y:\GL(S) \to \AAA(\hhat Y)\union \{\emptyset\}$ as
follows: Any lamination $\lambda$
that crosses $Y$ essentially lifts in $\hhat Y$ to a collection of
disjoint arcs,  some of which are essential. Hence we obtain a
collection of vertices of $\AAA(\hhat Y)$ any finite subset of which
form a simplex. 
Let $\pi_Y(\lambda)$ denote an arbitary choice of
vertex in this collection.

To simplify notation we often refer to $\AAA(\hhat Y)$ as
$\AAA(\alpha)$, where $\alpha$ is the core curve of $Y$.
We also let
$d_\alpha(\beta,\gamma)$ and $d_Y(\beta,\gamma)$
denote $d_{\AAA_1(\hhat Y)}(\pi_Y(\alpha),\pi_Y(\beta))$.
Similarly we denote
$\twist_{\hhat Y}(\pi_Y(\beta),\pi_Y(\gamma))$
by $\twist_\alpha(\beta,\gamma)$, and we also write $\pi_\alpha=\pi_Y$.

\subsubsection*{Projection bounds}
It is evident that the projections $\pi_W$ have the following
1-Lipschitz property: If $u$ and $v$ are vertices in $\AAA_0(S)$, 
both $\pi_W(u)\ne\emptyset$ and  $\pi_W(v)\ne\emptyset$, and 
$d_S(u,v) = 1$, then $d_W(u,v) \le 1$ as well. 

As a map from $\CC_0(S)$ to $\AAA_0(W)$, $\pi_W$ has the same
1-Lipschitz property when $\xi(S) > 4$. If $\xi(S)=4$ and $\xi(W)=2$
then the property holds with a Lipschitz constant of 3 (in
Lemma 2.3 \cite{masur-minsky:complex2} this is shown with a slight
error that leads to a constant of 2, but 3 is correct because two
curves in $S$ that intersect minimally can give rise to two arcs in
the annulus that intersect twice, and hence have distance 3. We are
grateful to Hideki Miyachi for pointing out this error).

A stronger contraction property applies to projection images of
geodesics, and plays an important role in the construction of
hierarchies in \cite{masur-minsky:complex2}:
\begin{lemma}{Geodesic Projection Bound} (Masur-Minsky
\cite{masur-minsky:complex2})  If $g$ is a (finite or infinite) geodesic
in $\CC_1(S)$ such that $\pi_W(v)\ne\emptyset$ for each vertex $v$ in
$g$, then
$$
\diam_W(g) \le A
$$
Where $A$ depends only on $S$. 
\end{lemma}

The complexes $\CC(W)$ and $\AAA(W)$ are in fact quasi-isometric when 
$\xi(W)\ge 4$. The inclusion $\iota:\CC_0(W)\to\AAA_0(W)$ has  a
quasi-inverse $\psi:\AAA_0(W)\to \CC_0(W)$ defined as follows:
On $\CC_0(W)$ let $\psi$ be the identity. 
If $a$ is a properly embedded arc in $W$ then the boundary of a
regular neighborhood  
of $a\union\boundary W$ contains either one or two essential curves, and we
let $\psi([a])$ be one of these (chosen arbitrarily).
In \cite[Lemma 2.2]{masur-minsky:complex2} we show that $\psi$ is a
2-Lipschitz map, with respect to $d_{\AAA_1(W)}$ and $d_{\CC_1(W)}$.

\subsection*{Hyperbolicity and Klarreich's theorem}
\label{hyperbolicity}

In Masur-Minsky \cite{masur-minsky:complex1}, we proved
\begin{theorem}{hyperbolicity thm}
$\CC(R)$ is $\delta$-hyperbolic.
\end{theorem}
A geodesic metric
space $X$ is $\delta$-hyperbolic if all triangles are
``$\delta$-thin''. That is, for any geodesic triangle 
$[xy]\union[yz]\union[xy]$, each side is contained in a
$\delta$-neighborhood of the union of the other two. 
The notion of $\delta$-hyperbolicity, due to Gromov \cite{gromov:hypgroups} 
and Cannon \cite{cannon:negative},
encapsulates some of the coarse properties of classical hyperbolic
space as well as metric trees and a variety of Cayley graphs of
groups. 
See also Alonso et al. \cite{short:notes}, Bowditch
\cite{bowditch:hyperbolicity}  and
Ghys-de la Harpe \cite{ghys-harpe}
for more about $\delta$-hyperbolicity.

In particular a $\delta$-hyperbolic space $X$ has  a well-defined {\em
boundary at infinity} $\boundary X$, which is roughly the set of
asymptotic classes of quasi-geodesic rays. There is a 
natural topology on $\bbar X\equiv X\union\boundary X$. When $X$
is proper (bounded sets are compact) $\bbar X$ is compact, but in our setting
$\CC(R)$ is not locally compact, and hence
$\boundary\CC(R)$ and $\bbar\CC(R)$ are not compact.
Since $\CC(R)$ and $\AAA(R)$ are naturally quasi-isometric,
$\boundary\CC(R)$ can be identified with $\boundary\AAA(R)$.

Klarreich showed in \cite{klarreich:boundary} that this boundary can
be identified with $\EL(R)$:

\begin{theorem}{EL is boundary}
{\rm (Klarreich \cite{klarreich:boundary})}
There is a homeomorphism 
$$k:\boundary\CC(R) \to \EL(R),$$
which is natural in the sense that a sequence
$\{\beta_i\in\CC_0(R)\}$ converges to $\beta\in\boundary\CC(R)$ if and
only if it converges to $k(\beta)$ in $\UML(R)$.
\end{theorem}
(Note that $\CC_0(R)$ can be considered as a subset of $\UML(R)$, and
hence the convergence $\beta_i\to k(\beta)$ makes sense.)

\section{Hierarchies}
\label{hierarchies}

In this section we will introduce the notion of hierarchies of tight
geodesics, and discuss their basic properties. Most of the material
here comes directly from \cite{masur-minsky:complex2}, although the
setting here is slightly more general in allowing infinite geodesics
in the hierarchy. Thus, although we will mostly state definitions and
facts, we will also need to indicate the changes to the arguments
in \cite{masur-minsky:complex2} needed to treat the infinite cases. 
In particular Lemma \lref{Existence of geodesics} is an existence
result for  infinite geodesics, and Lemma \lref{Resolution sweep} 
establishes a crucial property of ``resolutions by slices'' which was
essentially immediate in the finite case. 

\subsection{Definitions and Existence}
\label{hierarchy defs}

\subsection*{Generalized markings}
\label{markings}

The simplest kind of marking of a surface $S$ is a system of 
curves that makes a simplex in $\CC(S)$. In general we may want to include
twist information for each curve, and we also want to include
geodesic laminations instead of curves. 

Thus we define a {\em generalized marking} of $S$ 
as a lamination $\lambda\in\UML(S)$ together with a (possibly empty)
set of {\em transversals} $T$, where each element $t$ of $T$
is a vertex of $\AAA(\alpha)$ for a closed component $\alpha$ of
$\lambda$. For each closed component $\alpha$ of $\lambda$ there is at
most one transversal. 
The lamination $\lambda$ is called $\base(\mu)$, 
and the sublamination of $\lambda$ consisting of closed curves
is called $\simp(\mu)$. Note that it is a simplex of $\CC(S)$.

We further require that every non-closed component $\nu$ of $\base(\mu)$ 
is {\em filling} in the component $R$ of $S\setminus\collar(\simp(\mu))$
that contains it -- that is, $\nu\in\EL(R)$.

This is a generalization of the notion of
marking used in \cite{masur-minsky:complex2},
for which $\base(\mu)=\simp(\mu)$. We will call
such markings {\em finite}.

\subsubsection*{Types of markings}
When $\mu$ is a generalized marking let us call $\base(\mu)$ {\em
maximal} if it is not a proper subset of any element of $\UML(S)$.
Equivalently,
each complementary component $Y$ of $\simp(\mu)$ is either a
3-holed sphere or supports a component of $\base(\mu)$ which is in $\EL(Y)$. 
In particular the base of a finite marking is maximal 
if and only if it is a pants decomposition.

We call $\mu$ itself {\em maximal} if it is not properly contained in
any other marking -- equivalently, if its base is maximal, and if
every component of $\simp(\mu)$ has a transversal.

It will also be important to consider {\em clean markings}: a marking
$\mu$ is clean if $\base(\mu) = \simp(\mu)$, and 
if each transversal $t$ for a component $\alpha$ of $\base(\mu)$
has the form $\pi_\alpha(\bar t)$, where
$\bar t\in\CC_0(S)$ is disjoint from $\simp(\mu)\setminus\{\alpha\}$,
and where $\alpha\union\bar t$ fill a surface $W$
with $\xi(W)=4$ in which $\alpha$ and $\bar t$ are adjacent as
vertices of $\CC_1(W)$.
We will sometimes refer to $\bar t$ as the transversal to $\alpha$
in this case. 

If $Y\subset S$ is an essential subsurface and $\mu$ is a
marking in $S$ with $\base(\mu)\in\UML(Y)$, then we call $\mu$ a 
{\em marking in $Y$}. (Note that the transversals are not required to stay
in $Y$). If $Y$ is an annulus then a marking in $Y$ is just a simplex
of $\AAA(\hhat Y)$.

We can extend the definition of the subsurface projections $\pi_W$ to
markings as follows: We let $\pi_W(\mu) = \pi_W(\base(\mu))$ if the
latter is nonempty. If $W$ is an annulus and its core is a component
$\alpha$ of $\base(\mu)$ then we let $\pi_W(\mu)$ be the
transversal for $\alpha$, if one exists in $\mu$. In all other cases
$\pi_W(\mu)=\emptyset$. 

With these definitions, 
if $\base(\mu)$ is maximal, $\pi_W(\mu)$ is nonempty for any essential
subsurface  $W$ with $\xi(W) \ge 4$, and for any three-holed sphere
that is not a complementary component of $\simp(\mu)$.
If $\mu$ is maximal,
$\pi_W(\mu)$ is nonempty for all essential subsurfaces $W$ except
three-holed spheres that are components of $S\setminus \simp(\mu)$.

\subsection*{Tight geodesics.}
\label{tight}
A pair of simplices $\alpha,\beta$ in a  $\CC(Y)$ are
said to {\em fill } $Y$ if all non-trivial non-peripheral curves in $Y$ 
intersect at least one of $\gamma_\alpha$ or $\gamma_\beta$. If $Y$ is
an essential  subsurface of $S$ then it also holds that any curve
$\gamma$ in $S$ which essentially intersects a 
boundary component of $Y$ must intersect one of $\gamma_\alpha$ or
$\gamma_\beta$. 

Given arbitrary simplices $\alpha,\beta$ in $\CC(S)$, there
is a unique essential subsurface $F(\alpha,\beta)$ (up to isotopy)
which they fill: Namely, 
form a regular neighborhood of $\gamma_\alpha\union\gamma_\beta$,  and
and fill in all 
disks and one-holed disks. Note that $F$ is connected if and only
$\gamma_\alpha\union\gamma_\beta$ is connected. 

For a subsurface $X\subseteq Z$ let $\boundary_Z(X)$ denote the
{\em relative boundary} of $X$ in $Z$, i.e. those boundary components
of $X$ that are non-peripheral in $Z$.

\begin{definition}{tight seq def}
  Let $Y$ be an essential subsurface in $S$. If $\xi(Y)>4$,  a
  sequence of simplices 
  $\{v_i\}_{i\in\II} \subset \CC(Y) $ (where $\II$ is a finite or
  infinite interval in $\Z$) is called {\em tight} if
\begin{enumerate}
\item For any vertices $w_i$ of $v_i$ and $w_j$ of $v_j$ where $i\ne
  j$, $d_{\CC_1(Y)}(w_i,w_j) = |i-j|$,
\item Whenever $\{i-1,i,i+1\}\subset \II$, $v_i$ represents the relative
  boundary $\boundary_Y F(v_{i-1},v_{i+1})$.
\end{enumerate}

If $\xi(Y)=4$ then a tight sequence is just the vertex sequence of any
geodesic in $\CC_1(Y)$.

If $\xi(Y)=2$ then a tight sequence is the vertex sequence of any
geodesic in $\AAA(\hhat Y)$, with the added condition that the set of
endpoints on 
$\boundary\hhat Y$ of arcs representing the vertices equals the set of
endpoints of the first and last arc.
\end{definition}
Note that condition (1) of the definition specifies that given any
choice of components $w_i$ of $v_i$ the sequence $\{w_i\}$ is the
vertex sequence of a
geodesic in $\CC_1(Y)$. It also implies that $\gamma_{v_{i-1}}$ and
$\gamma_{v_{i+1}}$ always have connected union. 

In the annulus case, the restriction on endpoints of arcs is of little
importance, serving mainly to guarantee that there between any two
vertices there are only finitely many tight sequences.

With this in mind, a {\em tight geodesic} will
be a tight sequence together with some additional data:

\begin{definition}{tight geod def}
A {\em tight geodesic} $g$ in $\CC(Y)$ 
consists of a tight sequence
$\{v_i\}_\II$, and two generalized markings  in $Y$, $\I=\I(g)$ and $\T=\T(g)$,
called its {\em  initial} and {\em terminal} markings, such that:

If $i_0=\inf\II > -\infty$ then $v_{i_0}$  is
a vertex  of $\base(\I)$. If $i_\omega=\sup\II < \infty$ then
$v_{i_\omega}$ is a vertex of $\base(\T)$. 

If $\inf\II=-\infty$ then $\base(\I)$ is an element of $\EL(Y)$, and
$\lim_{i\to-\infty}v_i = \base(\I)$ in $\CC(Y)\union\boundary\CC(Y)$,
via the identification of Theorem \ref{EL is boundary}.
The corresponding limit holds for $\T$ if $\sup\II=\infty$.
\end{definition}
The length $|\II|\in[0,\infty]$
is called the length of $g$, usually written $|g|$.
We refer to each of 
the $v_i$ as {\em simplices} of $g$ (in \cite{masur-minsky:complex2}
we abused notation and called them ``vertices'', and in case $\xi(D(g))=4$
we may still do so).
$Y$ is called the {\em domain} or {\em support of $g$} and we write $Y=D(g)$.
We also say that $g$ is {\em supported in $D(g)$}.

If $Y$ is an annulus in $S$ then the markings $\I(g)$ and $\T(g)$ are
just simplices in $\AAA(\hhat Y)$. We can also define the signed
length $[g]$, as in (\ref{signed length}).

We denote the obvious linear order in $g$ as $v_i < v_j$ whenever
$i<j$. 

If $v_i$ is a simplex of $g$ define its {\em successor}
$$
\vsucc(v_i) = \begin{cases}
v_{i+1} & \text{$v_i$ is not the last simplex}\\
\T(g) & \text{$v_i$ is the last simplex}
\end{cases}
$$
and similarly define $\vpred(v_i)$ to be $v_{i-1}$ or $\I(g)$.

\subsection*{Subordinacy}
\label{subordinacy}

\subsubsection*{Restrictions of markings:}
If $W$ is an essential subsurface in $S$ and $\mu$ is a marking in $S$, then 
the {\em restriction}  of $\mu$ to $W$, which we write
$\mu\rest W$, is constructed from $\mu$ in the following way:
Suppose first that $\xi(W)\ge 3$. 
We let $\base(\mu\rest W)$ be the union of components of $\base(\mu)$
that meet $W$ essentially, and let the transversals of $\mu\rest W$ be
those transversals of $\mu$ that are associated to components of 
$\base(\mu\rest W)$.

If $W$ is an annulus ($\xi(W)=2$) then $\mu\rest W$ is just $\pi_W(\mu)$.

Note in particular that, 
if all the base components of $\mu$ which meet $W$ essentially are
actually contained 
in $W$, then $\mu\rest W$ is in fact a marking in $W$.
If $W$ is an annulus then $\mu\rest W$ is a marking in $W$ whenever it
is non-empty. 

\subsubsection*{Component domains:}
Given a surface $W$ with $\xi(W)\ge 4$ and a simplex $v$ in $\CC(W)$
we say that $Y$ is a 
{\em component domain of $(W,v)$} if either $Y$ is a component of
$W\setminus \collar(v)$, or $Y$ is a component of $\collar(v)$.
Note that in the latter case $Y$ is non-peripheral in $W$.

If $g$ is a tight geodesic with domain $D(g)$, 
we call $Y\subset S$ a {\em component domain of $g$} if
for some simplex $v_j$ of $g$, $Y$ is a component domain of 
$(D(g), v_j)$. 
We note that $g$ and $Y$ determine $v_j$ uniquely.
In such a case,
let
$$\I(Y,g) = \vpred(v_j)|_Y,$$
$$\T(Y,g) = \vsucc(v_j)|_Y.$$
Note in particular that these are indeed markings in $Y$ if nonempty,
except when 
$\xi(Y)=3$ (in which case they are just markings in $S$ whose bases
intersect $Y$) 

\medskip

If $Y$ is a component domain of $g$ and $\T(Y,g)\ne\emptyset$ then we
say that $Y$ is {\em directly forward subordinate} to $g$, or $Y\fsubd g$.
Similarly if
$\I(Y,g)\ne\emptyset$ we say that $Y$ is {\em directly backward
subordinate} to $g$, or $g\bsubd Y$. 
To clarify this idea let us note the following special cases when
$\T(Y,g)\ne \emptyset$:
\begin{enumerate}
\item $\xi(Y)\ge 4$: In this case $\vsucc(v_j)$ must have curves that
  are contained in $Y$. If $\vsucc(v_j) = \T(g)$ then $\T(Y,g)$ will
  contain the base components of $\T(g)$ that are contained in $Y$,
  together with their transversals if any.
\item $\xi(Y)=3$: Here we must have $\xi(D(g))=4$, and
   $v_j$ cannot be the last  vertex of $g$. $\T(Y,g)$ is the single
   vertex $v_{j+1}$.
\item $\xi(Y)=2$: If $v_j$ is not the last simplex of $g$ then 
    again we must have $\xi(D(g))=4$, and $\T(Y,g) = \pi_Y(v_{j+1})$.
    If  $v_j$ is the last simplex, then $\T(g)$  must
    contain a transversal for the core curve of $Y$. This transversal
    becomes $\T(Y,g)$. 
\end{enumerate}

\medskip

We can now define subordinacy for geodesics:

\begin{definition}{subordinate def}
If $k$ and $g$ are tight geodesics,
we say that $k$ is {\em directly forward subordinate} to $g$,
or $k\fsubd g$, provided
$D(k)\fsubd g$ and $\T(k) =
\T(D(k),g)$.  
Similarly we define  $g\bsubd k$ to mean $g \bsubd D(k)$ and $\I(k) =
\I(D(k),g)$. 
\end{definition}

We denote by {\em forward-subordinate}, or $\fsub$, 
the transitive closure of $\fsubd$,
and similarly for $\bsub$.
We let $h\fsubeq k$  denote the condition that $h=k$ or
$h\fsub k$, and similarly for $k\bsubeq h$.
We include the notation $Y\fsub f$ where $Y$ is a subsurface
to mean $Y\fsubd f'$ for some $f'$ such that
$f'\fsubeq f$, and similarly define
$b\bsub Y$.

\subsection*{Definition of Hierarchies}
\label{hierarchydefs}

\begin{definition}{hierarchy def}
A {\em hierarchy of geodesics} is a collection $H$ of
tight geodesics
in $S$ with the following properties:
\begin{enumerate}
\item
There is a distinguished {\em main geodesic} $g_H$ with domain $D(g_H)
= S$. The initial and terminal markings of $g_H$ are 
denoted also $\I(H), \T(H)$.

\item
Suppose $b,f\in H$, and  $Y\subset S$ is a subsurface with $\xi(Y)\ne
3$, such that  $b\bsubd Y$ and $Y\fsubd f$. Then $H$ contains a 
unique tight geodesic $k$ such that $D(k)=Y$, $b\bsubd k$ and $k\fsubd f$.

\item
For every geodesic $k$ in $H$ other than $g_H$, there are $b,f\in H$
such that  $b\bsubd k \fsubd f$.
\end{enumerate}
\end{definition}

\subsubsection*{Remark:} The notation here differs from that in
\cite{masur-minsky:complex2} only in the case of $\xi(Y)=3$. Here we
allow $\T(Y,g)$ and $\I(Y,g)$ to be nonempty, and hence $Y\fsubd f$
and $b\bsubd Y$ can occur, but we still explicitly
disallow $Y$ to be a domain of a geodesic in a hierarchy.

We will now investigate the structure of hierarchies, leaving the
question of their existence until \S\ref{hierarchies existence}.

\subsection{The descent theorem}
\label{descent theorem}

\subsubsection*{Footprints}
If $h$ is a tight geodesic in $D(h)$ and $Y$ is an essential subsurface of
$S$, we let the {\em footprint} of $Y$ on $h$, denoted $\phi_h(Y)$, 
be the set of simplices of $h$ that have no essential intersection
with $Y$. If $Y\subset D(h)$, the triangle inequality implies that
$\diam_{\CC_1(S)}(\phi_h(Y)) \le 2$. The condition of tightness implies
that, if $u,v,w$ are successive simplices of $h$ and
$u,w\in\phi_h(Y)$, then $v\in\phi_h(Y)$ as well. We remark that this is the
only place where the tightness assumption is used.

Thus $\phi_h(Y)$, if non-empty, is a subinterval of 1, 2 or 3
successive simplices in $h$. Let $\min\phi_h(Y)$ and
$\max \phi_h(Y)$ denote the first and last of these. 
Suppose $Y$ is a component domain of $(D(h),v)$
for a simplex $v$ of $h$. By definition, if 
$Y\fsubd h$, then $v = \max\phi_h(Y)$, and
similarly if $h\bsubd Y$ then $v = \min\phi_h(Y)$. 
An inductive argument then yields

\begin{lemma}{max footprint}
If $Y\fsub f \fsub h$ then
$$
\max \phi_h(Y) = \max \phi_h(D(f)).
$$
Similarly, if $h\bsub f \bsub Y$ then
$$
\min \phi_h(Y) = \min\phi_h(D(f)).
$$
\end{lemma}

This is used in \cite{masur-minsky:complex2} to control the structure
of the ``forward and backward sequences'' of a geodesic in a
hierarchy. That is, if $k\in H$ then by definition there exists 
$f\in H$ such that $k\fsubd f$. It can be shown that this $f$ is
unique, so that there is a unique sequence $k\fsubd f_0 \fsubd f_1
\fsubd \cdots\fsubd g_H$, and similarly in the backward direction.
The structure of these sequences is crucial to understanding and using
hierarchies. 

In particular it follows from Lemma \ref{max footprint} that
if $k\fsub f$ then $D(k)$ intersects $\T(f)$ nontrivially.
This condition can in fact capture all geodesics to which $k$ is
forward-subordinate.
Let us define, for any essential subsurface $Y$ of $S$,
\begin{equation}\label{Sigma plus def}
\Sigma^+_H(Y) = \{ f\in H: Y\subseteq D(f) \ \text{and}\  \T(f)|_{Y} \ne \emptyset\}
\end{equation}
and similarly
\begin{equation}\label{Sigma minus def}
\Sigma^-_H(Y) = \{ b\in H: Y\subseteq D(f) \ \text{and}\  \I(b)|_Y \ne
\emptyset\}.  
\end{equation}
(We write $\Sigma^\pm$ when $H$ is understood).
Theorem 4.7 of \cite{masur-minsky:complex2} is a central result in
that paper, and describes the structure of
$\Sigma^\pm(Y)$. We present here a slight extension of that theorem
which expands a little on the case of three-holed spheres.

\medskip

\begin{theorem+}{Descent Sequences}
Let $H$ be a hierarchy in $S$, and $Y$ any essential subsurface of $S$.

\begin{enumerate}
\sitem If $\Sigma^+_H(Y)$ is nonempty then it has the form 
$\{f_0,\ldots,f_n\}$ where $n\ge 0$ and 
$$f_0\fsubd\cdots\fsubd f_n=g_H.$$ 
Similarly, 
if $\Sigma^-_H(Y)$ is nonempty then it has the form 
$\{b_0,\ldots,b_m\}$ with $m\ge 0$, where
$$g_H=b_m\bsubd\cdots\bsubd b_0.$$

\sitem If  
$\Sigma^\pm(Y)$ are both nonempty and $\xi(Y)\ne 3$, then $b_0 = f_0$, and
$Y$ intersects every simplex of $f_0$ nontrivially.
\sitem If $Y$ is a component domain in any geodesic $k\in H$
$$f\in \Sigma^+(Y) \ \ \iff \ \ Y\fsub f,$$
and similarly, 
$$b\in \Sigma^-(Y) \ \ \iff \ \ b\bsub Y.$$

If, furthermore, $\Sigma^\pm(Y)$ are both nonempty,
then in fact $Y$ is the support of $b_0=f_0$. 

\sitem Geodesics in $H$ are determined by their supports. That is, if
  $D(h)=D(h')$ for $h,h'\in H$ then $h=h'$.

\end{enumerate}
\end{theorem+}

\begin{proof}
This theorem is taken verbatim from \cite{masur-minsky:complex2},
except that part (3) is stated there only for $\xi(Y)\ne 3$. 
When $\xi(Y)=3$, the definition in this paper of the relations
$Y\fsubd h$ and $h\bsubd Y$ is slightly different from the definition
in \cite{masur-minsky:complex2} (see \S\ref{hierarchy defs}), and
one can easily check that with this definition the proof 
in fact goes through verbatim in all cases (the relevant argument is
to be found in Lemma 4.21 of \cite{masur-minsky:complex2}).

\end{proof}

\subsubsection*{Completeness}
A hierarchy $H$ is {\em $k$-complete} if every
component domain $W$ in $H$ with $\xi(W)\ne 3$  and $\xi(W)\ge k$
is the domain of some geodesic in $H$. $H$ is {\em complete} if it  is
$2$-complete. 

If $\I(H)$ and $\T(H)$ are maximal markings, then 
${\I(H)}|_W$ and ${\T(H)}|_W$ are nonempty for every $W$ such that
$\xi(W)\ne 3$. Theorem \ref{Descent Sequences} part (3) then implies that 
$H$ is complete. 

If $\I(H)$ and $\T(H)$ only have maximal {\em bases} then 
${\I(H)}|_W$ is nonempty whenever $\xi(W)\ne 3$, 
except when $W=\collar(v)$ for a vertex $v\in\base(\I(H))$ which has
no transversal; and similarly for ${\T(H)}|_W$.  It follows that 
$H$ is $4$-complete.

\subsection*{Slices and resolutions}
\label{resolutions}

A {\em slice } of a hierarchy $H$ is a set $\tau$ of pairs $(h,v)$, 
where $h\in H$ and $v$ is a simplex of $h$, satisfying the following
properties: 
\begin{enumerate}
\item[S1:] A geodesic $h$ appears in at most one pair in $\tau$.
\item[S2:] There is a distinguished pair $(h_\tau,v_\tau)$ in $\tau$,
  called the bottom pair of $\tau$. We call $h_\tau$ the bottom geodesic. 
\item[S3:] For every $(k,w)\in \tau$ other than the bottom pair, $D(k) $
  is a component domain of $(D(h),v)$ for some $(h,v)\in\tau$.
\end{enumerate}

To a slice $\tau$ is associated a {\em marking}, $\mu_\tau$, whose
base is the union $\{v: (h,v)\in\tau, \xi(D(h))\ge 4\}$. The
transversal curves are the vertices $v$ for $(h,v)\in\tau$ with
$\xi(D(h)) = 2$. We say that $\tau$ is {\em saturated} if 
it satisfies
\begin{enumerate}
\item[S4:] Given $(h,v)\in\tau$,  if $k\in H$ has $D(k)$ equal to
a component domain of $(D(h),v)$ then
there is a pair $(k,w)\in\tau$.
\end{enumerate}
We say that $\tau$ is {\em full} if a stronger condition holds:
\begin{enumerate}
\item[S4':] Given $(h,v)\in\tau$,  if $Y$ is a
component domain of $(D(h),v)$ and $\xi(Y)\ne 3$ then
there is a pair $(k,w)\in\tau$ with $D(k)=Y$.
\end{enumerate}
(This terminology differs slightly from
\cite{masur-minsky:complex2}, but has the same content.
Also, the $\xi\ne 3$
condition was mistakenly left out of the condition corresponding to
(S4') in \cite{masur-minsky:complex2}.)

A full slice whose bottom geodesic is equal to the main geodesic
$g_H$ will be called {\em maximal}. Note that if $\tau$ is maximal
then $\mu_\tau$ is a maximal finite marking, and in particular
$\base(\mu_\tau)$ is a pants decomposition.

If $H$ is complete then every saturated slice is full.

\medskip
\subsubsection*{Elementary moves on slices}
We now define a {\em forward elementary
move} on a saturated slice $\tau$.
Say that a pair $(h,v)$ in $\tau$ is {\em forward movable} if:
\begin{itemize}
\item[M1:] $v$ is not the last simplex of $h$. Let $v'=\vsucc(v)$.
\item[M2:] For every $(k,w)\in \tau$ with $D(k)\subset D(h)$ and
$v'|_{D(k)}\ne \emptyset$, $w$ is the last simplex of $k$.
\end{itemize}
When this occurs we can obtain a slice $\tau'$  from $\tau$ by
replacing $(h,v)$ with 
$(h,v')$, erasing all the pairs $(k,w)$ that appear in condition (M2),
and inductively replacing them (starting with component domains of
$(D(h),v')$) so that the final $\tau'$ is saturated and satisfies
\begin{enumerate}
\item[M2':] For every $(k',w')\in \tau'$ with $D(k')\subset D(h)$ and
$v|_{D(k')}\ne \emptyset$, $w'$ is the first simplex of $k$.
\end{enumerate}
It is easy to see that $\tau'$ exists and is uniquely determined by
this rule. We write $\tau\to\tau'$, and say that the move
{\em advances} $(h,v)$ to $(h,v')$.

The definition can also be reversed: 
A pair $(h,v')$ in $\tau'$ is {\em backward movable} if it satifies 
\begin{itemize}
\item[M1':] $v'$ is not the first simplex of $h$
\end{itemize}
and the condition (M2'), with $v\equiv\vpred(v')$.
We can then construct $\tau$ by erasing pairs
appearing in (M2') and replacing them with pairs so that (M2) is
satisfied. 

These definitions come from 
\cite{masur-minsky:complex2}, where we consider elementary moves only for
maximal slices, but the definitions make sense in general.

\medskip

Let $V(H)$ denote the set of saturated slices of $H$ whose bottom
geodesic is $g_H$. We remark that $V(H)$ is nonempty, since starting
with any simplex of $g_H$ we can successively add pairs in component
domains to obtain a saturated slice. 
Note that if $H$ is complete then $V(H)$ is just the set of
maximal slices. If $H$ is $4$-complete, then every slice $\tau\in
V(H)$ has $\base(\mu_\tau)$ a pants decomposition (but $\mu_\tau$ may
be missing transversals).

For a complete finite hierarchy, we show in \cite{masur-minsky:complex2} the
existence of a {\em resolution}, which is a sequence 
$\{\tau_i\}_{i=0}^N$ in $V(H)$, such that $\tau_i \to \tau_{i+1}$ ,
$\tau_0$ is the initial slice of $H$, and $\tau_N$ is the terminal
slice. (The initial slice is the unique slice $\tau\in V(H)$ 
such that for every pair $(h,v)\in \tau$,
$v$ is the first simplex of $h$. The terminal slice is defined similarly).
To do this, we show that for each $\tau\in V(H)$
there is at least one elementary move $\tau\to\tau'$, unless $\tau$ is
the terminal slice. Beginning with the initial slice we
successively apply elementary moves until we obtain the terminal
slice. There is a certain partial  order $\sprec$ on $V(H)$ such that
$\tau\sprec\tau'$ whenever
$\tau\to\tau'$, and hence no $\tau$ can appear twice, and the process
terminates by finiteness. 

We wish to extend this idea to hierarchies that contain infinite
geodesics, and we also want to consider the case that $\I(H)$ and
$\T(H)$ may not be maximal.

In this setting we can consider
sequences $\{\tau_i\}_{i\in\II}$ in $V(H)$ where $\II\subset \Z$
is a possibly infinite interval. 
For every slice $\tau\in V(H)$, if it is  not terminal we will be able to find
at least one move $\tau \to \tau'$ and if it is not initial, at least one
$\tau''\to\tau$. Thus we can begin with any slice in $V(H)$
and apply moves in both the forward and backward direction.
There is only one pitfall: when there is an infinite geodesic whose
support is a proper subsurface, we must avoid making
infinitely many moves within that subsurface while the complement gets
``stuck''. Therefore we would like to consider sequences with the
following property:
\begin{itemize}
\item[R:] If $(h,v)$ is forward movable in $\tau_i$, then there
exists some
$j\ge i$ such that the move $\tau_{j}\to\tau_{j+1}$ advances
the pair $(h,v)$. Similarly if $(h,v)$ is backward movable
then there exists $j\le i$  such that the move
$\tau_{j-1}\to\tau_{j}$ advances from $(h,\vpred(v))$ to $(h,v)$.
\end{itemize}
We call an elementary move sequence $\{\tau_i\}$ satisfying this
condition a resolution. 

\begin{lemma}{Resolution exists}
For any hierarchy $H$ and any slice $\tau_0\in V(H)$, there is
a resolution of $H$ containing $\tau_0$. 
\end{lemma}

\begin{proof}
Fix $\tau_0$ in $V(H)$.  We define the rest of
the resolution sequence inductively.

Suppose $\tau_j$ has been defined for $0\le j \le i$. If
$(h,v)\in\tau_i$ is forward 
movable we have  seen that there exists a move $\tau_i\to \tau'$ that
advances $(h,v)$. If $(k,w)\in \tau$ is also forward movable, we
claim that $(k,w)$ is still in $\tau'$ and still forward
movable. $D(k)$ and $D(h)$ are either disjoint or
contained one in the other, since they are both in the same slice. If
they are disjoint then the claim is evident. 
Suppose $D(k)\subset D(h)$. Since $(k,w)$ is forward movable, $w$ is
not the last simplex of $k$, so 
since $(h,v)$ is forward movable and hence satisfies condition (M2),
we must have  
$\vsucc(v)|_{D(k)}=\emptyset$. Hence $(k,w)$ remains in $\tau'$, 
and so does any pair with domain in $D(k)$. It follows that $(k,w)$ is
forward movable in $\tau'$. Finally if $D(h)\subset D(k)$, the same argument 
tells us that $\vsucc(w)|_{D(h)}=\emptyset$. 
Since only domains in $D(h)$ are changed by the move, it follows that
$(k,w)\in\tau'$ and is still forward movable. Thus the claim is
estalished. 

We conclude that there is a sequence of elementary moves
$\tau_i\to\tau_{i+1}\to\cdots\to\tau_j$ which advances all of the
forward-movable pairs in $\tau$, one after the other. A similar
argument works in the 
backward direction, and we continue inductively, stopping (in each
direction) only if we reach a slice with no movable pairs.
The resulting sequence  (which may be infinite in either direction)
satisfies condition (R) and hence is a resolution.  
\end{proof}

The basic important property of the
resolution is that it sweeps through all parts of the hierarchy, in a
monotonic way:

\begin{lemma+}{Resolution sweep}
Let $H$ be a hierarchy and 
$\{\tau_i\}_{i\in\II}$ a resolution. For any pair $(h,v)$ with
$h\in H$ and $v$ a simplex of $h$, there is a slice $\tau_i$
containing $(h,v)$. 

Furthermore, if $v$ is not the last simplex of $h$
then there is exactly one elementary move $\tau_i\to\tau_{i+1}$ which
advances $(h,v)$.
\end{lemma+}

\begin{proof}
Let us first prove the following property of a resolution:
\begin{itemize}
\item[G:]  If $(k,w)\in\tau_i$ then there exists a finite interval
$\tau_{j}\ldots\tau_k$ in which all pairs $\{(k,w'): \text{$w'$ a simplex of
$k$}\}$ are obtained.
\end{itemize}
If $(k,w)$ is already forward movable
then condition (R) yields $j\ge i$ such that $\tau_j$ contains
$(k,\vsucc(w))$. Similarly if $(k,w)$ is  backward movable
then condition (R) yields $j\le i$ such that 
$\tau_j$ contains  $(k,\vpred(w))$.

Let us prove property (G) by induction on $\xi(D(k))$. If $\xi(D(k))$
is minimal among all $\{\xi(D(g)):g\in H\}$,
then $(k,w)$ is automatically forward movable whenever $w$ is not the
last simplex, and backward movable whenever $w$ is not the first
simplex (because conditions (M2) and (M2') hold vacuously). Thus we
can use condition (R) as in the previous paragraph to obtain
all simplices  of $k$ in finitely many steps.

Now in general,
suppose $w$ is not the last simplex of $k$, 
but that $(k,w)$ is not movable. Let $w'=\vsucc(w)$.
Then (M2) is violated so there is some $(p,u)\in\tau_i$ with
$D(p)\subset D(k)$, $w'|_{D(p)} \ne \emptyset$, 
and $u$ not the last simplex of $p$. 
For each such $p$, since $\phi_k(D(p))$ contains $w$ (by definition of
a slice), but not $w'$ (since $w'|_{D(p)} \ne \emptyset$), 
we must have $\T(k)|_{D(p)}\ne \emptyset$. Thus $k\in\Sigma^+(D(p))$,
and applying Theorem \ref{Descent Sequences} we have
$p\fsub k$. It follows that $p$ is not infinite in the
forward direction, for if it were, $\T(p)$ would be a lamination in
$\EL(D(p))$ and then, inductively, for each $f$ such that $p\fsub f$
we would have $\max \phi_f(D(p))$ equal to the last vertex of $f$ --
and this contradicts the fact that $w'\notin \phi_k(D(p))$. 

Thus, by the
inductive hypothesis we can advance such a $p$ in finitely many
moves to its last simplex. Apply this first to a $p$ with $D(p)$ a
component domain of $(D(k),w)$, reaching a slice $\tau_{i'}$
($i'\ge i$) containing 
$(p,u')$ with $u'$ the last vertex of $p$. Now for all $j\ge i'$, 
as long as $(k,w)$ has not been advanced to $(k,w')$, 
the pair $(p,u')$ must still be in
$\tau_j$. This is because no other move can affect $p$: since $D(p)$
is a component domain of $(D(k),w)$, there is no intervening $k'$ in
the slice with $D(p)\subset D(k') \subset d(k)$.
Thus we can apply the same argument to other component
domains of $(D(k),w)$ (if any) without removing $(p,u')$. Once these
component domains have had their pairs advanced to satisfy (M1), we 
repeat the argument
successively in {\em their} component domains. In finitely
many moves 
we will therefore reach a point where $(k,w)$ is forward movable, and we can
apply (R). The same
argument applies in the backward direction. 
Repeating this, we can obtain all the simplices of $k$, which
establishes (G).

\medskip

Now let $(h,v)$ be any pair, and let $\tau$ be any maximal slice (say,
$\tau=\tau_0$).
We claim that there is a unique pair $(k,w)\in\tau$ such that
$D(h)\subseteq D(k)$ and $\phi_k(D(h))$ does not contain $w$. We find
this pair by induction: let $(g,u)\in\tau$ be a pair with $D(h)\subseteq
D(g)$ (the bottom pair has this property). If $u\notin \phi_g(D(h))$
then $(g,u)$ satisfies our conditions, and no other pairs in $\tau$
with domain in $D(g)$ can  contain $D(h)$.  (Note, this case includes
$h=g$, in which case $\phi_g(D(h))=\emptyset$). If $u\in \phi_g(D(h))$
then let $W$ be the component domain of $(D(h),u)$ containing $D(h)$.
Since $\I(H)|_{D(h)}$ and $\T(H)|_{D(h)}$ are nonempty (by Theorem
\ref{Descent Sequences}), so are $\I(H)|_W$ and $\T(H)|_W$, and again
by Theorem \ref{Descent Sequences} $W$ must be the support of some
geodesic $g'\in H$. 
(A special case is when
$D(h)$ is an annulus with core a component of $u$, and then $g'=h$.)
Since $\tau$ is saturated, there is a pair
$(g',u')\in\tau$.
We can repeat the argument with $(g',u')$. Thus starting with the
bottom pair $(g_H,u)$ of $\tau$ we arrive at the unique $(k,w)$ as
claimed.

If $h=k$ then we can apply (G) to advance $\tau$ by forward or
backward elementary moves to a slice containing $(h,v)$, and we are done.

Suppose now that $D(h)$ is properly contained in $D(k)$.
Now since $w\notin \phi_k(D(h))$, either $\phi_k(D(h))=\emptyset$, or
(without loss of generality) assume $w<\min\phi_k(D(h))$. In either
case, $D(h)$ intersects the first vertex of $k$ so $k\in\Sigma^-(D(h))$
and hence $k\bsub h$ by Theorem \ref{Descent Sequences}. This in turn
rules out the possibility 
that $\phi_k(D(h)) = \emptyset$. Now again using property (G), we can
advance $\tau$ by elementary moves until we obtain a slice $\tau'$
containing $(k,x)$, where $x=\min\phi_k(D(h))$. The claim above
implies that there is a new pair $(k',w')$ in $\tau'$, with
$D(h)\subseteq D(k')\subset D(k)$, and $w'\notin\phi_{k'}(D(h))$. 
We can therefore repeat the argument with $(k',w')$ replacing
$(k,w)$. Since $\xi(D(k')) < \xi(D(k))$, the process must terminate in
a finite number of steps.

This proves the first statement of the Lemma, that every pair $(h,v)$
is obtained. Now consider a transition $(h,v)$ to $(h,v')$
where $v'=\vsucc(v)$. In
\cite{masur-minsky:complex2} we introduce a strict partial order
$\sprec$ on $V(H)$ (this is done there for a complete hierarchy, but
the proof carries through in our setting as well), which has the
following properties: First, if $\tau\sprec \tau'$ and
$(h,u)\in\tau$,
$(h,v)\in\tau'$, then $u$ precedes $v$ in $h$. Second, if
$\tau\to\tau'$ is an elementary move then $\tau\sprec \tau'$.

Now if the transition $(h,v)$ to $(h,v')$ occurs twice, then in
particular there are $i<j<k$ such that $(h,v)\in\tau_i$,
$(h,v')\in\tau_j$, and then again $(h,v)\in\tau_k$. But this is a
contradiction to the two properties of $\sprec$.
\end{proof}

\subsection{Projections and lengths}
\label{large link etc}

Lemma 6.2 of Masur-Minsky \cite{masur-minsky:complex2} states the
following, for any hierarchy $H$.

\begin{lemma}{Large Link}
If $Y$ is any essential subsurface in $S$ and 
$$d_Y(\I(H),\T(H)) > M_2$$
then $Y$ is the support of a geodesic $h$ in $H$.

Conversely if  $h\in H$ is any geodesic with $Y=D(h)$, 
$$
\left| |h| - d_Y(\I(H),\T(H)) \right | \le 2M_1.
$$
\end{lemma}
(In this and the next three lemmas, $M_i$ stand for constants
depending only on $S$.)

The proof of this lemma goes through in the infinite setting as
well. The implications of the lemma are that it is possible to
detect, just from the initial and terminal marking, which subsurfaces
in $S$ participate ``strongly'' in the hierarchy, and how long their
supported geodesics are. In fact a slightly more refined fact is shown
(in the course of the proof of this lemma):

\begin{lemma}{Endpoints determined}
IF $h\in H$ is any geodesic with $Y=D(h)$ then
$$
d_Y(\I(h),\I(H)) \le M_1,
$$
$$
d_Y(\T(h),\T(H)) \le M_1,
$$
\end{lemma}
Thus the endpoints of the geodesics are determined up to bounded
error, also. This implies, for annulus geodesics, a statement for {\em
signed} lengths (see (\ref{signed length})): 
\begin{lemma}{annulus signed estimate}
If $Y$ is an annulus supporting a geodesic $h$ in $H$ then
$$
| \tw_Y(\I(H),\T(H)) - [h] | \le M_3
$$
\end{lemma}

\subsection{Elementary moves on clean markings}
\label{del estimates}
Let $\mu$ be a maximal clean marking, which we recall is determined by
a pants decomposition $(u_i)$ and a curve $\bar t_i$ for each
$u_i$ which intersects $u_i$ in a standard way and meets none of
the other $u_j$. A {\em twist move} on $\mu$ changes one $\bar
t_i$ by a Dehn twist, or half-twist, around $u_i$, preserving all
the other curves. 

\realfig{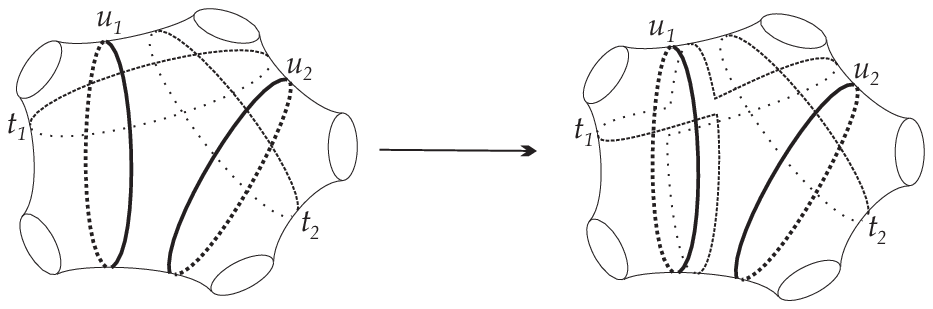}{A twist move in $S_{0,5}$}

A {\em flip move} interchanges a pair $u_i$ and
$\bar t_i$, and adjusts the remaining $\bar t_j$ by a surgery so that
they are disjoint from $\bar t_i$. Figures \ref{twistmove} and
\ref{flipmove} illustrate these moves in $S_{0,5}$. For more details
see \cite{masur-minsky:complex2}.

\realfig{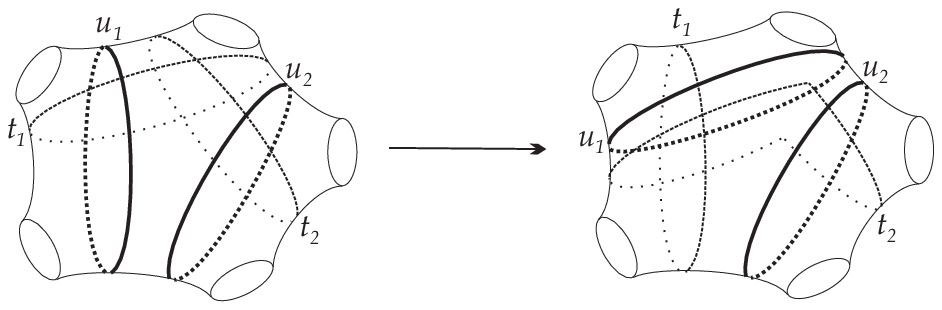}{A flip move in $S_{0,5}$}

Any two markings are related by a sequence
of elementary moves, and we let $d_{el}(\mu,\mu')$ be the length of
the shortest such sequence. 

In \cite{masur-minsky:complex2} we give ways of estimating 
$d_{el}(\mu,\mu')$ using hierarchies and subsurface projections. First
let us introduce some notation.  The relation
$$ x \Qle{a,b} y$$
will denote $x\le ay + b$, and $x\qle y$ denotes the same where $a,b$
are understood as independent of the situation. Similarly let
$$ x \Qeq{a,b} y$$
denote $x\Qle{a,b} y$ and $y\Qle{a,b} x$. 
Let
$$
\Tsh Kx = \begin{cases}
  x & \text{if $x\ge K$} \\
  0 & \text{if $x<K$}
\end{cases}
$$
be the ``threshold function''. 

\medskip

Now, every maximal slice $\tau$ of a hierarchy has an associated
maximal marking $\mu_\tau$. $\mu_\tau$ may not be clean, but there is always a
clean marking $\mu'$ with the same base and satisfying the property
that for any $u$ in $\base(\mu)$ with transversal $t$, 
$u$ has a transversal $t'$ in $\mu'$ such that $d_{u_i}(t_i,t'_i) \le
2$. We say 
that $\mu'$ is {\em compatible with $\tau$}. 
An elementary move $\tau_1\to\tau_2$ on maximal slices yields a bounded number
of elementary moves on compatible clean markings (roughly, a move
advancing a pair $(h,v)$ gives rise to twist moves if $\xi(D(h))=2$, 
a flip move if $\xi(D(h))=4$, and nothing at all if $\xi(D(h))>4$.)

Given two maximal clean markings $\mu,\nu$ there exists (by Theorem
4.6 of \cite{masur-minsky:complex2}) a complete hierarchy $H$ with $\I(H)=\mu$
and $\T(H)=\nu$. 
The length of any resolution  of $H$ is equal to the
sum of lengths 
$$|H| \equiv \sum_{h\in H} |h|,$$ 
and hence, by
considering the sequence of compatible clean markings, we obtain an
upper bound for $d_{el}(\mu,\nu)$. In fact in Theorem 6.10 of
\cite{masur-minsky:complex2} we show that
\begin{equation}\label{H and del}
d_{el}(\mu,\nu) \Qeq{a,b} |H|
\end{equation}
for $a,b$ depending only on $S$.

Using the ideas of Lemma \ref{Large Link} we
further deduce in \cite{masur-minsky:complex2} that
$|H|$ can be estimated in terms of the set of subsurface
projection distances $\{d_W(\mu,\nu)\}$, and in particular: 
\begin{lemma}{dW and del}{\rm (Thm. 6.12 of \cite{masur-minsky:complex2})}
Given $S$ there exists $K_0$ such that, for any $K\ge K_0$ there are
$a,b$ such that, for any pair of maximal clean markings $\mu,\nu$ and
hierarchy $H$ connecting them, 
\begin{equation}
\label{el dist from projections}
d_{el}(\mu,\nu) \Qeq{a,b} \sum_{W\subseteq S} \Tsh K{d_W(\mu,\nu)}
\end{equation}
\end{lemma}
This also follows from the counting results in \S\ref{counting arguments}.

\subsection{Existence of infinite hierarchies}
\label{hierarchies existence}

In Theorem 4.6 of
\cite{masur-minsky:complex2} we showed that 
a hierarchy exists connecting any two finite markings. 
In fact with a closer look we can now obtain

\begin{lemma}{Existence of Hierarchies}
  For any two generalized markings $\I,\T$ of $S$, which do not share
  any infinite-leaf components, there exists a
  hierarchy $H$ with $\I(H)=\I$ and $\T(H)=\T$. 
\end{lemma}

\begin{proof}
The main gap between the finite and infinite existence theorems is
filled by this lemma:

\begin{lemma+}{Existence of geodesics}
Let $X$ be a surface with $\xi(X)\ge 4$.
Let $\mu$ and $\nu$ be two distinct points in $\CC_0(X)\union\EL(X)$.
There exists a tight geodesic $g$ connecting $\mu$ and $\nu$.
\end{lemma+}

In the case $\xi(X)=2$ we need not consider vertices at infinity, and
the existence of a tight geodesic in that case is trivial. 

\begin{proof}[Proof of Lemma \ref{Existence of geodesics}]
The case where $\mu$ and $\nu$ are both vertices of $\CC(X)$ is
covered by \cite[Lemma 4.5]{masur-minsky:complex2}. 
(Of course in this case we need not require them to be distinct.)

In a locally compact $\delta$-hyperbolic space the other cases, which
correspond to endpoints at infinity,  would then follow
from the finite case by a limiting argument. This limiting
step requires a special argument in our setting.

\medskip

Consider first 
the case that $\mu$ is a finite point,
i.e. $\mu\in\CC_0(X)$, and $\nu\in\EL(X)$ is a point at infinity.
By definition there exists a sequence of points $\nu_i\in\CC_0(X)$
such that $\nu_i\to\nu$, and we may form tight geodesics
$[\mu,\nu_i]$ (not necessarily unique). Our goal is to extract a
limiting ray $[\mu,\nu)$.
By $\delta$-hyperbolicity and
the definition of the boundary at infinity, 
we have the following :
\begin{itemize}
\item[(*)] For each $R>0$ there exists $n$ so
that for $i,j\ge n$ the initial segment of length $R$ of
in $[\mu,\nu_i]$ is within $\delta$ of $[\mu,\nu_j]$, 
and vice versa.
\end{itemize}
Thus, let us extend $\mu$ to a maximal clean marking $\I$, and let
$\T_i$ be the marking consisting just of the endpoint
$\nu_i$. By \cite[Thm 4.6]{masur-minsky:complex2}, there exists a
hierarchy $H_i$ with $\I(H_i) = \I$ and $\T(H_i) = \T_i$. In fact we
may assume that the main geodesic of $H_i$ is $[\mu,\nu_i]$.

We now reprise an argument
used in Lemma 6.13 of \cite{masur-minsky:complex2}. Pick
$R>5\delta$, and let 
$n$ be as in (*). For $j\ge n$, let $v'$ be a simplex of $[\mu,\nu_j]$ within
distance $R/2$ of $\mu$, and let $v$ be a simplex of $[\mu,\nu_n]$
such that (via (*)) $d_{\CC(S)}(v,v') \le \delta$.  
Let $\tau$ and $\tau'$ be saturated slices of $H_n$ and $H_j$ whose
bottom simplices are $v$ and $v'$, respectively.
In fact $\tau$ and $\tau'$ must be maximal for large $n$: Since $\I$
is maximal, it intersects any subsurface, and since $\T_j$ or $\T_n$
are far away in $\CC(S)$, they must intersect any subsurface $W\subset
X$ 
for which $[\boundary W]$ is at $\CC_1(X)$-distance 1 from $v$ or $v'$. Hence (using
Theorem \ref{Descent Sequences}) every component domain that arises
in the construction of $\tau$ and $\tau'$ must support a geodesic, and
it follows that $\tau$ and $\tau'$ are maximal. 
Let $m$ and $m'$ be maximal clean markings compatible with $\tau$ and $\tau'$,
respectively; we wish to bound the elementary-move distance
$d_{el}(m,m')$.

Let $J$ be a hierarchy with $\I(J)=m$ and $\T(J)=m'$ ($J$ exists
again  by \cite[Thm 4.6]{masur-minsky:complex2}). If $W$ is a
subsurface that occurs in $J$ then $[\boundary W]$ is within $\delta$
of $v$. We claim that $d_W(\nu_j,\nu_n)$ is uniformly bounded, independently
of $j\ge n$. Let $w$ and $w'$ be points on $[\mu,\nu_n]$ and
$[\mu,\nu_j]$ that are at least $2\delta+2$ further from $\mu$ than $v$
and $v'$, respectively, and such that $d(w,w')\le \delta$ (this is
possible since $R$ is large enough). We can connect $w$ to $w'$ with a
geodesic all of whose vertices are distance at least 2 from
$[\boundary W]$, and hence by the Lipschitz property of $\pi_W$,
we have $d_W(w,w') \le \delta$ (if $W$ is an annulus and $\xi(S)=4$
then the bound is $3\delta$ -- see \S\ref{complexes} for a discussion
of the Lipschitz property.) 
Now $d_W(w,\nu_n)$ and $d_W(w',\nu_j)$
are each bounded by Lemma \ref{Geodesic Projection Bound}, and this
gives us a bound of the form 
\begin{equation}\label{dW for nu}
d_W(\nu_n,\nu_j) = O(1).
\end{equation}
It is a consequence of Lemmas 6.1 and 6.9 of
\cite{masur-minsky:complex2} (see the proof of Lemma 6.7 of that
paper) that given any 
hierarchy $H$ and any subsurface $W\subseteq S$, the projection
$\pi_W(v)$ for any vertex occuring in $H$, if nonempty,
is in an $M$-neighborhood of a geodesic in $\AAA(W)$ connecting
$\pi_W(\I(H))$ and $\pi_W(\T(H))$, where $M=M(S)$.
Applying this to the hierarchy $H_n$ we see that
$\pi_W(m)$ is in an $M$-neighborhood
of a geodesic in $\AAA(W)$ connecting $\pi_W(\mu)$ and $\pi_W(\nu_n)$ (and
similarly using $H_j$, for $\pi_W(m')$, $\pi_W(\mu)$ and $\pi_W(\nu_j)$), 
and we conclude using (\ref{dW for nu}) that 
$$
d_W(m,m') \le d_W(\mu,\nu_n) + O(1).
$$
This gives a uniform bound on the length of the hierarchy $J$ for
any $j\ge n$, and hence on $d_{el}(m,m')$ by (\ref{H and del}).

We conclude that, fixing $\tau$ and $m$, there are only finitely many
possibilities for $m'$, and in particular for $v'$. Thus, there are
only finitely many possibilities for the initial segment of length
$R/2$ of $[\mu,\nu_j]$, for all $j\ge n$. We can therefore take a
subsequence for which these initial segments are all the same, and
then increase $R$ and make the usual diagonalization argument.

The limiting sequence must be tight (since tightness is a local
property) and it must accumulate on $\nu$, since its initial segments
are equal to the initial segments of $[\mu,\nu_j]$ for
all large $j$ (in the subsequence). Hence this is the desired
geodesic $[\mu,\nu)$.

Now consider the case that both $\mu$ and $\nu$ are distinct points in
$\EL(W)$. Let $\mu_i\to \mu$, and $\nu_i\to\nu$ be sequences in
$\CC_0(W)$, and form segments $[\mu_i,\nu_i]$. In this case a similar
condition to (*) holds: There is a constant $\delta'$ and a sequence
of points $x_i\in 
[\mu_i,\nu_i]$  which remain in a bounded subset of $\CC(W)$, such that
\begin{enumerate}
\item[(**)] For any $R>0$ there
exists $n$ such that for all $i\ge n$ the subsegments of
$[\mu_i,\nu_i]$ of radius $R$ centered on $x_i$ are within $\delta'$ of
each other.
\end{enumerate}

The proof is left as an exercise in applying the notion of
$\delta$-thinness of triangles.

Now as before, we take hierarchies $H_i$ with base geodesic
$[\mu_i,\nu_i]$ and apply an argument similar to the previous case to
argue that a subsequence converges to a biinfinite geodesic
$(\mu,\nu)$. In fact this is the exact case treated in Lemma 6.13 of
\cite{masur-minsky:complex2}.

This concludes the proof of Lemma \ref{Existence of geodesics}.
\end{proof}

Continuing with the proof of Lemma \ref{Existence of Hierarchies},
given generalized markings $\I$ and $\T$ with no common infinite-leaf
components, we construct a hierarchy by following the argument of
\cite[Thm 4.6]{masur-minsky:complex2}: 
We begin by constructing a main geodesic $g$ in $\CC(S)$: 
If $\simp(\I)$ is nonempty we choose one of its vertices  to be the
initial vertex of 
the geodesic. If not then $\I$ is a lamination in $\EL(S)$ and gives a
point at infinity. The same holds for $\T$, and we connect the
resulting points with a tight geodesic via Lemma \ref{Existence of
geodesics}. We then inductively build up a sequence of ``partial
hierarchies'' $H_n$. At each stage
we find an ``unutilized
configuration'' in $H_n$, which is a triple $(W,b,f)$ where 
$W$ is a component domain in some
geodesic in $H_n$, $b,f\in H_n$ are tight geodesics such that 
$b\bsubd W\fsubd f$, but $W$ is not
the support of a geodesic in $H_n$. For such a domain 
we use Lemma \ref{Existence of geodesics} to get a 
tight geodesic $h$ with $\I(h)=\I(W,b)$ and $\T(h)=\T(W,f)$, and let
$H_{n+1} = H_n \union \{h\}$.

In the finite case, this was sufficient: the process was guaranteed to
terminate after finitely many steps and the last $H_n$ is a hierarchy,
because there are no more domains to fill in. 
In the general case, because of the possibility of infinite geodesics, the
process of filling in unutilized configurations may be infinite. Thus we must 
order the process in such a way that the {\em union} of the $H_n$ is a
hierarchy, i.e. so that every unutilized configuration is filled in a finite
number of steps.

To do this, we maintain an {\em order} on each  $H_n$, 
and in addition fix a choice of
basepoint $v_h$ for each geodesic
$h$. At the first step $H_1$ is just the main geodesic with an
arbitrary choice of basepoint. At each step of the process, 
consider the first geodesic $h$ in the order. If there are any 
unutilized configurations $(W,b,f)$ with $W$ a component domain of $h$,
choose one minimizing the distance of $[\boundary W]$ to $v_h$ in
$\CC(D(h))$ (there are finitely many such).
Construct a new tight geodesic $k$ with $D(k)=W$ and $b\bsubd k \fsubd f$, 
pick an arbitrary basepoint $v_k$, and add $k$ to
$H_n$ obtaining $H_{n+1}$, with the ordering unchanged except that $h$
and $k$ should now be the last two elements. 
If there are no such unutilized configurations, adjust the order so
that $h$ becomes the last element.

Repeating this, we see that every time a
geodesic is created or examined, it will be examined again in a finite
time. The component domains are examined by order of distance from the
basepoints, and hence every unutilized configuration will be filled after a
finite number of steps. Thus the union of the $H_n$
will have no unutilized configurations, and must therefore be a hierarchy.
\end{proof}

\subsection{Vertices, edges and 3-holed spheres}
\label{more hierarchies}

We can now  consider in a little more detail how vertices appear in a  
hierarchy. 
An edge $e=[vw]$ in a geodesic $h\in H$ with $\xi(D(h))=4$ is called 
a {\em 4-edge}. We write $v=e^-$ and $w=e^+$, where $v<w$ in the
natural order of $h$
(note in this case both $v$ and $w$ are actually
vertices, not general simplices).

We say that a vertex $v$ ``appears in $H$''  if it is part of a 
simplex in some geodesic in $H$. 
If $v$ is not a vertex of $\simp(\I(H))$ or $\simp(\T(H))$ then it is
called {\em internal}. A vertex of $\simp(\I(H))$ which does not have a
transversal in $\I(H)$ is called {\em parabolic} in $\I(H)$, and
similarly for $\T(H)$.

\begin{lemma}{Vertex configurations}
Let $v$ be a vertex appearing in $H$. 
Then $\gamma_v$ intersects $\base(\I(H))$ if and only if
there exists a 4-edge $e_1$ with $v=e_1^-$. Similarly
$\gamma_v$ intersects $\base(\T(H))$ if and only if
there exists a 4-edge $e_2$ with $v=e_2^+$. 

The edges $e_1$ and $e_2$ are unique.
\end{lemma}

\begin{proof}
Suppose $\gamma_v$ intersects $\base(\T(H))$.  We must find a geodesic $h$
with $\xi(D(h))=4$, such that $v$ is a vertex in $h$ but not the
last, and we must show that $h$ is unique.

The condition that $v$ appears in $H$ is equivalent
to the condition that the annulus $A=\collar(\gamma_v)$ is a component
domain in 
$H$, and the condition that $v$ intersects $\base(\T(H))$ means that
$\Sigma^+(A)$ is nonempty, since it contains the main geodesic $g_H$.
By part (3) of Theorem \ref{Descent Sequences}, 
there exists a geodesic $f$ such that $A\fsub f$, so by definition
there must be some $h\in H$ with $A\fsubd h$. 

The condition $A\fsubd h$ means that $v$ appears as a vertex in a simplex $w$
of $h$,
and that $\vsucc(w)|_A\ne \emptyset$. If
$w$ is not the last 
simplex then, if $\xi(D(h))>4$, 
$\vsucc(w)$ would be a simplex disjoint from $w$ and hence from $v$, a
contradiction. 
Thus $\xi(D(h)) = 4$, $v=w$, and $\vsucc(w)$ is a vertex $w'$.
This gives us our 4-edge, $e_2 = [vw']$.

We claim that 
$w$ cannot be the last simplex of $h$. If, by way of contradiction, it
is, then $\vsucc(w)=\T(h)$ and $v\in\base(\T(h))$,
so $\T(h)$ must contain a
transversal to $v$. If $h=g_H$ then $\gamma_v$ is a base curve of $\T(H)$,
contradicting the hypothesis of the lemma. 
Hence $h\fsubd h'$ for some $h'$. Now since $\T(h)$ contains
a transversal, $D(h)$ must be a component domain for the last simplex
of $h'$, and $\T(h) = \T(h')|_{D(h)}$. Hence
$v$ must be a component of $\base(\T(h'))$ with a transversal,
and we repeat this by induction until we terminate with the main
geodesic $g_H$, and obtain our contradiction.

Uniqueness of $e_2$ is seen as follows. If there were two such
edges, there would be two geodesics $h,h'$ with $\xi=4$, and such
that $A\fsubd h$, $A\fsubd h'$. Hence both are in $\Sigma^+(A)$, but
this contradicts part (1) of Theorem \ref{Descent Sequences}, which
says that  $\Sigma^+(A)$ is a sequence
$f_0\fsubd f_1 \fsubd \cdots$.

In the converse direction, if a 4-edge $e$ exists with $v=e^-$
then $v$ appears in a geodesic $h_1$ with $\xi=4$, where $v$ is not the
last vertex. 
Thus $\gamma_v$ intersects the last vertex of $h_1$, and in
particular it intersects $\base(\T(h_1))$. 
Note that $A\fsubd h_1$, and Theorem 
\ref{Descent Sequences} gives us a sequence
$h_1\fsubd h_2 \fsubd \cdots \fsubd
g_H$. We claim by induction that $\gamma_v$ intersects
$\base(\T(h_i))$ for each $i$. Let $w_i =
\max\phi_{h_i}(D(h_{i-1}))$, so that $D(h_{i-1})$ is a 
component domain of $(D(h_i),w_i)$. 
By Lemma \ref{max footprint}, $w_i$ must equal 
$\max\phi_{h_i}(A)$.
If $w_i$ is not the last simplex then 
the last simplex of $h_i$ intersects $\gamma_v$, and hence does
$\base(\T(h_i))$. If $w_i$ is the last 
simplex, then $\T(h_{i-1}) = \T(h_i)|_{D(h_{i-1})}$, so that the
statement for $h_i$ follows from the statement for $h_{i-1}$.
Hence $\gamma_v$ intersects $\base(\T(H))$, which is what we wanted to show.

\medskip

The argument for $e_1$ is identical, with directions reversed. 
\end{proof}

Now given $v$ appearing in a hierarchy $H$, and fixing a resolution
$\{\tau_i\}_{i\in\II}$ of $H$, with $\II$ a subinterval of $\Z$, let
$$
J(v) = \{ i\in \II:  v \in \base(\mu_{\tau_i})\}.
$$
The following fact will be fundamental for us:
\begin{lemma}{vertex interval}
Let $H$ be a $4$-complete hierarchy  and 
$\{\tau_i\}_{i\in\II}$ a resolution. 
If $v$ is a vertex in $H$, then
$J(v)$ is an interval in $\Z$.
\end{lemma}
Recall that an interval can be finite, one-sided infinite, or bi-infinite.

\begin{proof}
  The only type of elementary move that can change the vertex set of a
  slice is a move along a 4-edge $e$, which replaces one vertex
  $e^-$ by its successor $e^+$. By Lemma \ref{Resolution sweep},
  there is exactly one step in the resolution which advances along any
  given edge. 
  By Lemma \ref{Vertex configurations}, there is at most one $4$-edge,
  $e_1$, for which $e_1^+=v$, 
  and at most one $4$-edge, $e_2$, for which $e_2^- = v$. Thus there is at
  most one move $\tau_{i-1}\to\tau_{i}$ at which $v$ appears, and at
  most one   move $\tau_j\to\tau_{j+1}$ at which $v$ disappears. 
  The possibilities for $J(v)$, depending on the existence of $e_1$
  and $e_2$, are therefore 
  $[i,\sup \II]$, $[\inf\II,j]$, $\II$, $[i,j]$ (if $i\le j$),
 or   $[\inf\II,j]\union [i,\sup\II]$ (if $j<i$).  Only the last of
  these fails to be an interval and must be ruled out (note that
  $\inf\II$ may be $-\infty$ and $\sup\II$ may be $\infty$, yielding
  infinite intervals).

 In the last case where  $j<i$, $v$ must be in all slices after
 $\tau_i$. We claim that $v$ is therefore in $\simp(\T(H))$.
Let $(g_H,w_k)$ be the bottom pair of $\tau_k$. For $k>i$ we have
$d_{\CC_1(S)}(w_k,v)\le 1$,  and on the other hand 
Lemma \ref{Resolution sweep} implies that all vertices of $g_H$ are
 obtained in the resolution --  hence $g_H$ is finite in the
 forward direction, and for sufficiently large $k$ we must have $w_k$
 equal to the last vertex of $g_H$, which is in  $\simp(\T(g_H)) =
 \simp(\T(H))$. If $v\ne w_k$, then it is contained in a component $W$ of
 $S\setminus \collar(w_k)$, which must therefore support a geodesic
 $g'$ (otherwise the slice $\tau_k$ cannot have any vertices in $W$),
and we have $\T(g') = \T(g_H)|_W$. The same argument implies that $g'$ is
 finite in the forward direction, and we proceed inductively. The
 process must terminate, and at that point $v$ must be a vertex of
 $\simp(\T(H))$. 

 However Lemma \ref{Vertex configurations} then implies that $v$ never
 appears as $e^-$ for a $4$-edge $e$, and this is a
 contradiction. This rules out the case $j<i$. 
\end{proof}

\subsection*{Three-holed spheres}
\label{xithree}
The following result for three-holed spheres is analogous to 
Lemma \ref{Vertex configurations}.

\begin{lemma}{xithree configurations}
Let $Y$ a component domain of a geodesic in $H$, and suppose
$\xi(Y)=3$. If $\T(H)|_Y \ne \emptyset$ then there exists a unique geodesic
$f\in H$ with $\xi(D(f))=4$ and $Y\fsubd f$. Similarly if
$\I(H)|_Y \ne \emptyset$ then there exists a unique geodesic
$b\in H$ with $\xi(D(b))=4$ and $b\bsubd Y$.
\end{lemma}

Note that $\T(H)|_Y\ne\emptyset$ just means that $Y$ has an essential
intersection with a base curve of $\T(H)$. In particular if $Y$
intersects both $\base(\T(H))$ and $\base(\I(H))$ then the lemma gives
$b,f$ with $\xi=4$ and $b\bsubd Y \fsubd f$. We call this a {\em
gluing configuration}, and it will be used in the model manifold
construction in Section \ref{model}.

\begin{proof}
If $\T(H)|_Y \ne \emptyset$, part (3) of Theorem \ref{Descent Sequences}
gives some $f$ with $Y\fsubd f$. Thus it only remains to
verify that $\xi(D(f))=4$. 

The condition $Y\fsubd f$ means that
there is a simplex $w$ of $f$, with $Y$ a component domain of
$(D(f),w)$, and such that $\vsucc(w)|_Y\ne \emptyset$. 
Suppose that $w$ is not the last vertex of $f$. 
If $\xi(D(f))>4$ then $\vsucc(w)$ is disjoint from $w$, and since
a three-holed sphere can have no nontrivial nonperipheral simple
curves, $\vsucc(w)|_Y =  \emptyset$. This contradiction implies
$\xi(D(f))=4$. 

If $w$ is the last vertex then $w\subset \base(\T(f))$ and
$\vsucc(w) = \T(f)$. Since $Y$ is a three holed sphere bounded
by $w$ in $D(f)$,  it cannot
contain other components of $\base(\T(f))$, so we obtain
$\vsucc(w)|_Y = \emptyset$, again a contradiction.

Uniqueness of $f$ follows from the sequential structure of
$\Sigma^+(Y)$ (part (1) of Theorem \ref{Descent Sequences}).
The backwards case is proved in the same way.
\end{proof}

\section{The coarse projection property}
\label{coarse projection}

In this section we will begin the geometric argument that connects the
geometry of a Kleinian surface group $\rho:\pi_1(S)\to\PSL 2(\C)$
with the structure of the complex of curves $\CC(S)$.
The main tool for this is the 
``short-curve projection'' $\Pi_\rho$, which takes any element
$a\in\AAA(S)$ to the set of short curves on pleated surfaces in $N$
mapping the curves associated to $a$ geodesically. This projection was
studied in 
\cite{minsky:kgcc}  and \cite{minsky:boundgeom}, where it was shown to
have properties analogous to a nearest-point projection to a convex
set. Here we will refine these to obtain {\em relative} statements
about projections $\pi_Y$ to subsurfaces. In Section \ref{projection
bounds}, we will use $\Pi_\rho$ to generate our main geometric
estimates. 

We proceed to define the projection map.
Fix $L\ge L_0$, where $L_0$ is Bers' constant (see \S\ref{thick thin}).
As in \cite{minsky:boundgeom}, we define for a hyperbolic structure
$\sigma$ on $int(S)$
$$
\short_L(\sigma) = \{\alpha\in\CC_0(S): \ell_\sigma(\alpha)\le L\}.
$$

For any lamination $\lambda\in\GL(S)$, we define
$$
\Pi_{\rho,L}(\lambda) = \bigcup_{f\in\pleat_\rho(\lambda)} \short_L(\sigma_f).
$$
(If $L$ is understood we write just $\Pi_\rho$).

In particular any $x\in\AAA(S)$ determines a lamination $\lambda_x$
consisting of leaves representing the vertices of the smallest simplex
containing $x$, and we define $\Pi_{\rho,L}(x) \equiv
\Pi_{\rho,L}(\lambda_x)$. This gives us a map
 $$
 \Pi_{\rho,L} : \AAA(S) \to \PP(\CC(\rho,L))
 $$
where $\PP(X)$ denotes the set of subsets of $X$. 

We will establish the following ``coarse projection'' property, generalizing
Lemma 3.2 of \cite{minsky:boundgeom}:

\begin{theorem+}{Relative Coarse Projection}
For a surface $S$ and a constant $L\ge L_0$ there exist $D_0,D_1>0$ such
that,  if $\rho\in\DD(S)$ and $Y\subseteq S$ is an essential
subsurface with $\xi(Y)\ne 3$, then:

\begin{enumerate}
\item[(P1)] (Coarse definition)
If $v$ is any vertex in $\AAA(S)$ with $\pi_Y(v)\ne \emptyset$, then
\begin{equation*}
\diam_Y(\Pi_{\rho,L}(v)) \le D_0.
\end{equation*}

\item[(P2)] (Relative Coarse Lipschitz)
If $v,w$ are adjacent vertices in $\AAA(S)$ with 
$\pi_Y(v)\ne \emptyset$ and
$\pi_Y(w)\ne \emptyset$, then
\begin{equation*}
\diam_Y(\Pi_{\rho,L}(v)\union\Pi_{\rho,L}(w)) \le D_1.
\end{equation*}

\item[(P3)] (Relative Coarse Idempotence)
If $v\in\CC(\rho,L)$ and $\pi_Y(v)\ne \emptyset$, then
$$
d_Y(v,\Pi_{\rho,L}(v)) = 0.
$$
\end{enumerate}
\end{theorem+}

Note that property (P2) translates to a weak sort of Lipschitz
condition for the map $\pi_Y\circ\Pi_{\rho,L}$: If
$v_0,v_1,\ldots,v_n$ is a sequence of vertices in $\AAA(S)$ 
with $v_i,v_{i+1}$ adjacent
such that $\pi_Y(v_i)\ne\emptyset$, and in addition
$\pi_Y(\Pi_{\rho,L}(v_i))\ne \emptyset$, then repeated application of
(P2) yields
$$
\diam_Y(\Pi_{\rho,L}(v_0)\union\Pi_{\rho,L}(v_n)) \le D_1 n.
$$
The issue of when $\Pi_{\rho,L}$ is nonempty will be addressed in
Lemma \ref{nonempty Pi}.

\subsection{Train tracks}
\label{tracks}

Before proceeding with the proof of Theorem \ref{Relative Coarse
Projection}, let us digress a bit to discuss the construction and
properties of train tracks. 
Thurston \cite{wpt:notes}
first introduced the notion of a train track in a surface as
a finite approximation to a geodesic lamination. 
Although it has by now become a
standard tool in the field, we will attempt to be careful with details
here, especially because of the presence of thin parts. 

\subsubsection*{Definitions}
(See Penner-Harer \cite{penner-harer}).
A {\em train track} in a surface $S$ is an embedded 1-complex $\tau$ 
with a special structure at the vertices. It is convenient to describe
this in terms of a foliation of a small neighborhood of $\tau$ by
leaves (called ``ties'') that are transverse to $\tau$ at every
point. A tie passing through a vertex locally divides the ends of
adjacent edges according to which side of the tie they are on. One
thinks of these as ``incoming'' and ``outgoing'' sides, and both sets
must be nonempty. Vertices of a train track are called ``switches''
and edges are called ``branches.''

Two branch ends coming in to the same side of a switch bound between
them a corner of a complementary region which we call a ``cusp''. 
If a switch meets $n$ branch ends then a small neighborhood of the
switch is cut by $\tau$ into $n$ corners, $n-2$ of which are
cusps. The other two are called smooth.
A further condition usually imposed on the train track is that
each complementary region is {\em hyperbolic}, in the sense that it is
not a disk with 2 or fewer boundary cusps, a once punctured
disk with no boundary cusps, or an annulus with no boundary cusps.  

When $\boundary S\ne \emptyset$,
we will also allow train tracks to have branches which terminate in 
the boundary. A regular neighborhood of a boundary component is cut by
$\tau$ into regions which we think of as cusps
for the purpose of the hyperbolicity condition.

A {\em train route} in $\tau$ is an immersion of a 1-manifold into
$\tau$ which always traverses
switches from one side to the other, and whose endpoints (if any) map
to $\boundary S$. 
An element of $\AAA(S)$ is
is {\em carried} by $\tau$ if it can be represented by a 
train route. 
If $\alpha$ is carried in $\tau$ it imposes
a {\em measure} on its branches counting how many times each is
traversed by $\alpha$. For each branch $b$ we denote this measure by
$\alpha(b)$. We let 
$$
\ell_\tau(\alpha) = \sum_b \alpha(b)
$$
denote the {\em combinatorial length} of $\alpha$.

A train track in $int(S)$ is just the restriction of a train track in
$S$ to the interior. In particular branches that terminate in
$\boundary S$ become branches that exit cusps of $int(S)$. 

\subsubsection*{Collapsing curves to train tracks}

Fix a Margulis constant $\bar\ep>0$, and let $\sigma$ be a complete
(finite-area) hyperbolic
metric on $int(S)$.  If $\gamma$ is a geodesic representative of a vertex
in $\AAA(S)$, or in general any geodesic lamination, and
$0<\ep<\bar\ep$, we define an
{\em $\ep$-collapse} of $\gamma$ to a train track to be a map
$$q:int(S)\to  S'$$
that is homotopic to a homeomorphism, and such that:
\begin{enumerate}
\item $q(\gamma)$ is a train track $\tau$ in $S'$, and 
$q$ restricted to any leaf of $\gamma$ is a train route.
\item There is a metric on $\tau$ such that $q|_\gamma$ is a local
  isometry.
\item For each $x\in\tau$ the preimage $q^{-1}(x)$ is a point or arc.
If $q^{-1}(x)$ touches  $S_{thick(\bar\ep)}$ then its length
is at most $\ep$.
\end{enumerate}

\begin{lemma}{collapse exists}
Given a surface $S$ and $\bar\ep>0$, for each $\ep<\bar\ep$, 
for any hyperbolic metric $\sigma$ on $S$ and any
geodesic lamination $\gamma$ there is an $\ep$-collapse of $\gamma$ to
a train track. 
The total length of 
$$\tau_{\bar\ep} \equiv q(\gamma\intersect S_{thick(\ep)})$$
is bounded by a constant $K$ depending only on $\ep$. 
\end{lemma}

\begin{proof}
A version of this construction is discussed in Thurston
\cite[\S8.9]{wpt:notes}. 
See Brock \cite{brock:continuity} for a complete
discussion in the case without thin parts.

First, for simplicity let us add enough leaves to
$\gamma$ to obtain a lamination $\lambda$ whose complement 
is a union of (interiors of) ideal hyperbolic triangles. 
Fixing a small $\ep$, let $\FF_\ep$ be the foliation of the ends of
each of these triangles by horocyclic arcs of length less than $\ep$. 
Thus $\FF_\ep$ is supported in an open subset of $S\setminus\lambda$
we'll call $U_\ep$.
Because the tangent directions of the horocycles form a Lipschitz
line field, $\FF_\ep$ may be extended across the leaves of $\lambda$
to be a foliation of the interior of the closure $\overline U_\ep$. Let us
continue to call this foliation $\FF_\ep$. In fact $\FF_\ep$ extends
to give a decomposition of all of $\overline U_\ep$ into leaves, where
the boundary consists of endpoints of leaves or boundary horocycles of
length $\ep$.  

This foliation inherits a transverse measure from the length measure
along leaves of $\lambda$. This is evident from the fact that, in each
end of a triangle, the horocyclic flow preserves length for its family
of orthogonal geodesics. 

Let us establish the following claims:

\begin{enumerate}
\item
Each cusp of $S$ has a
neighborhood that, if it intersects $\FF_\ep$ at all, is foliated by
closed leaves of $\FF_\ep$.
\item
All leaves of $\FF_\ep$ are compact. 
\item There are constants $C,C'$ such that, 
if $\ep<C\bar\ep$, then any
leaf meeting the $\bar\ep$-thick part of $S$ is an arc, and has length
bounded by $C'\ep$. 
\end{enumerate}

\begin{proof}
Claim (1)
is fairly clear: There is a uniform $\delta$
such that the $\delta$-Margulis tube associated to a cusp either avoids
$\lambda$ altogether or contains leaves of $\lambda$ which go all the
way out the cusp. In the latter case, if we choose $\delta<\ep/2$,
say, the leaves of $\FF_\ep$ in this neighborhood are complete
horocycles going around the cusp.

Let a {\em foliated rectangle} of $\FF_\ep$ denote a 
region on which $\FF_\ep$ is equivalent to the foliation of $[0,1]^2$
by vertical arcs, and for which the horizontal arcs are geodesics
orthogonal to $\FF_\ep$ (e.g. leaves of $\lambda$).
The leaves of $\FF_\ep$ identified with $\{0\}\times[0,1]$ and
$\{1\}\times[0,1]$ are its {\em boundary leaves}, and 
the distance along $\lambda$ between them is called the {\em width} 
$w(R)$. The lengths of the leaves vary by a factor of at
most $e^w$. Note that the area of $R$ is at least the smaller leaf
length times the width. 

We can make the same definitions for the lift $\til\FF_\ep$ of $\FF_\ep$ to the
universal cover $\Hyp^2$. 
Let $R$ be a foliated rectangle in the univeral cover, with boundary
leaves $l_0,l_1$.
Suppose that $l_0$ and $l_1$ can be extended indefinitely in one 
direction (say the upward  vertical direction in the identification
with $[0,1]^2$). There is an $a$ depending only on $\FF_\ep$ such that,
provided $w(R)<a$, the whole rectangle can be extended indefinitely in
the upward direction -- that is the region between the extensions of
$l_0$ and $l_1$ can be foliated to make a foliated rectangle in
$\Hyp^2$ for
arbitrarily long extensions.  To see this, choose $a$ less than the
minimal length of an arc of $\lambda$ in the boundary of a
complementary region of $\til \FF_\ep$. Suppose $R$ has been extended to a
rectangle $R'$ with boundary leaves $l'_0$ and $l'_1$. Since
$w(R')=w(R)<a$, and by assumption $l'_0$ and $l'_1$ can be extended
further, the top boundary of $R'$ does not lie on the boundary of a
complementary region, and hence the rectangle can be extended
further. 

This gives at least an immersed foliated rectangle. However this must
in fact be an embedding (in $\Hyp^2$), since a self-intersection would give
rise to a disk whose boundary consists of one arc of $\til \FF_\ep$ and
one arc transverse to $\til\FF_\ep$ -- violating the index formula for
line fields since $\FF_\ep$ has an extension to a foliation of $S$
with negative-index singularities. 

Now suppose that $l$ is a noncompact leaf of $\FF_\ep$. 
In view of (1), $l$ must accumulate somewhere in $S$. Hence there must
be two segments $l_0$ and $l_1$ close enough together to bound a
foliated rectangle $R$, and such that $l_0$ and $l_1$ extend
infinitely in one direction. In the universal cover, the lift of $R$
must extend indefinitely by the above argument. The extension
$R_\infty\homeo [0,1]\times[0,\infty)$
has infinite area, so its immersed image downstairs cannot be
embedded. It follows that it is an annulus in which the half-leaf
$\{0\}\times[0,\infty)$ (which is part of $l$) either spirals or maps
to a closed leaf. In the former case the annulus has infinite area,
and in the latter case $l$ is a closed leaf after all. We conclude
that there are no noncompact leaves.  This proves claim (2).

For claim (3),
if $x$ is in the $\bar\ep$ thick part it is the center of an embedded
$\bar\ep/2$-disk. Consider the intersection $m$ of 
a leaf of $\FF_\ep$ passing through
$x$ with a ball of radius
$\bar\ep/4$  around $x$. Every interval of $m\setminus \lambda$
is the boundary leaf of a foliated rectangle on one side or the other,
of width at least $\bar\ep/4$, these rectangles are all disjoint, and
hence $m$
is in the boundary of a region of $\FF_\ep$ with area at least
$C|m|$, for a uniform $C$ (see Brock \cite{brock:continuity} for
details). On the other hand the total area of $\FF_\ep$ is at most
$N\ep$, where $N$ is the number of vertices of ideal triangles in the
lamination, by an easy computation (and using the fact that
$area(\lambda)=0$ -- see \cite{casson-bleiler}). It follows that 
$|m|$ is bounded by $C'\ep$ where $C'$ depends
on $\bar\ep$. In particular if $\ep$ is sufficiently small then the
leaf containing $m$ must
terminate before reaching the boundary of the $\bar\ep/4$ disk,
and claim (3) follows.

\end{proof}

In order to get rid of the closed leaves of $\FF_\ep$, we make the
following adjustment. If $A$ is a cusp neighborhood foliated by closed
leaves, then $A$ is cut by leaves of $\lambda$ into one or more
cusp regions, each equal to an end of one of the complementary
triangles of $\lambda$. Pick one of these regions and erase all the
arcs of $\FF_\ep$ inside the corresponding end. This 
will cut all the loops in $A$ into arcs. 

If $A$ is the collar of a geodesic with length less than $\bar\ep$,
and contains closed leaves of $\FF_\ep$, then there are leaves of
$\lambda$ passing  all the way through $A$.
Choose one complementary region $R$ of $A\setminus\lambda$ -- it is
contained in 
the cusp end of a triangle in $S\setminus \lambda$. For that end
choose $\ep'$ sufficiently small that, if we remove all the $\FF_\ep$
arcs of length at least $\ep'$ in the cusp end containing $R$ then
this will remove all the arcs in 
$R$, again cutting all the closed leaves. 

By claim (3) there are no other closed leaves in $\FF_\ep$.

\medskip 

Call the new foliation $\FF'_\ep$. 
Now consider the quotient  $q:S\to S' = S/\FF'_\ep$ obtained by
identifying each leaf of $\FF'_\ep$ to a point. We claim this is our
desired map, if $\ep$ is taken sufficiently small. Since all leaves of
$\FF'_\ep$ are arcs, one can see by Moore's theorem
\cite{moore:collapse,moore:planetop} that $S'$ is a surface
homeomorphic to $S$ and that $q$ is homotopic to a homeomorphism.
In fact if $x$ is a point on a leaf of $\FF'_\ep$,
we can explicitly describe a neighborhood of $q(x)$ by considering a
neighborhood of the leaves passing near $x$ (Figure \ref{localfol}). 

\realfig{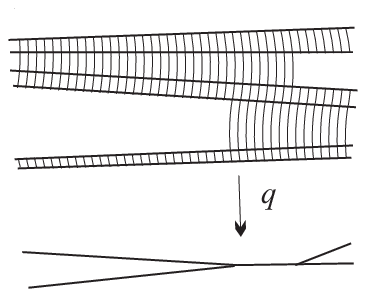}{The local picture of
$\FF'_\ep$, collapsing to a train track.}

The image $\tau=q(\lambda)$ in $S'$ is a finite graph (with some edges
going out the punctures).
To see the train track structure on $\tau$, extend the foliation
$\FF'_\ep$ slightly to obtain, in the image, a foliation in a
neighborhood of $\tau$ transverse to $\tau$. 
An interior leaf $a$ of
$\FF'_\ep$ maps to an interior point of $\tau$. This is because $a$
terminates on two leaf segments of $\lambda$, and leaves sufficiently
close to $a$ terminate on the same segments; hence producing a
foliated rectangle that maps to a segment in $\tau$. 
A boundary leaf $b$ of $\FF_\ep$ maps to a switch of $\tau$: The leaves on
either side of $b$ map to arcs of $\tau$ on either side of $q(b)$,
hence producing the switch structure. 

The complementary regions of $\tau$ are all three-cusped disks, obtained
directly from the triangular components of $S\setminus\lambda$
by the quotient map.

A branch $b$ of $\tau$ is the $q$-image of a foliated rectangle $R_b$
in $S$. The lengths of $\lambda$-leaves in $R_b$ induce a metric on
$b$ with total length $w(R_b)$.
Two leaves of $\lambda$
in $R_b$ lie on the boundary, and are also boundary arcs of triangles in
$S\setminus  \lambda$.
The intersection of these boundary arcs with $S_{thick(\ep)}$
has total length bounded by some $K(\ep)$ (proportional to $\log 1/\ep$).
This gives a bound on the total length of $\tau_{\bar\ep}$.

Now any point $x$ of $\tau$ has preimage $q^{-1}(x)$ which is an arc
of $\FF'_\ep$, with both endpoints on $\lambda$. If $q^{-1}(x)$ meets
the $\bar\ep$-thick part of $S$ then its length is bounded by $C'\ep$,
as  in the proof of claim (3).

Thus to obtain an $\ep$-collapse, we take the quotient
$S/\FF_{\ep/C'}$,  and then restrict to the sub-track 
$\tau_\gamma = q(\gamma)$, noting that all the properties of a
collapse are inherited by this track. 
\end{proof}

\subsubsection*{Nice representatives}
If $\tau$ is obtained from $\gamma$ by an $\ep$-collapse $q$, and
$\beta$ is a curve carried on $\tau$ (more accurately, $q(\beta)$ is
carried on $\tau$), then $\beta$ can be realized after homotopy as a
chain of segments
$$
\beta^\tau = \beta_1 * t_1 * \cdots * \beta_n * t_n
$$
where $n=\ell_\tau(\beta)$, $\beta_i$ are arcs of $\gamma$ and each
$t_i$ is a subarc of 
$q^{-1}(s)$ for a switch $s$. This is done by first choosing, for each
branch $b$ of $\tau$, an arc of $\gamma$ in $q^{-1}(b)$ which maps
isometrically to $b$, and then for branches meeting at a switch,
connecting their preimage arcs with the preimage of the switch.
The total length of $\beta^\tau$ in the $\bar\ep$-thick part of $S$ is
bounded by
\begin{equation}\label{thick beta bound}
\ell(\beta^\tau\intersect S_{thick(\bar\ep)}) \le (K(\ep)+\ep)
\ell_\tau(\beta)
\end{equation}
since each $t_i$ has length at most $\ep$ and each $\beta_i$ has
length at most $K(\ep)$ in the thick part.

\subsubsection*{Short curves with prescribed intersection properties}

This lemma will be used in the final stages of the proof in
\S\ref{coarse proof}:
\begin{lemma}{filling track}
For any compact surface $S$ 
there is a number $A$  such that, for any train track
$\tau$ in $S$, if 
$\beta\in\CC_0(S)$ has essential intersection with
some $\gamma\in\AAA_0(S)$ carried on $\tau$, then there exists
$\alpha\in\AAA_0(S)$ carried on $\tau$ with
$$\ell_\tau(\alpha) \le A,$$ that also intersects
$\beta$ essentially.
\end{lemma}

\begin{proof}
Since, up to homeomorphisms of $S$, there are only finitely many train
  tracks in $S$ (see e.g. Penner-Harer \cite{penner-harer}), and since
  for any finite set of 
  curves or arcs 
  carried in $\tau$ there is an upper bound on $\ell_\tau$, it will
  suffice to prove the following statement: For any train track $\tau$
  there is a finite set $\{\alpha_1,\ldots,\alpha_k\}\subset\AAA_0(S)$
  carried in $\tau$ so that, for any $\beta\in\CC_0(S)$, if
  $i(\beta,\alpha_j) = 0$ for each $j$ then $i(\beta,\gamma)=0$ for
  every $\gamma\in\AAA_0(S)$ carried on $\tau$.

  For any set $X$ in $\AAA_0(S)$ there is an essential subsurface
  $S_X$ {\em filled } by $X$, which is unique up to isotopy. To form
  $S_X$ place representatives of $X$ in minimal position (for example
  using geodesic representatives in a hyperbolic metric on $S$), form
  a regular neighborhood, and adjoin all complementary components which
  are disks or annuli with one boundary component in $\boundary S$. 
  If $\beta\in\CC_0(S)$ intersects $S_X$ essentially then $\beta$ must
  intersect   one of the elements of $X$ essentially. 

Let $X(\tau)$ be the set of elements of $\AAA_0(S)$ carried on $\tau$,
  and let $X_1\subset X_2 \subset \cdots $ be an exhaustion of
  $X(\tau)$ by finite sets.  An Euler characteristic argument implies
  that the sequence $S_{X_1}\subseteq S_{X_2}\subseteq\cdots$ is
  eventually constant (up to isotopy), and hence there exists some $k$
  for which $S_{X_k}=S_{X(\tau)}$.
  
  $X_k$ is therefore our desired finite set $\{\alpha_i\}$.
\end{proof}
(Remark: there ought to be a more concrete proof of this lemma, which
gives estimates on $A$, but I have not found it).

\subsection{Proof of Coarse Projection}
\label{coarse proof}

We can now proceed with the proof of Theorem \ref{Relative Coarse
Projection}. 

\subsubsection*{Property (P3):} This property is immediate: $v$ is
realized with its minimal $\rho$-length 
in any $f\in \pleat_\rho(v)$, and since this is at most $L$ we have
$v\in\short_L(\sigma_f)\subset \Pi_\rho(v)$. Since
$\pi_Y(v)$ is not empty, it is
also a component of $\pi_Y(\Pi_\rho(v))$ which establishes the
property. 

\subsubsection*{Property (P1) $\implies$ Property (P2):}
Suppose first that $\xi(Y)\ge 4$. Then $Y$ has essential intersection
with any pants decomposition. Since $L\ge L_0$, $\Pi_\rho(x)$ contains
a pants decomposition, and hence $\pi_Y(\Pi_\rho(x))$ is nonempty
for any lamination $x$. Now if 
$v$ and $w$ are adjacent in $\AAA_0(S)$, the disjoint union of their
representatives is a lamination which we denote $v\union w$. We have
$$\pi_Y(\Pi_{\rho,L}(v\union w)) \subset \pi_Y(\Pi_{\rho,L}(v))\intersect
\pi_Y(\Pi_{\rho,L} (w)),$$
and in particular this intersection is non-empty. 
Thus Property (P1) for $v$ and $w$
immediately yields (P2), with $D_1=2D_0$.

Now suppose $\xi(Y)=2$, so $Y$ is an annulus. Let $u\in\CC_0(S)$
represent the core curve of $Y$. If either $\pi_Y(\Pi_{\rho,L}(v))$ or
$\pi_Y(\Pi_{\rho,L}(w))$ are empty there is nothing to prove; so
assume they are nonempty.
In particular for some $f\in\pleat_\rho(v)$, $\gamma_u$ is
intersected by a curve of $\sigma_f$-length at most $L$, and hence
$\ell_{\sigma_f}(u) \ge \neck(L)$. 
It follows that $f(S)$ must stay out of the tube
$\MT_{\ep'}(u)$, where $\ep'=\epshrink(\neck(L))$ (see \S\ref{margulis
tubes} for the definitions of $\neck$ and $\epshrink$).
Now for any $h\in\pleat_\rho(v\union w)$, since it also maps
$\gamma_v$ to its geodesic representative, it follows that
$\ell_{\sigma_h}(u) \ge \ep'$. Therefore
there is a curve intersecting $\gamma_u$ of $\sigma_h$-length at most
$L'\equiv\max(L,\bersfcn(\ep'))$  (see \S\ref{margulis tubes}),
and thus
$\pi_Y(\Pi_{\rho,L'}(v\union w))$ is nonempty. As in the $\xi(Y)\ge 4$ case 
we therefore obtain a bound on 
$$
\diam_Y(\Pi_{\rho,L'}(v)\union \Pi_{\rho,L'}(w))
$$
and since  this contains $\Pi_{\rho,L}(v)\union \Pi_{\rho,L}(w)$, we
obtain Property  (P2).

\subsubsection*{Proof of Property (P1)}
Let $\gamma=\gamma_v$.
Since $\Pi_\rho(v)$ is the union of $\short(\sigma_f)$ for
$f\in\pleat_\rho(\gamma)$, it suffices to bound
$$
\diam_Y(\short(\sigma_f) \union \short(\sigma_g))
$$
for any two $f,g\in\pleat_\rho(\gamma)$.

Note first that there is a $K(L)$ such that
\begin{equation}\label{bound  short}
\diam_Y(\short_L(\sigma))\le K
\end{equation}
for any hyperbolic metric $\sigma$ on $S$. This is because two curves
of bounded length in the same metric have a 
bound on their intersection number, and this readily gives a bound
on the distance between the corresponding curve systems in
$\AAA(Y)$.
From now on we will assume that $\pi_Y(\short_L(\sigma_f))$ and
$\pi_Y(\short_L(\sigma_g))$ are nonempty, since otherwise there is
nothing to prove (this is of consequence only if $Y$ is an annulus --
in other cases it holds automatically).

Thus our goal will be to find a curve $\alpha$, intersecting $Y$ essentially,
whose length  in {\em both} $\sigma_f$ and $\sigma_g$ is
bounded by some a-priori $L''\ge L$. This would give us a nonempty
intersection $\pi_Y(\short_{L''}(\sigma_f))\intersect
\pi_Y(\short_{L''}(\sigma_g))$, and hence by (\ref{bound short}) bound
the diameter of the union.

\subsubsection*{Bounding bridge arcs}
Let us invoke some machinery developed in \cite{minsky:kgcc}.
First, by Lemma 3.3 of \cite{minsky:kgcc}, up to precomposing
$g$ with a homeomorphism 
isotopic to the identity, we may (and will) assume that $f$ and $g$ are 
{\em homotopic rel $\gamma$}: this means that $f|_\gamma = g|_\gamma$
and $f$ and $g$ are homotopic keeping the points of $\gamma$ fixed. 
A {\em bridge arc} for $\gamma$ is an arc with endpoints on
$\gamma$, which is not homotopic rel endpoints into $\gamma$. 

We have the following lemma, which is part (1) of Lemma 3.4 of
\cite{minsky:kgcc}: 

\begin{lemma}{Short bridge arcs}
Fixing $\hat\ep$, for any $\delta_1$ there exists
$\delta_0\in(0,\delta_1)$ such that, if $f,g \in
\pleat_\rho(\gamma)$ and are homotopic rel $\gamma$, then 
for any bridge arc $t$ for $\gamma$ in the $\hat\ep$-thick part of
$\sigma_f$ we have 
$$
\ell_{\sigma_f}(t) \le \delta_0 \implies \ell_{\sigma_g}(t) \le
\delta_1.
$$
\end{lemma}

This lemma is a fairly direct consequence of Thurston's Uniform
Injectivity Theorem \cite{wpt:I}, which uses a compactness argument on
the space of all pleated maps to conclude that different leaves of a
pleating lamination cannot line up too closely in the image, in a
uniform sense.

\medskip

Now let $\ep_2 = \min(\ep_1/2,\neck(L))$.
Let $\Gamma$ denote the system of simple closed
$\sigma_f$-geodesics in $S$ whose 
$\sigma_f$-lengths are less than $\ep_2$, and which intersect $\gamma$
essentially (on a first reading one should consider the case that
$\sigma_f$ is $\ep_2$-thick, and in particular $\Gamma = \emptyset$).
For each component $\alpha$ of $\Gamma$,  its $f$-image
must be in the $\ep_2$ thin part of $N$. Since $\alpha$ crosses  $\gamma$,
and $f$ and $g$ are homotopic rel $\gamma$, $g(\alpha)$ meets the
$\ep_2$ thin part of $N$ as well. 
By the choice of $\ep_1$ (see \S\ref{margulis tubes}),
the $\sigma_g$-length of the $\sigma_g$-geodesic representative of
each component of $\Gamma$ is bounded by 
$\ep_0$. Hence if $Y$ crosses $\Gamma$ essentially, we may let
$\alpha$ be a a component of $\Gamma$ intersecting $Y$, and we are done.

Assume therefore that $Y$ has no essential intersections with
$\Gamma$. If $\xi(Y)\ge 4$, $Y$ must be an essential subsurface in 
a component $R$ of $S\setminus(\collar(\Gamma,\sigma_f)\union \collar(\boundary
S,\sigma_f)$.
If $\xi(Y)=2$ then there is the possibility that $Y$ is the collar of
a component of $\Gamma$. However, since $\pi_Y(\short_L(\sigma_f))\ne
\emptyset$, there is a curve of length at most $L$ crossing $Y$, and
because $\ep_2\le\neck(L)$ this is not true for any component of
$\Gamma$. Thus  $Y$ is a nonperipheral annulus in $R$.

\subsubsection*{Train track}
We now apply Lemma \ref{collapse exists} from 
\S\ref{tracks} to find a good train-track approximation of $\gamma$. 

Setting $\hat\ep = \ep_2/2$ and $\delta_1=\ep_2$, let $\delta_0$ be
the constant provided by Lemma \ref{Short bridge arcs}.
Now let $\ep = \min(\ep_2/4,\delta_0)$. Setting
the Margulis constant $\bar\ep$ in Lemma \ref{collapse exists}
to $\ep_2$, we obtain an
$\ep$-collapse $q:int(S)\to S'$ (with respect to the metric
$\sigma_f$) where $S'$ is homeomorphic to
$int(S)$, taking $\gamma$ to  a train track $\tau$.

Now since $\tau$ carries $\gamma$ and $\gamma$ has an essential
intersection with some curve $\beta$ in $Y$ (possibly a
component of $\boundary Y$), it follows from Lemma \ref{filling track}
that there exists $\alpha'$ carried in $\tau$ with 
$\ell_\tau(\alpha')\le A$ (with $A$ a bound depending only on $S$), such
that $\alpha'$ also intersects $\beta$ essentially. We will use
$\alpha'$ to construct our curve of bounded length in both $\sigma_f$
and $\sigma_g$.

\subsubsection*{Bridge arcs and nice representatives}
Using the discussion in 
\S\ref{tracks}, $\alpha'$ is homotopic to 
$$
 \alpha^\tau = \alpha_1 * t_1 * \cdots * \alpha_n * t_n
$$
where $n= \ell_\tau(\alpha')\le A$, each $\alpha_i$ runs along $\gamma$,
and each $t_i$ is in the $q$-preimage of a switch.
Lemma \ref{collapse exists} gives us a uniform
bound $K(\ep)$ on the length of the intersection of each $\alpha_i$ with the
$\ep_2$-thick part of $(S,\sigma_f)$, and a bound of $\ep$ on each
$t_i$ that meets the 
$\ep_2$-thick part. We also note that, by definition of $\Gamma$, no
$\alpha_i$ crosses a thin 
collar of a curve of $\sigma_f$-length $\le \ep_2$, unless that curve is in
$\Gamma$.
It follows that $\alpha_R=\alpha^\tau\intersect R$ remains in the
$\ep_2$-thick part,  and hence has uniformly
bounded length. The endpoints of $\alpha_R$ lie in $\boundary R$,
which is part of $\boundary \collar(\Gamma,\sigma_f) \union \boundary
\collar(\boundary S,\sigma_f)$. 
Since $\alpha'$ intersects $Y$ essentially and $Y$ is an essential
subsurface of $R$, there must be some arc $a$ of $\alpha_R$ that
intersects $Y$ essentially.

We already have a uniform bound on $\ell_{\sigma_f}(a)$. To view
this arc in $\sigma_g$, note the following: Each arc $\alpha_i$ has
exactly the same length in $\sigma_g$ as in $\sigma_f$. Each $t_i$
is a bridge arc for $\gamma$, and since the ones occuring in $a$ 
touch the $\ep_2$-thick part of $\sigma_f$,
their length is bounded by $\ep$. Since
$\ep<\ep_2/4$, the entire $t_i$ is contained in the $\ep_2/2$
thick part. Now since $\ep\le\delta_0$,
then we may apply Lemma \ref{Short bridge arcs} to 
conclude that each  $t_i$
may be deformed rel endpoints to an arc of $\sigma_g$-length at most
$\ep_2$.  Thus, $a$ may be deformed (rel endpoints if it is an arc) to
have uniformly bounded $\sigma_g$-length. 

If $a$ is a  closed curve, this is our desired curve $\alpha$, and we
are done. 

If $a$ is an arc, 
its endpoints lie in standard collars of $\Gamma$, whose boundaries
have $\sigma_f$-length $\ep_1$.
We may assume (possibly adjusting the curve $\alpha^\tau$ slightly)
that these endpoints lie on $\gamma$. Thus
their $g$-images agree with their $f$-images, and lie in the
$\ep_1$-thin parts of $N$ -- hence the endpoings lie in the 
$\ep_0$-thin parts of $\sigma_g$.
Now we can do the same mild surgery in both surfaces:
form a small regular neighborhood of the union of $a$ with the 
standard collar (or collars) meeting its endpoints. 
Let $\alpha$ be the boundary component of this neighborhood 
which passes close to $a$ 
once or twice, and makes one or two additional trips along the collar
boundaries. These arcs again are bounded in both metrics, and
$\alpha$ still intersects $\beta$ essentially. (See Figure \ref{surgery}).
This concludes the proof of Theorem \ref{Relative Coarse
Projection}.

\realfig{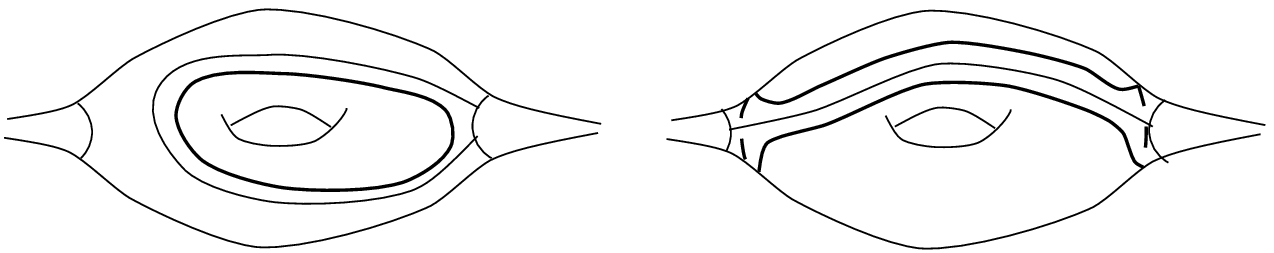}{The surgery yielding the final
curve $\alpha$ from the arc $a$. The two possible configurations are shown.}

\subsubsection{A variation on Theorem \ref{Relative Coarse
Projection}}
\label{projection variant}

It will be useful to have the following minor variant of the Coarse
Lipschitz property (P2) of 
Theorem \ref{Relative Coarse Projection}, for a special case in the
proof of Theorem \ref{Projection Bounds}.  If $W\subset S$
and $a$ is a vertex in $\AAA(W)$ represented by an arc, let $\lambda_a$
be the lamination obtained by spinning this arc leftward around
$\boundary W$ (as in Figure \ref{spin}) and including $\boundary
W$. Using this, we can define 
$\pleat_\rho(a)$ to be $\pleat_\rho(\lambda_a)$, and this allows us
to define $\Pi_{\rho,L}(a)$. 

\begin{lemma}{RCP for arcs}
Let $Y\subseteq W \subseteq S$, with $\boundary W \ne \emptyset$.
Let $\rho\in\DD(S)$ be  a Kleinian surface group such that 
$\ell_\rho(\boundary W) < \ep_0$.
If $v,w\in\AAA_0(W)$ and $d_W(v,w) \le 1$, then
\begin{equation*}
\diam_Y(\Pi_{\rho,L}(v)\union\Pi_{\rho,L}(w)) \le D_2
\end{equation*}
where $D_2$ depends only on $S$ and $L$.
\end{lemma}

\begin{proof}
The proof proceeds as in Theorem \ref{Relative Coarse Projection}, 
where the main point is to bound 
$$
\diam_Y(\short_L(\sigma_f) \union \short_L(\sigma_g))
$$
for $f,g\in \pleat_\rho(v)$, where $v$ is now a vertex of
$\AAA(W)$. We construct a train track $\tau$ from $\lambda_v$ in the
metric $\sigma_f$, and
note that $\tau$ has arcs that enter the collar of each component of
$\boundary W$ that contains an endpoint of $v$. We can consider
$\tau\intersect W$ as a train track in $W$,
and find as before an element of $\AAA(W)$ carried on $\tau$
which has bounded length, outside thin parts, in both $\sigma_f$ and
$\sigma_g$. To obtain a closed curve in $\CC_0(W)$, we do the same
surgery construction as in Figure \ref{surgery}, but using the collars
of $\boundary W$ rather than those of $\boundary S$. Because
$\ell_\rho(\boundary W) < \ep_0$, the surgered curve still has bounded
length.

\end{proof}


\section{The projection bound theorem and consequences}
\label{projection bounds}

In this section we will associate to a pair $\nu$ of end invariants 
a hierarchy $H_\nu$, and prove Theorem
\lref{Projection Bounds}, which in particular states that
$$
d_Y(v,\Pi_{\rho,L}(v)),
$$
where defined, is uniformly bounded above for all vertices $v$ in
$H_{\nu(\rho)}$ and essential subsurfaces $Y\subseteq S$.
Theorem \ref{Projection Bounds} generalizes Theorem 3.1 of \cite{minsky:boundgeom},
which only applies to  the case $Y=S$.

This bound should be taken as an indication that the vertices in
$H_{\nu(\rho)}$ and the 
bounded-length curves in $N_\rho$ are somehow close to each other
in a combinatorial sense. Indeed we will deduce the following
two corollaries of this theorem: 

The Tube Penetration Lemma \ref{Tube Penetration} controls which
pleated surfaces can penetrate deeply into a Margulis tube in $N$. In
particular it states if we pleat along the curves of a pants
decomposition coming from a slice of $H_\nu$, then the resulting surface
cannot enter any $\ep$-tubes (for a certain $\ep$) except those
corresponding to the pants curves.

The Upper Bound Lemma \ref{Upper Bounds} then shows that there is a
uniform upper bound on the length $\ell_\rho(v)$ for {\em every} $v$
appearing in $H$. This is the main step to obtaining Lipschitz bounds
on the model map  in Section \ref{lipschitz}.

\subsection{From end invariants to hierarchy}
\label{define Hnu}

Given a pair $\nu=\nu_\pm$ of end invariants we now produce a hierarchy
$H_\nu$. This is done by associating to $\nu_+$ and $\nu_-$ a pair of
generalized markings $\mu_+$ and $\mu_-$, and then applying Lemma
\lref{Existence of Hierarchies}. 

The ending laminations in $\nu_\pm$ will be part of the base of the
markings, so that what is left to do is encode the Teichm\"uller
data. Note that $\nu_\pm$ are not uniquely recoverable from $\mu_\pm$.
In what follows let $\nu_s$ denote $\nu_+$ or $\nu_-$. 

If $R$ is a component of $R^T_s$ (see \S\ref{cores and ends} for notation)
then the hyperbolic structure $\nu_R\in\Teich(R)$ admits a pants
decomposition of total length at most $L_0$ (see \ref{thick thin}),
which includes all curves of length bounded by $\ep_0$. 
Let $\mu_s(R)$ be a maximal clean
marking in $R$ whose base is such a pants decomposition, and whose
transversals $\bar t_i$ are taken to be  of minimal possible length
(There is a bounded number of choices).
Note that, by Bers' inequality relating lengths on the conformal
boundary with lengths in the interior of a hyperbolic manifold
(see \cite{bers:boundaries}), we
have $\ell_\rho(\base(\mu_s(R))) \le 2L_0$.

Recall that $\nu^L_s$ is the union of parabolic
curves $p_s$ and  ending laminations associated to the $s$ side.
Define $\mu_s$ to be the clean marking whose base is the union of
$\nu^L_s$ and $\base(\mu_s(R))$ for all components $R$ of $R^T_s$, and
whose transversals are the transversals of the markings $\mu_s(R)$.

By the discussion in Section \ref{end invariants}, $\mu_+$ and $\mu_-$
do not share any infinite-leaf components. Thus we can apply Lemma
\ref{Existence of Hierarchies} to conclude that there exists a
hierarchy $H_\nu$ with $\I(H_\nu) = \mu_-$ and $\T(H_\nu) = \mu_+$.

Note that $\base(\mu_\pm)$ is maximal, and therefore $H_\nu$ is {\em
$4$-complete} and in particular the base of every slice of $H_\nu$ is
a pants decomposition (see discussion in \S\ref{descent theorem}).

Note also that $H_\nu$ is not uniquely defined by $\nu$, as there were
choices in the construction of $\mu_\pm$, and there are choices in the
construction of a hierarchy. However, our
results will hold for any choice of $H_\nu$, and we emphasize that no
properties of the representation $\rho$ other than its end invariants
are used in the construction.

\medskip

The main theorem of this section can now be stated: 

\begin{theorem+}{Projection Bounds}
Fix a surface $S$. 
There exists $L_1\ge L_0$ such that for every 
$L\ge L_1$ there exist $B,D_2>0$ such that, 
given $\rho\in\DD(S)$, a hierarchy $H=H_{\nu(\rho)}$,
and an essential  subsurface $Y$ in $S$ with $\xi(Y)\ne 3$, the set
$$
\pi_Y(\CC(\rho,L))
$$
is $B$-quasiconvex in $\AAA(Y)$.  Furthermore,
\begin{equation}\label{eq:projection bound}
d_Y(v,\Pi_{\rho,L}(v)) \le D_2
\end{equation}
for every vertex $v$ appearing in $H$ such that the left-hand
side is defined.
\end{theorem+}
Note that the left-hand side of (\ref{eq:projection bound}) is defined
provided both $\pi_Y(v)$ and $\pi_Y(\Pi_{\rho,L}(v))$ are non-empty.
The former is satisfied whenever $\gamma_v$ intersects $Y$ essentially.
For the latter, since $\Pi_{\rho,L}(v)$ contains a pants
decomposition, its $\pi_Y$-image is always nonempty provided $\xi(Y)\ge 4$.

\subsection{Quasiconvexity}
We first consider the bound (\ref{eq:projection bound}) in a case
that is only a slight perturbation of the result proved in
\cite{minsky:boundgeom}.
After proving this version, we will prove 
Theorem \ref{Projection Bounds} in the non-annulus case
in \S\ref{sec:non annulus proj bound}. As a consequence of 
this we will obtain the Tube Penetration Lemma in \S\ref{sec:tube
penetration}, 
and then in \S\ref{sec:annulus proj bound} we will 
complete the proof of Theorem \ref{Projection Bounds} for the case of annuli.

\begin{lemma+}{Quasiconvexity}
Fix $L\ge L_0$, and 
suppose that $h\in H$ with $\xi(D(h))\ge 4$ satisfies: 
\begin{equation}\label{endpoint bound}
d_{D(h)}(u,\Pi_{\rho,L}(u)) \le d
\end{equation}
whenever $u$ is a vertex of 
$\simp(\I(h))$ or $\simp(\T(h))$. Then
\begin{equation}\label{base bound}
d_{D(h)}(v,\Pi_{\rho,L}(v)) \le d'
\end{equation}
for all simplices $v$ in $h$, where $d'$ depends on $d$ and $L$.
\end{lemma+}

\begin{proof}
We will need the following lemma, which we proved in 
\cite{minsky:boundgeom}:
\begin{lemma}{Lipschitz implies Quasiconvex}
{\rm (Lemma 3.3 of \cite{minsky:boundgeom})}
Let $\XX$ be a $\delta$-hyperbolic geodesic metric space and $\YY\subset
\XX$ a subset 
admitting a map $\Pi:\XX\to \YY$ which is coarse-Lipschitz and
coarse-idempotent. That is, there exists $C>0$ such that 
\begin{itemize}
\item[(Q1)] If $d(x,x') \le 1$ then $d(\Pi(x),\Pi(x'))\le C$, and
\item[(Q2)] If $y\in \YY$ then $d(y,\Pi(y))\le C.$
\end{itemize}
Then $\YY$ is quasi-convex, and furthermore if $g$ is a geodesic in
$\XX$  whose 
endpoints are within distance $a$ of $\YY$ then 
$$
d(x,\Pi(x)) \le b
$$
for some $b=b(a,\delta,C)$, and every $x\in g$.
\end{lemma}

The proof uses a variation of the ``stability of quasigeodesics''
argument originating in Mostow's rigidity theorem.

Proceeding with the proof of Lemma \ref{Quasiconvexity}, 
let us consider first the case that $h$ is a finite geodesic, with
endpoints $u,u'$ for which condition (\ref{endpoint bound}) holds.
Let $Z=D(h)$. In this case we will not use the fact that $h$
is in the hierarchy -- the Lemma will hold for any finite geodesic $h$.

\subsubsection*{Short boundary case}
Consider now the case that $\ell_\rho(\boundary Z) < \ep_0$. 
If $x\in\CC_0(Z)$ then, letting $f\in\pleat_\rho(x\union\boundary Z)$,
the induced metric $\sigma_f$ satisfies $\ell_{\sigma_f}(\boundary Z)
< \ep_0$, and therefore 
by our choice of $L_0$ (\S\ref{margulis tubes}),
there is a pants decomposition of $S$ of length at most $L_0\le L$
which contains $\boundary Z$ as components, and hence also contains
elements of $\CC_0(Z)$. Hence
\begin{equation}\label{short curves in Z}
\Pi_{\rho,L}(x) \intersect \CC_0(Z) \ne \emptyset.
\end{equation}

In order to apply Lemma \ref{Lipschitz implies Quasiconvex}, set
\begin{align*}
\XX&= \CC_1(Z),\\
\YY&= \CC(\rho,L) \intersect \CC_0(Z),
\end{align*}
and define a map $\Pi:\XX\to\YY$ by
letting $\Pi(x)$ be an arbitrary choice of vertex in
$\Pi_{\rho,L}(x) \intersect \CC_0(Z).$ Note that $\XX$ is
$\delta$-hyperbolic by Theorem \ref{hyperbolicity thm}.

Hypothesis (Q2) of Lemma \ref{Lipschitz implies Quasiconvex} follows
from property (P3) of Theorem \ref{Relative Coarse Projection}, noting
that any vertex $y$ of $\YY$ is also a  vertex of $\CC(\rho,L)$ that satisfies
$\pi_Z(y) \ne \emptyset$. 

If $\xi(Z)>4$ then hypothesis (Q1) follows from property (P2)
of Theorem \ref{Relative Coarse Projection}, since two adjacent
vertices of $\CC(Z)$ are also adjacent vertices of $\CC(S)$.

If $\xi(Z)=4$, then two adjacent vertices $x,x'$ of $\CC(Z)$
represent  curves that intersect, and hence are {\em not} adjacent
in $\AAA(S)$, so we cannot apply Theorem
\ref{Relative Coarse Projection} directly.  However, there are
vertices $a,a'\in\AAA_0(Z)$ (represented by arcs) such that
$[x,a],[a,a']$ and $[a',x']$ are edges in $\AAA(Z)$. We can make sense
of $\Pi_{\rho,L}(a)$ and $\Pi_{\rho,L}(a')$ as in \S\ref{projection
variant}, and apply Lemma \ref{RCP for arcs} (with $Y=W=Z$) to bound
$\diam_Z(\Pi_{\rho,L}(v)\union \Pi_{\rho,L}(w))$
for $(v,w) = (x,a), (a,a')$ and $(a',x')$. We then conclude (Q1) for
$x$ and $x'$ via the triangle inequality.

For each endpoint $w$ of  $h$, property (P1) of 
Theorem \ref{Relative Coarse Projection} bounds
$\diam_Z(\Pi_{\rho,L}(w))$, and the hypotheses of Lemma
\ref{Quasiconvexity} tell us that $d_{Z}(w,\Pi_{\rho,L}(w)) \le d$. 
Together this bounds $d_\XX(w,\Pi(w))$, and hence we can apply 
Lemma \ref{Lipschitz implies Quasiconvex} to get the desired bound
(\ref{base bound}) on $d_Z(v,\Pi_{\rho,L}(v))$ for each simplex $v$ in $h$.

\subsubsection*{Long boundary case}
Now consider the case that $\ell_\rho(\boundary Z) \ge \ep_0$. 
The main theorem of \cite{minsky:kgcc} states that there is a constant
$K$ depending only on $L$, $\ep_0$ and $S$ such that
$$
\diam_Z{\CC(\rho,L)} \le K.
$$
This means, since the endpoints of $h$ are within $d$ of their
$\pi_Z\circ\Pi_{\rho,L}$ images which are in $\pi_Z(\CC(\rho,L))$, that the
length of $h$ is at most $2d+K$. The desired bound (\ref{base bound})
now follows from the relative coarse Lipschitz property (P2) of 
Theorem \ref{Relative Coarse Projection} (note that $\pi_Z\circ
\Pi_{\rho,L}$ is never empty since $\xi(Z)\ge 4$).
This concludes the proof of the lemma when $h$ is a finite geodesic.

\medskip

Suppose that 
$\base(\T(h))$ has no finite vertices, meaning it is an element of $\EL(Z)$.
Then $\base(\T(h))$ is a component of  $\base(\T(H))$, and in particular
an ending lamination component of $\nu_+(\rho)$. The structure of
ending laminations (\S\ref{cores and ends}) implies that 
$\boundary Z$ must be parabolic in $\rho$, 
and hence we are in the short boundary case above. 
Thurston's Theorem \ref{EL exist} gives us
a sequence $\{\alpha_i\}_{i=1}^{\infty}$ 
in $\CC_0(Z)\intersect \CC(\rho,L_0)$, 
such that $\alpha_i \to \T(h)$ as $i\to \infty$.
A similar statement is true if $\I(h)$ has no finite vertices. 

Note that, since the parabolics in $\nu_+$  and $\nu_-$ must be
distinct, if $\T(h)$ has no finite vertices then $\I(h)$ has a finite
vertex, and vice versa,  unless $h=g_H$ and $Z=S$.
Thus let us assume now that $\I(h)$ has a finite vertex $\alpha_0$.

The geodesics
$h_i=[\alpha_0,\alpha_i]$ then satisfy the conditions of Lemma
\ref{Lipschitz implies Quasiconvex}, and hence the bound (\ref{base
bound}) holds for all vertices of $h_i$. Because $\CC(Z)$ is
$\delta$-hyperbolic and $\base(\T(h))$ is a point in $\boundary\CC(Z)$ by
Klarreich's theorem \ref{EL is boundary}, the $h_i$
are fellow-travelers of $h$ on larger and larger subsets. That is, for
any simplex $u$ in $h$, for large enough $i$, $u$ is at most 
$\delta$ from a vertex $u'$ in $h_i$. 

Connecting $u$ to $u'$ by a geodesic in $\CC_1(Z)$, we can apply
the Coarse Lipschitz property (P2) of $\Pi_{\rho,L}$ to bound
$d_Z(\Pi(u),\Pi(u'))$. We already have the bound (\ref{base bound}) on
$d_Z(u',\Pi(u'))$, so the triangle inequality then bounds
$d_Z(u,\Pi(u))$. 

The case that $\I(h)$ and $\T(h)$ are both infinite (and $Z=S$) is
similar: there is a biinfinite sequence
$\{\alpha_i\}_{i=-\infty}^{\infty}$ so that
$h_i=[\alpha_{-i},\alpha_{i}]$ are fellow travelers of $h$ on larger
and larger segments, and $\alpha_i\in \CC(\rho,L_0)$ for all $i$. The
bound is then obtained in the same way.
\end{proof}

\subsection{Proof of the projection bound theorem: non-annulus case}
\label{sec:non annulus proj bound}
To proceed with the proof of Theorem \ref{Projection Bounds}, 
we will need the following lemma:
\begin{lemma}{nonempty Pi}
  There exists $L$ such that, for any essential subsurface $Y$ with
  $\xi(Y)\ne 3$,  and any
  vertex $v$ in the hierarchy $H_{\nu(\rho)}$ with $\pi_Y(v)\ne \emptyset$, we have
$$
\pi_Y(\Pi_{\rho,L}(v)) \ne \emptyset.
$$
\end{lemma}

This will allow us to effectively use the Coarse Lipschitz statement of
Theorem \lref{Relative Coarse Projection}.

If $\xi(Y)\ge 4$, the lemma evidently holds with
$L\ge L_0$, since as we have already observed, 
$\Pi_{\rho,L}(v)$ always contains a pants decomposition, which must
intersect $Y$ essentially. The proof of the lemma for $\xi(Y)=2$ will 
be postponed to Section \ref{sec:annulus proj bound}.

Thus, although the rest of this section is written to be valid for any
$Y$, we may only apply it for $\xi(Y)=2$ after Lemma \ref{nonempty Pi}
has been established in that case. 
Let us henceforth assume that $L$ has been given at least as large as
the constant in Lemma \ref{nonempty Pi}, and denote $\Pi_\rho= \Pi_{\rho,L}$.

\subsubsection*{Bounds between levels}
We next consider conditions on a geodesic $h$ that allow us to
establish (\ref{eq:projection bound}) for $Y\subseteq D(h)$
and for simplices appearing in $h$.

\begin{lemma}{nesting step}
Suppose $h\in H$ with $\xi(D(h))\ge 4$, $Y\subseteq D(h)$, and
\begin{enumerate}
\item $\displaystyle d_{D(h)}(v,\Pi_\rho(v)) \le d $
for all simplices $v$ in $h$, 
\item $\displaystyle d_Y(u,\Pi_\rho(u)) \le d $
if $u$ is a vertex of $\simp(\I(h))$ or $\simp(\T(h))$ that intersects $Y$.
\end{enumerate}

Then
\begin{equation}\label{inductive conclusion}
d_Y(v,\Pi_\rho(v)) \le d'
\end{equation}
for all simplices $v$ of $h$ which intersect $Y$, where $d'$ depends
on $d$. 
\end{lemma}

\begin{proof}
Let $v$ be a simplex of $h$ that intersects $Y$. Then $v$ is not in
$\phi_h(Y)$, and let us assume without loss of generality that
$\max\phi_h(Y)<v$.

Suppose that the distance from $v$ to the last simplex $\omega$ of $h$
is no more than $2d+1$.  
Let $\alpha$ be the segment of $h$ from $v$ to $\omega$. 
Since every simplex
in $\alpha$ crosses $Y$, the 1-Lipschitz property for $\pi_Y$ implies
$$
\diam_Y(\alpha) \le  2d+1.
$$
The Relative Coarse Lipschitz property (P2) of $\Pi_\rho$ in Theorem
\ref{Relative Coarse Projection} says that for every two successive
simplices $x$,$x'$ of $\alpha$ we have 
$$
\diam_Y(\Pi_\rho(x) \union \Pi_\rho(x')) \le D_1.
$$
By Lemma \ref{nonempty Pi},
$\pi_Y(\Pi_\rho(x))$ is nonempty for each
$x$ in $\alpha$, so we can sum this over  $\alpha$ to obtain
$$
\diam_Y(\Pi_\rho(\alpha)) \le D_1(2d+1).
$$
Finally, we have by hypothesis (2) the  bound
$$d_Y(\omega,\Pi_\rho(\omega))\le d.$$
Noting that $\pi_Y(\omega)\in\pi_Y(\alpha)$ and
$\pi_Y(\Pi_\rho(\omega))\subset 
\pi_Y(\Pi_\rho(\alpha))$, we can then put these three bounds together
to obtain
$$
\diam_Y(\alpha\union \Pi_\rho(\alpha)) \le d+(D_1+1)(2d+1).
$$
In particular this bounds $d_Y(v,\Pi_\rho(v))$, and we have the
desired statement.

Alternatively, suppose that the distance from  $v$ to the last simplex
of $h$ is
at least $2d+2$ (including the possibility that it is infinite).
Let $\alpha$ be the segment of length $2d+2$ beginning with $v$, and
let $w$ be the other endpoint of $\alpha$. By hypothesis we have
$$d_{D(h)}(w,\Pi_\rho(w)) \le d$$
and hence we can join $w$ with $\pi_{D(h)}(\Pi_\rho(w))$ with a path
in $\AAA_1(D(h))$ of length at most $d$, and whose first and last
vertex (by Lemma \ref{nonempty Pi}) intersect $Y$.
Using the 2-Lipschitz map
$\psi:\AAA_0(D(h)) \to \CC_0(D(h))$ described in Section
\ref{complexes}, we can replace this with a path $\beta$ with the same
endpoints and whose
interior vertices are in $\CC(D(h))$, of length at most $2d$. 
Now recalling that
$\max\phi_h(Y)<v$,  we find that $d_{\CC_1(D(h))}(w,\max\phi_h(Y)) \ge 2d+3$.
By the triangle inequality in $\CC_1(D(h))$, every point in $\beta$ is
at least distance 3 from
$\max\phi_h(Y)$ and hence at least distance 2 from $[\boundary Y]$. 
In particular every interior vertex in $\beta$ has nontrivial
intersection with $Y$. Now the 1-Lipschitz property of $\pi_Y$ again
applies, to give us $\diam_Y(\beta) \le 2d$. In particular
$$
d_Y(w,\Pi_\rho(w)) \le 2d.
$$ 
Now we can apply exactly the same argument as in the previous case. 
\end{proof}

\subsubsection*{Inductive argument}
We will now establish, by induction, the following claim:
\begin{lemma}{inductive claim}
Let $h\in H$ be a geodesic with $\xi(D(h))\ge 4$, and $Y\subset D(h)$.
If $v$ is a simplex of $h$ or $\simp(\I(h))$ or $\simp(\T(h))$,
and $\pi_Y(v)\ne \emptyset$, then
$$
d_Y(v,\Pi_\rho(v)) \le d
$$
where $d$ depends on $\xi(D(h))$ and $\xi(S)$. 
\end{lemma}

\begin{proof}
Consider first the case where $h=g_H$ and $Y=S$, and let us apply
Lemma \ref{Quasiconvexity}. 
If $g_H$ is bi-infinite there are no conditions to check. If not, suppose that
$\I(g_H)=\I(H)$ contains a finite vertex $u$. 
Then by definition of $\I(H)$ (\S\ref{define Hnu}) we have
$\ell_\rho(u)\le 2L_0$. This uniformly bounds $d_S(u,\Pi_\rho(u))$, since 
$u$ and any element of $\Pi_\rho(u)$ have lengths bounded by $2L_0$
and $L_0$ on the same pleated surface $f\in\pleat_\rho(u)$. 
Hence the condition of Lemma \ref{Quasiconvexity}
holds for $u$. The same is true for $\T(g_H)$, and we therefore have
\begin{equation}\label{main geodesic bound}
d_{S}(v,\Pi_\rho(v)) \le d
\end{equation}
for every simplex $v$ of $g_H$,
where $d$ is a constant depending only on $S$.

We may now apply
Lemma \ref{nesting step} to $h=g_H$, with
every subsurface $Y$:
Condition (1) of the lemma follows from (\ref{main geodesic bound})
which we have just proved. 
For condition (2), note as above that for any finite vertex $u$ in
$\I(H)$ or $\T(H)$, $u\in\Pi_\rho(u)$ and hence if $u$ intersects $Y$
we have $d_Y(u,\Pi_\rho(u))=0$.
Thus the conclusion of Lemma \ref{inductive claim} holds for the
geodesic $g_H$, and we have established the base case. 

\medskip

Now let $h\in H$ be any geodesic other than $g_H$, and suppose that Lemma 
\ref{inductive claim} holds for all $h'$ with $\xi(D(h'))>\xi(D(h))$,
for some constant $d$. In order to apply Lemma \ref{Quasiconvexity}
to $h$, we must consider the condition on its endpoints.
Let $b,f$ be such that $b\bsubd h \fsubd f$.
If $\T(h)$ is in $\EL(D(h))$ ($h$ is infinite in the forward
direction) then there is nothing to check. Otherwise
$\T(h)$ has a finite simplex $w$, which by definition of $h\fsubd f$ appears
either in a simplex of $f$ or in $\T(f)$.
By the induction hypothesis applied to $f$ we have
$d_{D(h)}(w,\Pi_\rho(w)) \le d$. 
The same reasoning applies to $\I(h)$, using $b$. 
Thus, we may apply Lemma \ref{Quasiconvexity} to obtain the bound 
(\ref{base bound}) for all simplices of $h$.

This then establishes condition (1) of Lemma \ref{nesting step} for
$h$. To obtain condition (2) for any $Y\subseteq D(h)$, we again use
the fact that the endpoints of $h$ are contained in $b$ and $f$, and
hence the inductive hypothesis for $b$ and $f$ yields this bound
also (for suitable constants). Thus, we apply Lemma \ref{nesting step}
to give us the statement of Lemma \ref{inductive claim} for $h$. 
\end{proof}

\subsubsection*{Reduction to nested case}
So far we have proved the bound on $d_Y(v,\Pi_\rho(v))$ in the case
where $v$ is in a geodesic $h$ with $Y\subseteq D(h)$. It remains to
show that the general case reduces to this one. 

For any geodesic $h\in
H$ with $\xi(D(h))\ge 4$ and simplex $v$ in $h$, let
$v\union \boundary D(h)$ denote the simplex of $\CC(S)$ corresponding
to the disjoint union of curves $\gamma_v$ and (the nonperipheral
compoents of) $\boundary D(h)$. 
Let us fix $Y$ and prove, by induction on $\xi(D(h))$, that  
\begin{equation}\label{simplex with boundary}
d_Y(v\union \boundary
D(h),\Pi_\rho(v\union \boundary D(h))) \le C,
\end{equation}
provided the left hand side
is defined (where $C$ is a uniform constant). 
By Lemma \ref{nonempty Pi}, all that is necessary
for the left hand side to be defined is that $v \union\boundary D(h)$ have
non-trivial intersection with $Y$. 

For any $h$ such that $Y\subseteq D(h)$, (\ref{simplex with boundary})
reduces to what we have already proven. This applies in particular to $h=g_H$.
Now given any other $h\in H$,
suppose the claim is true for all $h'$ with
$\xi(D(h'))>\xi(D(h))$. Let $h\fsubd f$. Then $D(h)$ is a component
domain of $(D(f),w)$ for some simplex $w$ in $f$, and hence
$\boundary D(h)$ is contained in $w\union \boundary D(f)$.

If $\boundary D(h)$ intersects $Y$ essentially
then so does $w\union\boundary D(f)$, and 
hence the inductive hypothesis bounds $d_Y(w\union\boundary
D(f),\Pi_\rho(w\union\boundary D(f)))$. 
Now note that $v\union\boundary D(h)$ and $w\union\boundary D(f)$ have
intersection $\boundary D(h)$, and since this has non-empty
projection $\pi_Y(\boundary D(h))$ we have a bound on 
$\diam_Y(v\union\boundary D(h) \union w\union\boundary D(f))$.
Furthermore $\Pi_\rho(\boundary D(h))$ {\em contains}
$\Pi_\rho(v\union\boundary D(h))$ and $\Pi_\rho(w\union\boundary D(f))$.
Thus, invoking Theorem \ref{Relative Coarse Projection} we may
deduce a bound on 
$$
d_Y(v\union\boundary D(h), \Pi_\rho(v\union\boundary D(h)))
$$

If $\boundary D(h)$ does not intersect $Y$ essentially, but $v$ does,
then $Y\subseteq D(h)$, and we have already proven the bound. 

This concludes the proof of Theorem \ref{Projection Bounds}, in the
non-annulus case.

\subsection{Tube penetration lemma}
\label{sec:tube penetration}

Applying the non-annulus case of Theorem \ref{Projection Bounds},
we are now going to control which pleated surfaces arising from a
hierarchy can meet Margulis tubes in $N_\rho$.

\begin{lemma+}{Tube Penetration}
There exists $\ep_3>0$ with the following property:
Let $u$ be a vertex appearing in $H_{\nu(\rho)}$ and
let $\alpha\in\pi_1(S)$ satisfy   $\ell_\rho(\alpha) < \ep_3$.
A map $$f\in\pleat_\rho(u)$$
meets the Margulis tube $\MT_{\ep_3}(\alpha)$
only if $\alpha$ represents a simple curve which has no
essential intersection with $\gamma_u$.
\end{lemma+}

\begin{corollary}{Penetration by slices}
 Let $\tau$ be a slice of $H$, and $\alpha\in\pi_1(S)$ a primitive
 element with
 $\ell_\rho(\alpha) < \ep_3$. A map 
$$
f\in \pleat_\rho(\base(\mu_\tau))
$$
meets the Margulis tube $\MT_{\ep_3}(\alpha)$ if and only if $\alpha$
represents an element of $\base(\mu_\tau)$.
\end{corollary}

The ``if'' direction of the corollary is obvious. For the ``only if''
direction, if $f$ meets $\MT_{\ep_3}(\alpha)$ then Lemma \ref{Tube
Penetration} implies that every component of $\base(\mu_\tau)$ has no
essential intersection with the simple curve $\alpha$. Since
$\base(\mu_\tau)$ is a pants decomposition, $\alpha$ must be one of
the components.

\begin{pf*}{Proof of the Tube Penetration Lemma}
We will choose $\ep_3 \le \ep_1$, and assume that 
$f(S)$ meets $\MT_{\ep_3}(\alpha)$.
By the choice of 
$\ep_1$ in \S\ref{margulis tubes} we know that 
$\alpha$ has length at most
$\ep_0$ in $\sigma_f$. Hence it is a multiple of the core of a thin collar
and since we assumed it was primitive it must in fact be the core, 
and in particular represents a vertex of $\CC(S)$. 

The vertex $u$ must appear in some simplex $v$ of a geodesic $h\in H$, with
$\xi(D(h))\ge 4$. Assuming $\gamma_v$ intersects  $\alpha$
essentially, we will prove 
that $f(S)$ can penetrate no further than a distance $d$ into
$\MT_{\ep_1}(\alpha)$, where $d$ depends on $\xi(D(h))$. The proof
will be by downward induction on $\xi(D(h))$ with $D(h)=S$ being the base
case. 

Since $v\notin \phi_h(\alpha)$, we may assume without loss of
generality that $v> \max\phi_h(\alpha)$ (otherwise we reverse the
directions in the rest of the argument). Let $J$ denote the largest
segment of $h$ beginning with $v$ so that for every vertex $x$ of
$J$ its geodesic representative $x^*$ in $N_\rho$ meets
$\MT_{\ep_1}(\alpha)$.
For any $x\in J$, every $f\in \pleat_\rho(x)$ has nontrivial intersection
with $\MT_{\ep_1}(\alpha)$. 
As above this means $\ell_{\sigma_f}(\alpha)\le\ep_0$, and in particular
$$
\alpha\in\Pi_\rho(x).
$$ 
Since $\xi(D(h))\ge 4$ we may
apply the non-annulus case of Theorem \ref{Projection Bounds} to 
obtain
\begin{equation}
  \label{eq:all near alpha}
d_{D(h)}(x,\alpha) \le b
\end{equation}
where $b$ is a constant derived from that lemma and the bound $D_0$
on $\diam_{D(h)}(\Pi_\rho(x))$ in Theorem \ref{Relative Coarse
  Projection}.
We conclude, since $J$ is geodesic in $\CC_1(D(h))$, a bound of the form
\begin{equation}
|J| \le b'.
\end{equation}
If $J$ is the entire portion of $h$ following $v$, then let $w$ denote
the last simplex of $J$ and hence the last of $h$. 
If $h=g_H$ (the base case of the induction)
then $w\subset \base(\T(H))$ and hence represents curves in
$N_\rho$ of length at most $2L_0$. This bounds by $L_0$ the distance
which $w$ can penetrate into $\MT_{\ep_1}(\alpha)$, since $w^*$ cannot be
completely contained in this Margulis tube (if it were it would be
$\alpha$ itself, but we have $w> \max \phi_h(\alpha)$).

If $h\ne g_H$ then there exists $f\in H$ such that $h\fsubd f$, and hence
$w$ is contained in a simplex of $f$. By the induction,
$w^*$ cannot penetrate further than $d_{\xi(D(f))}$ into
$\MT_{\ep_1}(\alpha)$.  

If the last simplex of $J$ is not the last simplex of $h$, then let
$w$ denote its successor, so by definition $w^*$ does not meet
$\MT_{\ep_1}(\alpha)$. Define $J' = J\union \{w\}$ (so $|J'|\le |J|+1$).

\medskip

We will now construct a path in $N_\rho$ joining any point
$q_v\in v^*\intersect \MT_{\ep_1}(\alpha)$ to the boundary of
$\MT_{\ep_1}(\alpha)$, whose length will be 
bounded. There are two possible cases: 

\subsubsection*{Case 1: $\xi(D(h))>4$}
Consider any two successive simplices $x,x'$ in $J'$,
and let $q$ be a point of $x^*\intersect
\MT_{\ep_1}(\alpha)$. Since $x$ and $x'$ determine two disjoint curves
in $D(h)$, there is a pleated surface $f\in \pleat_\rho([xx'])$. 

Again by the choice of $\ep_1$, 
$q$ is the $f$-image of a
point $p\in \gamma_x$ in the $\ep_0$-thin part of $\sigma_f$ associated to
$\alpha$. By our construction both $x$ and $x'$ cross $\alpha$, so
there is a point $p'$ on $\gamma_{x'}$ in this thin part, a
$\sigma_f$-distance of at most 
$\ep_0$ from $p$. Thus $q'=f(p')$ is connected to $q$ by a path of
length at most $\ep_0$. 

Apply this successively to all the vertices of $J'$,
beginning with $q_v$. Each new point is either outside
$\MT_{\ep_1}(\alpha)$, in which case we stop, or inside it, in which
case we continue, and reach $w^*$ in
at most $|J'|$ steps. Since $w^*$ is either disjoint from
$\MT_{\ep_1}(\alpha)$ or penetrates no more than $d_{\xi(D(h))+1}$
into it, we conclude that the distance from $q$ to
the boundary of $\MT_{\ep_1}(\alpha)$ is at most $|J'|\ep_0 +
d_{\xi(D(h))+1}$. 

\subsubsection*{Case 2: $\xi(D(h))=4$}
Now successive vertices $x,x'$ in $J'$
are {\em not} disjoint and hence we cannot form $\pleat_\rho(x\union
x')$.
Instead, extend $x\union \boundary D(h)$ and
$x'\union \boundary D(h)$ to pants decompositions $p,p'$ which differ by
an elementary move, and
consider the corresponding surfaces $g\in\pleat_\rho(p)$, 
$g'\in\pleat_\rho(p')$, and the halfway surface
$f=f_{p,p'}$ (see \S\ref{halfway}). 
The halfway surface $f$ is pleated along a lamination
containing
a leaf (or two leaves) $l$ that is part of the pleating locus of $g$
in $D(H)$, 
and a leaf or two $l'$ that 
is part of the pleating locus of $g'$ in $D(h)$. The discussion in
\S\ref{halfway} 
tells us that $\alpha$ intersects $l$ essentially since it intersects
$\gamma_x$, and intersects $l'$ since it intersects $\gamma_{x'}$.

Thus in $\sigma_g$ both $\gamma_x$ and $l$ pass through the collar of
$\alpha$, in $\sigma_f$ both $l$ and $l'$ pass through it, and in
$\sigma_{g'}$ both $l'$ and $\gamma_{x'}$ do. It follows that we
can apply the argument of case 1 to produce a path as before, but
with three times the number of steps.

\medskip

Thus we have shown that $v^*$ cannot penetrate more than a certain
$d_{\xi(D(h))}$ into $\MT_{\ep_1}(\alpha)$. 
Now using (\ref{definite tube nesting})
this implies that for a certain uniform
$\ep_3$, $v^*$ cannot meet the tube $\MT_{\ep_3}(\alpha)$. 
\end{pf*}

\subsection{Projection bounds in the annulus case}
\label{sec:annulus proj bound}

We are now ready to 
complete the proof of Theorem \ref{Projection Bounds} in the case that
$Y$ is an annulus. It suffices to give the proof of 
Lemma \ref{nonempty Pi} in this case:

\begin{pf}
Let $\alpha$ denote the core of the annulus $Y$. 
Let $\ep_3$ be the constant in Lemma \lref{Tube Penetration}, and take
$L\ge \bersfcn(\ep_3)$. Then $\pi_Y(v)\ne \emptyset$ implies that, for
any $f\in\pleat_\rho(v)$, $f(S)$ does not meet
$\MT_{\ep_3}(\alpha)$. Thus, $\ell_{\sigma_f}(\alpha) \ge \ep_3$ and it
follows by definition of the function $\bersfcn()$ (see \S \ref{thick
  thin}) that there is a pants decomposition of total length at most
$L$ which crosses $\alpha$ essentially. 
Thus, $\pi_Y(\Pi_{\rho,L}(v)) \ne \emptyset$, which gives the statement
of Lemma \ref{nonempty Pi}.
\end{pf}

The rest of the proof of Theorem \ref{Projection Bounds} proceeds
exactly as in Section \ref{sec:non annulus proj bound}.

\subsection{Upper  length bounds}
\label{sec:upper bounds}

Our final application of the projection mechanism will be to obtain an
a priori upper bound on the length of every curve that appears in the
hierarchy.

\begin{lemma+}{Upper Bounds}
For every vertex $v$ in the hierarchy $H_{\nu(\rho)}$,
$$
\ell_\rho(v) \le D_3
$$
where $D_3$ is a constant depending only on the topological type of $S$.
\end{lemma+}

\begin{proof}
Any vertex $v$ in $H$ is contained in some slice $\tau$, by Lemma
\ref{Resolution sweep}. Let $\mu =
\mu_\tau$ be the associated marking, and fix
$$
f\in\pleat_\rho(\base(\mu)) \subset \pleat_\rho(v).
$$
Let $\ep_3$ be the constant in Lemma \lref{Tube Penetration} and
Corollary \ref{Penetration by slices}. Applying 
Corollary \ref{Penetration by slices}, we find that all curves
$\alpha$ in $S$ with $\ell_{\sigma_f}(\alpha) \le \ep_3$ must be
components of $\base(\mu)$. If $\ell_{\sigma_f}(v) \le \ep_3$ then we
already have our desired bound, and we are done. Thus assume
$\ell_{\sigma_f}(v) >  \ep_3$. Since $\gamma_v$ is disjoint from the
other curves of $\base(\mu)$, it must be contained, non-peripherally,
in a component $R$ of the $\ep_3$-thick part of $\sigma_f$. 

Since $\ep_3<\ep_0$, there is a pants decomposition $\eta$ of $R$ of
$\sigma_f$-length at most $L_0$, and in particular
$$\eta\subset \Pi_{\rho,L_0}(\base(\mu)).$$
Since $R$ is $\ep_3$-thick, there is a constant $L$ depending on
$\ep_3$ so that
each component of $\eta$ is crossed by a transversal curve of
$\sigma_f$-length at most $L$, which misses the rest of $\eta$. 
Let $\mu_2$ be the clean marking with  
$\base(\mu_2)=\eta$ and these bounded-length transversal curves.

By Theorem \lref{Projection Bounds}, if $Y$ is any subsurface of $R$
that has essential intersection with both $\eta$ and $\base(\mu)$
(including $Y=R$) then 
$$
d_Y(\base(\mu),\eta) \le D_2.
$$
The only subsurfaces of $R$ excluded by this are annuli associated
with components of $\base(\mu)$ or $\eta$.  Let $Y$ be an annulus
whose core is a component of $\eta$, but not a component of
$\base(\mu)$. The transversal $t$ in
$\mu_2$ crossing $Y$, since it has $\sigma_f$-length at most $L$, 
is contained in $\Pi_{\rho,L}(\base(\mu))$. Thus applying
Theorem \ref{Projection Bounds} again, there is a $D'_2$ such that
$$
d_Y(\base(\mu),t) \le D'_2.
$$

To deal with the remaining annuli,
extend $\base(\mu)$ to a
new clean marking $\mu'$ by adding for each component $\alpha$ a transversal
$t$ satisfying $d_\alpha(t,\mu_2) \le 2$ (this is always possible by
picking some transversal and applying Dehn twists.)
Thus for all annuli $Y$ with cores in $\base(\mu)=\base(\mu')$ we also
have
$$
d_Y(\mu',\mu_2) \le 2.
$$
In other words we have shown that the two clean markings satisfy
$$
d_Y(\mu',\mu_2)\le D
$$ 
for a uniform $D$, and 
{\em all} domains $Y\subseteq R$. 
By Lemma \ref{dW and del}, 
this gives an upper bound $E$ on the elementary move distance
$d_{el}(\mu_2,\mu')$ (to obtain the bound, set $K=\max(K_0,D+1)$ in
the lemma).

Since the base curves and transversals of $\mu_2$ have 
$\sigma_f$ lengths bounded by 
$L$, and since an elementary move can change the lengths of components
of a maximal clean marking by at most a bounded factor, 
we obtain a bound on the $\sigma_f$-lengths of the curves of
$\mu'$, in terms of $L$ and $E$. In particular we have an upper bound
on $\ell_\rho(v)$, and we are done. 
\end{proof}

\section{The model manifold}
\label{model}

In this section we will construct an oriented metric 3-manifold $\modl$ associated to 
the end invariants $\nu=\nu_\pm(\rho)$. $\modl$ is intended to be a
model for the geometry of the augmented convex core $\hhat C_N$ of $N_\rho$.
The following is a summary of the structure of $\modl$.

\begin{enumerate}
\item
$\modl$ is properly embedded in $\hhat S\times\R$, and 
is homeomorphic to $\hhat C_N$.
\item
An open subset $\UU\subset \modl$ is called the set of {\em tubes} of
$\modl$. Each component $U\subset \UU$ is of the form
$$\collar(v)\times I$$
where $v$ is either
a vertex of $H_\nu$ or a boundary component of $S$,
and $I\subset \R$ is an open interval.
$I$ is bounded except for the finitely many $v$
corresponding to parabolics.
The correspondence $U \leftrightarrow v$ is bijective.

\item

Let $\modl[0]$ be $\modl\setminus\UU$. $\modl[0]$ is a
union of standard ``blocks'' of a  finite number of topological types.

\item
Except for finitely many blocks adjacent to $\boundary \modl$, 
all blocks fall into a predetermined finite number of isometry types. 

\item
Each tube $U$ is isometric to a hyperbolic or a parabolic tube, and
with respect to a 
natural marking has boundary parameters $(\Momega(U),\ep_1)$. (See
\S\ref{thick thin}). We call $\Momega(U)$ the {\em meridian
coefficient} of $U$.
\end{enumerate}

Let $H=H_\nu$ be associated to $\nu$ as in \S\ref{define Hnu}.
We will begin by  constructing $\modl[0]$ abstractly as a union of
blocks. We will then show how to embed it in $\hhat S\times\R$, and in its
complement we will find the tubes $\UU$ which we adjoin to obtain $\modl$.

\subsection{Blocks and gluing}
\label{blocks and gluing}

The typical blocks from which we build $\modl[0]$ are called {\em internal
  blocks}. 
There are some special cases of blocks associated to the
boundary of the convex core, but these can be ignored on a first
reading (and do not appear, for example, in the doubly degenerate
case).

Given a 4-edge $e$ in $H$, let $g$ be the $4$-geodesic 
containing it, and let $D(e)$ be the domain, $D(g)$. Recall
(\S\ref{more hierarchies}) that $e^-$ and $e^+$ denote the initial and
terminal vertices of $e$. 

To each $e$ we will associate a block $B(e)$, defined as follows:
\begin{align*}
B(e) =  (D(e)\times [-1,1]) \setminus & \left(\
\collar(e^-)\times[-1,-1/2)\right. \union\\ 
& \ \ \left. \collar(e^+)\times(1/2,1]\ \right).
\end{align*}
That is, $B(e)$ is the product $D(e)\times [-1,1]$, with solid-torus
trenches dug out of its top and bottom boundaries, corresponding to
the two vertices of $e$. See Figure \ref{block}.

We remark that, in this construction, we think of blocks of distinct
edges as disjoint (for example think of the intervals $[-1,1]$ as disjoint
copies of a standard interval). Afterwards we will glue them together
using specific rules, and embed the resulting manifold in
$\hhat S\times\R$. 

\realfig{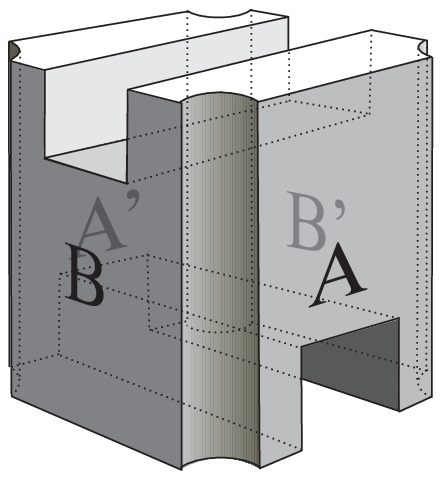}{Constructing an internal block $B(e)$.
If $D(e)$ is a one-holed torus then $B(e)$ is obtained by gluing face
$A$ to face $A'$, and $B$ to $B'$. The curved vertical faces become
$\boundary D(e)\times [-1,1]$. If $D(e)$ is a 4-holed sphere then $B(e)$ is
obtained by doubling this object along $A,A', B$ and $B'$. }

We break up the boundary of $B(e)$ into several parts (see Figure
\ref{blockschem} for a schematic).
The {\em gluing boundary} of $B(e)$ is 
$$
\boundary_\pm B(e) \equiv (D(e)\setminus\collar(e^\pm))
\times \{\pm 1\}.
$$
Note that the gluing boundary is always a union of three-holed
spheres.

The rest of the boundary is a union of annuli, with
$$
\boundvert B(e)\equiv \boundary D(e) \times [-1,1]
$$
being the {\em outer annuli}.

The {\em inner annuli} $\boundary^\pm_i B(e)$ are the boundaries of
the removed solid tori. That is, 
$$
\boundary^\pm_i B(e) = \boundary B(e) \intersect
        \boundary(\collar(e^\pm)\times\pm(1/2,1])
$$
(where $+(a,b]$ denotes $(a,b]$ and $-(a,b]=[-b,-a)$).
These annuli break up into a {\em horizontal part}
$$
\boundary^\pm_{ih} B(e) = \bcollar(e^\pm)\times\{\pm 1/2\}
$$
and a {\em vertical part}
$$
\boundary^\pm_{iv} B(e) = \boundary\collar(e^\pm)\times\pm [1/2,1].
$$

\realfig{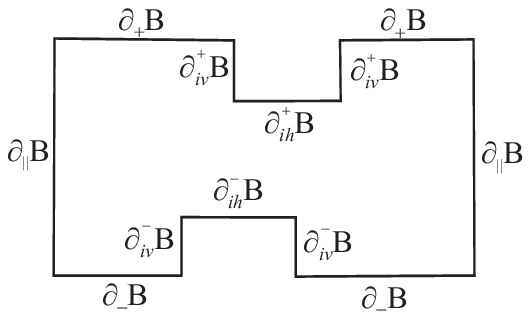}{Schematic diagram of the
  different pieces of the boundary of a block.}

\subsection*{Boundary blocks}

Recall from \S\ref{cores and ends} that $R^T_+$ denotes the union of
subsurfaces in the top of the relative compact core that face
geometrically finite ends. Let
$R$ be a subsurface of $S$ homotopic to a 
component of $R^T_+$, and let
$\nu_R$ be the associated component of $\nu_+^T$ in $\Teich(R)$.
We construct a block $B\sbtop(\nu_R)$ as follows:
Let $\T_R$ be 
the set of curves of $\base(\T(H_\nu))=\base(\mu_+)$ that are
contained in $R$. Define 
$$
B'\sbtop(\nu_R) = R\times[-1,0] \setminus \left(\collar(\T_R) \times[-1,-1/2)\right)
$$
and let
$$
B\sbtop(\nu_R) = B'\sbtop(\nu_R) \union \boundary R\times[0,\infty).
$$
This is called a {\em top boundary block}. Its {\em outer boundary}
$\boundary_o B\sbtop(\nu_R)$ is
$R\times\{0\}\union \boundary R\times[0,\infty)$, which we note is
homeomorphic to $int(R)$. This will correspond to a boundary component
of $\hhat C_N$.
The gluing boundary of this block lies on its bottom: it is  
$$
\boundary_-B\sbtop(\nu_R) = (R\setminus \collar(\T_R))\times\{-1\}.
$$
Similarly if $R$ is a  component of $R^T_-$ we let $\I_R =
\I(H_\nu)\intersect R$ and define 
$$
B'\sbot(\nu_R) = R\times[0,1] \setminus \collar(\I_R) \times(1/2,1].
$$
and the corresponding {\em bottom boundary block}
$$
B\sbot(\nu_R) = B'\sbot(\nu_R) \union \boundary R\times(-\infty,0].
$$
The gluing boundary here is $\boundary_+B\sbot(\nu_R) = 
(R\setminus \collar(\I_R))\times\{1\}.$ 

The vertical annulus boundaries are now 
$\boundvert B\sbtop(\nu_R) = \boundary R \times [-1,\infty)$
and the internal annuli
$\boundary_i^\pm$ are are a union of possibly several
component annuli, one for each component of $\T_R$ or $\I_R$.

Since each block is a subset of $Z\times[-1,1]$ for  a subsurface $Z$
of $S$, it inherits a natural orientation from the fixed orientation
of $S$ and $\R$. 

\subsection*{Gluing instructions}

We obtain $\modl[0]$ by taking the disjoint union of all blocks and identifying
them along the three-holed spheres in their gluing
boundaries. The rule is that whenever two blocks $B$ and $B'$ have the same
three-holed sphere $Y$ appearing in both
$\boundary^+B$ and $\boundary^- B'$, we identify these boundaries
using the identity on $Y$. The hierarchy will serve to organize these
gluings and insure that they are consistent.

By definition, all these three-holed spheres are component domains in
the hierarchy. Conversely, we will now check
that every component domain $Y$ in $H$ with $\xi(Y)=3$ must occur in
the gluing boundary of exactly two blocks. 

Lemma \ref{xithree configurations} tells us that
$\T(H)|_Y\ne\emptyset$ if and only if there exists a 4-geodesic $f\in
H$ with $Y\fsubd f$, and $f$ is unique if it exists. When $f$ exists, 
there is an edge $e_f$ in $f$ with $Y$ a component
domain of $(D(f),e_f^-)$, and hence $Y\times\{-1\}$ is a component
of $\boundary_- B(e_f)$. 

If $f$ does not exist and $\T(H)|_Y = \emptyset$, $Y$ must be a component
domain of $\base(\T(H))$, and hence $Y\times\{-1\}$ occurs on the
gluing boundary $\boundary_- B\sbtop(\nu_R)$ for some top boundary
block. 

Similarly if $\I(H)|_Y \ne \emptyset$ then $b\bsubd Y$ for a
unique $4$-geodesic $b$, and 
$Y\times\{1\}$ is in $\boundary_+ B(e_b)$, and if $\I(H)|_Y=\emptyset$
then $Y\times\{1\}$ appears in the gluing boundary of a bottom
boundary block. 

We conclude that each $Y$ serves to glue exactly two blocks, and in
particular $\modl[0]$ is a manifold. The orientations of the blocks
extend consistently to an orientation of $\modl[0]$.

\subsection{Embedding in $\hhat S\times\R$}
\label{embedding the model}

The interior of
each block $B(e)$ inherits a 2-dimensional ``horizontal foliation'' from the 
foliation of the product
$D(e)\times [-1,1]$ by surfaces $D(e)\times\{t\}$.
Let us call the connected leaves
of this foliation the {\em level surfaces} of $B(e)$.
(For boundary blocks we do similarly, and for the added
vertical annuli of the form $\boundary
R\times[0,\infty)$ or $\boundary R \times (-\infty,0]$, we also
include the level circles as leaves of this foliation.)
An embedding
$f:\modl[0]\to S\times\R$ will be called {\em flat} if each connected
leaf of the horizontal foliation
$Y\times\{t\}$ in $\modl[0]$ is mapped to a level set $Y\times\{s\}$
in the image, 
with the map on the first factor being the identity.

\begin{theorem}{embed model}
$\modl[0]$ admits a proper flat orientation-preserving embedding
$\Psi:\modl[0]\to S\times\R$. 
\end{theorem}

\subsubsection*{Remark} Once we have fixed the embedding $\Psi$ we will 
adopt notation where $\modl[0]$ is identified with its image in
$S\times\R$, and $\Psi$ is the identity.

\begin{proof}
In the course of the proof we will fix a certain exhaustion of
$\modl[0]$ (minus its boundary blocks)
by subsets $M_i^j$, where $M_i^j$ and $M_i^{j+1}$ differ by the
addition of one block on the ``top'', and $M_{i-1}^j$ and $M_i^j$
differ by the addition of a block on the bottom. The map will be built
inductively on $M_i^j$. 

Let $\{\tau_i\}_{i\in \II}$ be a resolution of $H$, as in Lemma
\ref{Resolution sweep}. 

The pants decompositions $\base(\mu_{\tau_i})$ form a sequence
that may have adjacent repetitions -- that is, $\base(\mu_{\tau_i})$ and 
$\base(\mu_{\tau_{i+1}})$ may be the same because the move
$\tau_i\to\tau_{i+1}$ involves a $\xi=2$ geodesic (twists in an
annulus complex) or a $\xi>4$ geodesic (reorganization moves).
If we remove such repetitions we obtain a sequence of pants
decompositions $\{\eta_i\}_{i\in\II'}$, where $\II'$ is a new index
set, so that each step $\eta_i\to\eta_{i+1}$ corresponds to an edge in
a $\xi=4$ geodesic. For a vertex $v$ appearing in $H$, define
$$
J'(v) = \{i\in\II': v\in \eta_i\}.
$$
Lemma \ref{vertex interval} implies that $J(v)$ is an interval in $\Z$, and it
follows that $J'(v)$ is an interval as well.

To each $\eta_i$ we associate a subsurface $F_i$ of $\modl[0]$ whose
components are level surfaces of blocks, as follows. 
A complementary component $Y$ of $\eta_i$ is necessarily a three-holed
sphere which appears as a component domain in $H$. 
Thus there are blocks $B_1$ and $B_2$ such that $Y$ is isotopic
to a component of $\boundary_+B_1$ and to a 
component of $\boundary_-B_2$, which are identified in $\modl[0]$. Let
this identified subsurface be a 
component of $F_i$. Repeating this for all complementary components of $\eta_i$
we obtain all of $F_i$, which we call 
a {\em split-level surface}. 

If $\eta_i$ and $\eta_{i+1}$ differ
by a move in a $\xi=4$ geodesic $k$, then there is a block $B_i$ with
domain $D(k)$, so that $F_{i+1}$ is obtained from $F_i$ by removing 
$\boundary_-B_i$ and replacing it with $\boundary_+B_i$. 

Now define $M_0^0 = F_0$, and inductively define 
$$M_i^{j+1} = M_i^j \union B_j$$
and
$$M_{i-1}^j = M_i^j \union B_{i-1}.$$
Thus we are building up $\modl[0]$ by successively adding blocks above and
below.  Since the resolution $\{\tau_i\}$ contains an elementary move
for every 4-edge of $H$ by Lemma \ref{Resolution sweep},
every internal block is included in $\{M_i^j\}$.  

Now the map $\Psi$ can easily be defined
inductively on $M_{-\infty}^\infty = \union_{i,j\in\Z} M_i^j$. 
Begin by mapping $F_0=M_0^0$ to $S\times \{0\}$ by the map that
restricts to the identity on the surface factors. Now suppose we wish
to add a block $B_{j+1}$ to $M_i^j$. By induction
the boundary components of
$\boundary_-B_{j+1}$, which are part of $F_j$, are already mapped
flatly by $\Psi$. 
We  extend $\Psi$ to $B_{j+1}$ so that it is a orientation-preserving flat embedding.
Since there may be two components of $\boundary_-B_{j+1}$ which are mapped
to different heights, we may have to stretch the two ``legs'' of
$B_{j+1}$ by different factors, but the map can be piecewise-affine on
the vertical directions, and the identity on the surface factors.
(See Figure \ref{S05model} for a schematic example).
We always choose $\Psi$ on $B$ so that the image of every annulus in
$\boundvert B$ has height at least 1. This will guarantee that
$\Psi$ is proper.

To define the map on the boundary blocks, note that the gluing
boundary of a boundary block must be part of 
the boundary of
$M^\infty_{-\infty}$, and is mapped by $\Psi$ so that each component
is mapped flatly. Hence we can extend $\Psi$ as an
orientation-perserving flat embedding on each boundary block, making
sure that the map on the added annuli $\boundary R\times[0,\infty)$
(or $\boundary R\times(-\infty,0]$) is proper.

We should now verify that $\Psi$ is an embedding, and
keep track along the
way of how the solid tori $\{\storus v\}$ arise in the complement of
$\Psi(\modl[0])$. 

\begin{lemma}{torus of v}
Let $v$ be vertex of $H$ and let $\FF(v)$ be the union of annuli in
the boundaries of blocks of $\modl[0]$ that are in the homotopy class
of $v$. Then $\FF(v)$ is a torus or an annulus, and $\Psi|_{\FF(v)}$
is an embedding with image
$$
\boundary(\collar(v)\times[s_1,s_2])
$$ 
if $v$ is not parabolic in either $\I(H)$ or $\T(H)$,
$$
\boundary(\collar(v)\times[s_1,\infty))
$$
if $v$ is parabolic in $\T(H)$, and
$$
\boundary(\collar(v)\times(-\infty,s_2])
$$
if $v$ is parabolic in $\I(H)$.
\end{lemma}

\begin{proof}
If $v$ intersects  $\base(\I(H))$, then Lemma \ref{Vertex
  configurations} gives us a unique 4-edge $e_1$ with $v=e_1^+$. Thus
  the block $B(e_1)$ has inner annulus $\boundary^+_i B(e_1)$ in the
  homotopy class of $v$.  Let us call the 
  horizontal subannulus $\boundary^+_{ih}B(e_1)$ the {\em bottom
  annulus of $v$}.
  By the definition, $\Psi$ maps
  this  annulus to $\collar(v)\times\{s_1\}$ for some $s_1(v)\in\R$,
  and the vertical part $\boundary^+_{iv} B(e_1)$ is mapped to annuli
  of the form $\alpha\times[s_1,t]$ where $\alpha$ is a boundary
  component of $\collar(v)$ and $t>s_1$.  

The elementary move associated to $e_1$ introduces $v$ into the
resolution of the hierarchy, so that all other annuli in this homotopy
class occur after $e_1$. 

If $v$ is a component of $\base(\I(H))$, and if it has a
  transversal, then a similar description holds, where $B(e_1)$ is
  replaced by a bottom boundary block.  If $v\in\base(\I(H))$ but has
  no transversal, it is parabolic in $\I(H)$, and there is no block $B$ with an
  inner annulus homotopic to $v$. 

The same discussion holds with regard to $\T(H)$, yielding a unique
{\em top annulus} for $v$ in the bottom boundary
of an appropriate block, unless $v$ is parabolic in $\T(H)$.
We note that $v$ cannot be parabolic in both $\I(H)$ and
$\T(H)$. All other annuli in the homotopy class of $v$ must occur as
outer annuli, i.e. in the sides of blocks whose domains have boundary components
homotopic to $v$. These blocks occur in the resolution in the interval
$J(v)$, that is between the move that introduces $v$ in the
resolution and the one that takes it out. Since $\Psi$ is locally an
embedding (by the orientation-preserving condition) these annuli must
fit together into a torus or annulus, as described in the statement of
the lemma. 
\end{proof}

We note also that these annuli and tori, together with the outer
boundaries of boundary blocks, form the entire boundary of $\modl[0]$.
From the description of $\Psi(\FF(v))$ it clearly bounds an open solid torus
in $S\times\R$, namely $\collar(v)\times(s_1,s_2)$,
$\collar(v)\times(s_1,\infty)$, or
$\collar(v)\times(-\infty,s_2)$ in the various cases, and we denote
this solid torus $U(v)$.

Define the ``top''  $top(M_i^j)$ to be the union of $F_j$ with the
bottom annulus of $v$ for each vertex $v$ in $\eta_j$ (these are the annuli that
$\Psi$ maps to $\collar(v)\times\{s_1(v)\}$, as in Lemma \ref{torus of
  v}), provided this annulus is in $M_i^j$.
We similarly define the ``bottom''  $bot(M_i^j)$ using 
$F_i$ and the corresponding top annuli, if any. Note that
$M_0^0$ has both top and bottom just $F_0$, with no added annuli, but
the annuli appear when we add new blocks.

We can now establish the inductive hypothesis that the
$\Psi|_{M_i^j}$ is an embedding, and its
image $\Psi(M_i^j)$ lies below the image of $top(M_i^j)$ (and a
similar statement for the bottom, which we omit). That is,
for each $(p,t)\in \Psi(top(M_i^j))$, the ray $p\times(t,\infty)$ is
disjoint from $\Psi(M_i^j)$. Thus, since the
block we add to obtain $M_i^{j+1}$ is mapped by an orientation
preserving flat embedding, it is immediate that its image lies above
$\Psi(top(M_i^j))$, and the inductive property is preserved. 
The same argument
applies to adding a block from below to get $M_{i-1}^j$.

It follows that $\Psi$ is a proper embedding on
$M_{-\infty}^\infty$. The boundary blocks can then be treated in the
same way, noting that the gluing boundary of a top boundary block is
eventually in $top(M_i^j)$ for large enough $j$ (and similarly for
bottom boundary blocks). Thus $\Psi$ is an embedding on all of $\modl[0]$.

Finally, we note that each solid torus $U(v)$ is in fact disjoint from
$\Psi(\modl[0])$: For $j=\inf J'(v)$, the block $B_{j-1}$
contains the bottom annulus of $v$, and hence this annulus is in
$top(M_i^j)$.  Thus
$\Psi(M_i^j)$ is below the image annulus
$\collar(v)\times\{s_1(v)\}$, and in particular disjoint from $U(v)$
which lies above this annulus. 
For $\inf J'(v) <  j \le \sup J'(v)$, $B_{j-1}$
is disjoint from $\collar(v)\times\R$ and in particular from $U(v)$. 
For all $j>\sup J'(v)$, $\Psi(top(M_i^j)$ no longer contains
$\collar(v)\times\{s_1(v)\}$ and all subsequent block are placed above
the image of the top annulus, 
$\collar(v)\times\{s_2(v)\}$, and hence disjoint from $U(v)$.

\realfig{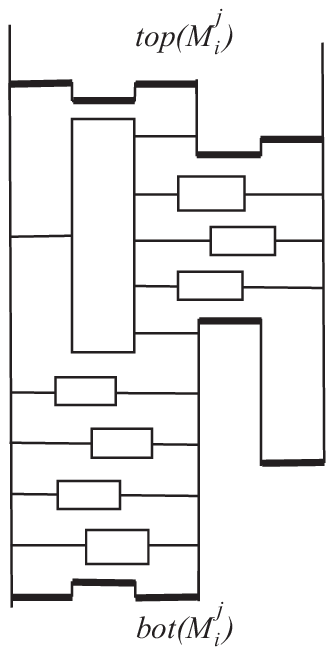}{A schematic of a flat
  embedding of $M_i^j$ into $S\times\R$. Note that some portions 
  of blocks are stretched vertically. The top and bottom are indicated
  with thickened lines. This picture is a fairly accurate rendition of
  the case $S=S_{0,5}$.}

\subsubsection*{Peripheral tubes}
It remains to describe the tubes $U(v)$ for $v$ a component of
$\boundary S$. So far, both $\modl[0]$  and the tubes $U(v)$ for
vertices $v$ have embedded, disjointly, into $S\times \R$. Now for a boundary
component $v$, the annulus $\collar(v)$ is a
component of $\hhat S \setminus S$. We let $U(v)$ be simply 
$\collar(v) \times \R$. 

The intersection of $\boundary U(v)$ with $\Psi(M_\nu)$ is the union
of all vertical annuli of blocks which are in the homotopy class of
$v$. As we observed above for vertices, these annuli cover
all of $\boundary U(v)$.
\end{proof}

Let $\UU$ be the union of the open tubes $U(v)$ we have described, so
that $\UU$ is disjoint from $\Psi(\modl[0])$ and 
$\boundary\UU$ is contained in $\boundary\Psi(\modl[0])$. From now on
we identify $\modl[0]$ with its $\Psi$-image, and define 
$\modl = \modl[0] \union\UU$.

\subsection{The model metric}
\label{define metric}

We will describe a metric for each of the finitely many block types of
which $\modl[0]$ is constructed, and use this to piece together a 
metric on all of $\modl[0]$. Then we will discuss the structure of the 
tubes $\UU$, and extend the metric to them as well. 

\subsubsection*{Internal blocks} 
Fix one copy $W_1$ of a four-holed sphere, and one copy $W_2$ of a
one-holed torus. Mark $W_1$ with a pair $v^-_1,v^+_1$ of adjacent vertices
in $\CC(W_1)$, and similarly $v^-_2,v^+_2$ for $\CC(W_2)$. This
determines two blocks $\hhat B_k = B([v^-_kv^+_k])$, $k=1,2$. That is,
$\hhat B_k$ is constructed from $W_k\times[-1,1]$ by removing solid torus
neighborhoods of $\gamma_{v^-_k}\times\{-1\}$ and
$\gamma_{v^+_k}\times\{1\}$, as in the construction at the 
beginning of the section. 

Let us define some {\em standard metrics} on the surfaces
 $S_{0,2}$, $S_{0,3}$, $S_{0,4}$ and $S_{1,1}$: 

Call an annulus {\em standard} if it is isometric to $\bbar\collar(\gamma)$
for a geodesic $\gamma$ of length $\ep_1/2$. 
Fix a three-holed sphere $Y'$ with a hyperbolic metric so that
$\boundary Y'$ is geodesic with components of length $\ep_1/2$. 
call a three-holed sphere {\em standard} if it is isometric to 
$Y'\setminus \collar(\boundary Y')$.

Now fix a surface $W'_k$ homeomorphic to $W_k$, and
endowed with a fixed hyperbolic metric $\sigma$ for which $\boundary W'_k$ is 
geodesic with components of length $\ep_1/2$, 
and fix an identification of
$W_k$ with $W'_k \setminus \collar(\boundary W'_k)$. We can do this, for
specificity, in such a way that the 
curves $\gamma_{v^\pm_k}$ are identified with a
pair of orthogonal geodesics of equal length. 
Finally, choose the identification so that the collars
$\collar(v_k^\pm)$ in $W_k$ (fixed by our global convention in
\S\ref{isotopy convention})
are identified with $\collar(v_k^\pm,\sigma)$ in $W'_k$.
We call this metric  {\em standard} on $W_k$. 

Now we may fix a metric on $\hhat B_k$ with the following properties: 
\begin{enumerate}
\item  The metric restricts to standard metrics on
 $W_k\times\{0\}$, and on each 3-holed sphere in 
$\boundary_{\pm}\hhat B_k$.
\item  Each annulus component of $\boundvert\hhat B_k$
(respectively $\boundary_i \hhat B_k$)
is isometric to
$S^1\times[0,\ep_1]$ (resp. $S^1\times[0,\ep_1/2]$), 
with $S^1$ normalized to length $\ep_1$,
and this product structure agrees with
the product structure imposed by the inclusion in $W_k\times[0,1]$.
\end{enumerate}
(The details of the construction do not matter, just that these
properties hold and that a fixed choice is made. The length of $\ep_1$ for 
the cores of all the Euclidean annuli is made possible by the
definition of collars in \S\ref{collar defs}.)

Note that, by definition of a standard metric, each component of the
gluing boundary $\boundary_\pm\hhat B_k$ admits an
orientation-preserving isometry group
realizing all six permutations of the three boundary components. This will
enable us to glue the blocks via isometries. 

We'll call these two specific blocks the ``standard blocks.'' Every
block $B(e)$ in $\modl[0]$ 
associated to a 4-edge $e$ can be identified with one of $\hhat B_1$ or
$\hhat B_2$, depending on the homeomorphism type of $D(e)$, by a map
that takes $e_\pm$ to $v^\pm_k$. This identification is unique up to
isotopy preserving the various parts of the boundary, and
any such identification yields what we call a standard metric on $B(e)$.

\subsubsection*{Boundary blocks} 

Consider a top boundary block $B=B\sbtop(\nu_R)$.
Recall that the outer boundary $\boundary_o B$
was constructed as $R\times\{0\} \union \boundary
R\times[0,\infty)$, which is homeomorphic to $int(R)$. 
Endow it with the Poincar\'e metric $\sigma_\infty$ representing
$\nu_R$, in such a way that 
$\collar(\boundary R,\sigma_\infty)$ is identified with
$\boundary R\times(0,\infty)$, and
$\collar(\T_R,\sigma_\infty)$ is identified with
$\collar(\T_R)\times\{0\}$. 
Let $\sigma^m$ be the conformal rescaling described in
\S\ref{augmented core}, which makes the collars of curves of length
less than $\ep_1$ into Euclidean cylinders. Note that (by definition
of $\T_R$) these collars
are either components of $\collar(\T_R)$ or of $\collar(\boundary R)$.

Let $\sigma'_\infty$ be a hyperbolic metric on $int(R)$ for which
{\em every} component of $\T_R$ has length less than $\ep_1/2$, and
which differs from $\sigma_\infty$ by a uniformly bilipschitz
distortion (the constant can be chosen to depend only on $\ep_1/L_0$).
Let $\sigma^{m'}$ be the conformal
rescaling of $\sigma'_\infty$ that makes component of 
$\collar(\T_R\union \boundary R)$ a Euclidean cylinder, and equals
$\sigma'_\infty$ elsewhere. Then $\sigma^{m'}/\sigma^m$ is uniformly
bounded above and below.

Now we can transport $\sigma^{m'}|_{R\times\{0\}}$
to $R\times\{-1/2\}$ via the identity on the first factor. 
Extend these to a metric on
on $R\times[-1/2,0]$ which is in uniformly bounded ratio with
the product metric of $\sigma^m$  and 
$dt$ (here $t\in [-1,0]$). The rest of the block, $(R\setminus
\collar(\T_R)) \times 
[-1,-1/2]$, may be metrized as follows: 
The restriction of $\sigma^{m'}$ to each component $Y$ of
$R\setminus\collar(\T_R)$ 
is uniformly bilipschitz equivalent to a standard metric on a
three-holed sphere, since all the curves of $\T_R$ have
$\sigma_\infty$-length at most $L_0$ by the choice of $\T_R$
(see \S\ref{define Hnu}).
Thus we may place a standard metric on 
the corresponding gluing surface $Y\times\{-1\}$, and interpolate
between them on $Y\times[-1,-1/2]$ so that the resulting metric is
uniformly bilipschitz equivalent to the product metric of $\sigma_\infty|_Y$
and $dt$. 

We may do this, as for the internal blocks, in such a way that the
annulus boundaries are given a flat Euclidean structure, and the
natural product structure agrees with that inherited from $S\times\R$.

The metric on bottom boundary blocks is defined analogously.

\begin{lemma}{Build metric}
There is a metric on $\modl[0]$ which restricts on each block to a
standard metric. 
\end{lemma}

\begin{proof}
The main issue is to choose the identifications of each block with the
standard blocks in a consistent way. We will use the exhaustion of
$\modl[0]$ by subsets $M_i^j$, and the fact that all the gluing boundaries
of the standard blocks are isometric to a specific three-holed sphere
$Y_0$, which itself admits all possible symmetries.

Beginning with $M_0^0$, which is just the
split-level surface $F_0$, choose an identification of each component
$Y$ with $Y_0$. Now inductively suppose we have a metric on $M_i^j$
which is standard on each block. Then each gluing boundary component
$Y$ has been identified to $Y_0$ in some way. To obtain $M_i^{j+1}$ we
adjoin a block $B$ by attaching its lower gluing boundary
$\boundary_-B$ to one or two components of the upper  gluing boundary
of $M_i^j$. Consistency of this gluing determines the metric on
$\boundary_- B$. Now as above, we have a unique equivalence class of
identifications of $B$ with the appropriate standard block $\hat B_k$
($k=1$ or $2$), which determines an isotopy class of identifications
of $\boundary_- B$ with $\boundary_- \hat B_k$.
the metric we already have on $\boundary_-B$ determines a unique
representative of this isotopy class. 

We therefore only have to extend the map on $\boundary_- B$ to all of
$B$. This is an elementary consequence of the homotopy extension
theorem. In fact we can do this in such a way that 
map on the annuli $\boundvert B$ preserves the natural product
structure.

We argue similarly in the downward direction, adding a block from
below to get $M_{i-1}^j$. Thus by induction we have a metric on all of
$M_{-\infty}^\infty$. A similar argument for  the boundary blocks
extends the metric to all of $\modl[0]$.
\end{proof}

\subsection*{Meridian coefficients}

Let $v$ be a non-parabolic vertex in $H$. The torus $\boundary U(v)$
inherits a Euclidean metric from $\modl$,  and a boundary orientation from
$\storus v$. It also has a natural marking $(\alpha,\mu)$,
where $\alpha$ is the homotopy class of 
the cores of the annuli making up $\boundary\storus v$,
homotopic to $\gamma_v$ in $S$, and $\mu$
is the meridian class of $\boundary\storus v$.
To describe the meridian
explicitly, recall that we represent $\storus v$ as a product
$\collar(v)\times (s,t)$. If $a$ is any simple arc in
$\bcollar(v)$ connecting the boundaries then
\begin{equation}\label{meridian def}
 \boundary (a\times[s,t])
\end{equation}
is a meridian. 
Note that the choice of $a$ does not affect the isotopy class of this curve
in $\boundary\storus v$.

As in \S\ref{thick thin}, we can describe the geometry of this
oriented marked Euclidean torus with parameters $(\omega,t)$. In this
case $t=\ep_1$, since the circumference of the annuli in $\boundary\storus v$
is $\ep_1$ 
by construction. The Teichm\"uller parameter $\omega\in\Hyp^2$
will be called the meridian coefficient of $v$, and denoted
$\Momega(v)$.

Note that $\ep_1|\Momega(v)|$ is the length of the meridian, and also
that the imaginary part $\ep_1\Im \omega$ is simply the sum of the heights of
the annuli that make up $\boundary\storus v$. We have $\Im \omega \ge 1$ since
$\boundary\storus v$ contains at least the annuli from its bottom and top blocks. 

If $U(v)$ is a parabolic tube, we  define $\omega_M = i\infty$.

\subsection*{Metrizing the tubes}
For each non-parabolic tube $\storus v$, 
Lemma \ref{torus boundary data} gives us a unique
hyperbolic tube $\MT(\lambda,r)$ whose boundary parameters with
respect to the natural marking, are
$(\Momega(v),\ep_1)$, so we identify this tube with $\storus v$ via a 
marking preserving isometry on the boundary. There is clearly a unique
way to do this up to isotopy of $\storus v$.

$\boundary U(v)$ for a parabolic $v$ or a boundary component of $S$
is an infinite Euclidean cylinder of
circumference $\ep_1$, and there is a unique (up to isometry) rank-1 parabolic
tube with circumference $\ep_1$, so we impose that metric on $U(v)$.

This completes the definition of the model metric on all of $\modl$.

We remark that our eventual goal (in \cite{brock-canary-minsky:ELCII})
is to show that these hyperbolic tubes are in fact bilipschitz
equivalent to the tubes in the hyperbolic manifold $N_\rho$. In this
paper we will only obtain a Lipschitz map, and only for those tubes
with $|\Momega|$ bounded above by a certain constant.

\subsubsection*{Fillings of $\modl[0]$}
As in the introduction, we define
$$
\UU[k] = \bigcup_{|\Momega(v)|\ge k} \storus v
$$
and
$$
\modl[k] \equiv \modl[0] \union \bigcup_{|\Momega(v)|<k} \storus v.
$$
Note that  $\modl[\infty]$ indicates the inclusion of all the
non-parabolic tubes.

\section{Comparing meridian coefficients}
\label{count}

We introduced the meridian coefficient $\Momega(v)$ in Section
\ref{model}, to describe the geometry of the model torus associated to
a vertex $v$ of a hierarchy $H$. Now we will describe two more ways of
estimating this coefficient: one, $\Homega(v)$, will be defined using
the data in the hierarchy $H_\nu$, and the third, $\nomega(v)$,
using the subsurface projection maps $\pi_W$ applied directly to the
end invariants $\nu$.

We will prove that these invariants are close in the following sense:

\begin{theorem}{Omegas close}
  For any pair of end invariants $\nu$ and associated hierarchy $H$
  and model manifold $M$, and for any non-parabolic vertex $v$ in
  $H$, we have bounds
  \begin{equation}
    \label{HM bound}
    d_{\Hyp^2}(\Homega(v),\Momega(v)) \le D
  \end{equation}
and
\begin{equation}
  \label{Hn bound}
  d_{\Hyp^2}(\Homega(v),\nomega(v)) \le D
\end{equation}
where $D$ depends only on the topological type
of $S$. Here $d_{\Hyp^2}$ refers to the Poincar\'e metric in the upper
half plane.
\end{theorem}

Note that we are interpreting the $\omega$'s as
Teichm\"uller parameters, so this estimate is natural since
$d_{\Hyp^2}$ can be identified with the Teichm\"uller distance
in  Teichm\"uller space of the torus.

If $v$ is a vertex of $\CC(S)$, let us denote the two boundary
components of $\collar(v)$ (arbitrarily) as $v\sll$ and $v\srr$.
Now for $\alpha = v\srr$ or $v\sll$, define
\begin{equation}\label{Xalpha def}
X_{\alpha} = \{h\in H: \alpha\subset \boundary D(h)\}.
\end{equation}
Note that the same $h$ can be in both $X_{v\sll}$ and $X_{v\srr}$, if
its domain borders $\collar(v)$ from both sides. We similarly define
\begin{equation}\label{Xalpha4 def}
X_{\alpha,k} = \{h\in H: h\in X_\alpha, \xi(D(h))=k\}
\end{equation}

and
\begin{equation}\label{Xalpha4+ def}
X_{\alpha,k+} = \{h\in H: h\in X_\alpha, \xi(D(h))\ge k\}
\end{equation}

\subsubsection*{Coefficients for internal vertices}
From now until \S\ref{noninternal omega}, we will
assume that $v$ is an {\em internal} vertex of $H$ --
that is, it is not a vertex of $\simp(\I(H))$ or $\simp(\T(H))$.
Let $h_v$ be the annulus geodesic in $H$ with domain $\collar(v)$.
We can then define:

\begin{equation}\label{internal Homega}
\Homega(v) = [h_v] + i\biggl(1 + 
\sum_{\alpha=v\sll,v\srr}\sum_{h\in X_{\alpha,4+}} |h| \biggr).
\end{equation}
(Here $[h_v]$ is the 
{\em signed} length of $h_v$, defined as in (\ref{signed length})).

For the next definition, again with $\alpha=v\sll$ or $v\srr$, let
\begin{equation}\label{Yalpha def}
\YY_\alpha = \{Y\subset S: \alpha\subset \boundary Y\}
\end{equation}
(where we recall our convention of only taking the standard representatives
of isotopy classes of surfaces, whose boundaries are collar
boundaries.)
Define also  $\YY_{\alpha,4}$ and $\YY_{\alpha,4+}$ in analogy with
(\ref{Xalpha4 def}) and (\ref{Xalpha4+ def}).

We now define
\begin{equation}
  \label{internal nomega}
\nomega(v) = \tw_v(\nu_-,\nu_+) + i\biggl(1 +
  \sum_{\alpha=v\sll,v\srr}\sum_{Y\in\YY_{\alpha,4+}} 
                                    \Tsh K{d_Y(\nu_-,\nu_+)}
    \biggr)
\end{equation}
where $K$ will be determined later, and $\Tsh K{x}$ is the threshold
function, defined in \S\ref{del estimates}. Here for convenience we
have defined $d_Y(\nu_-,\nu_+) \equiv d_Y(\mu_-,\mu_+)$ and similarly
for $\tw_Y$, where $\mu_\pm$ are the generalized markings derived from
$\nu_\pm$ in \S\ref{define Hnu}.

\subsection{Shearing}

Let us recall from \cite{minsky:boundgeom} the notion of {\em shearing
outside a collar} for two hyperbolic metrics on a surface, and some of
its properties. Let
$\gamma\in\CC_0(S)$, and suppose that $\sigma_1$ and $\sigma_2$ are
two hyperbolic metrics on $S$ for which 
$\collar(\gamma,\sigma_1) =\collar(\gamma,\sigma_2)$ (therefore
denote both as $\collar(\gamma)$). 
We emphasize that $\sigma_i$ are {\em actual metrics} rather than
isotopy classes. 

We define a quantity
$$
\shear_\gamma(\sigma_1,\sigma_2) 
$$
as follows. Let 
$\hhat Y$ denote the compactified annular cover of $S$ associated to
$\gamma$ (as in Section \ref{complexes}), let
$\hhat B$ be the annular lift of $\collar(\gamma)$ to this cover, 
and let $\hhat\sigma_1$ and $\hhat\sigma_2$ denote the lifts of
$\sigma_1$ and $\sigma_2$ to $int(\hhat Y)$.
Define $\shear_\gamma(\sigma_1,\sigma_2) $ to be
$$
\sup_{E,g_1,g_2}
d_{\AAA(E)}(g_1\intersect 
E,g_2\intersect E).
$$
Here $E$ varies over the two 
complementary annuli $E_1,E_2$ of $\hhat B$ in $\hhat Y$, and
$g_i$ (for $i=1,2$) varies 
over all arcs in $\hhat Y$ connecting the two boundaries
which are $\hhat\sigma_i$-geodesic in $int(\hhat Y)$.
Thus we are measuring the relative twisting, outside $\hhat B$, of any
two geodesics in the two metrics.

A bound on the shear allows us to measure twisting by restricting to a
collar.  If $\alpha_i$ are simple closed $\sigma_i$-geodesics in $S$
crossing $\gamma$ (for $i=1,2$)
and $a_i$ are components of $\alpha_i\intersect
\bcollar(\gamma)$ that cross $\collar(\gamma)$
then  we will show
\begin{equation}\label{twists in collar}
|\twist_{\bcollar(\gamma)}(a_1,a_2)
- \twist_{\gamma}(\alpha_1,\alpha_2)| 
\le 2\shear_\gamma(\sigma_1,\sigma_2) + 2
\end{equation}
where we recall from \S\ref{complexes} that
$\twist_{\bcollar(\gamma)}(a_1,a_2)$ is a 
quantity that depends on the exact arcs of intersection with the
collar, whereas $\twist_\gamma(\alpha_1,\alpha_2)$ depends only on the
homotopy classes of $\alpha_1$ and $\alpha_2$.

\begin{proof}[Proof of (\ref{twists in collar})]
Let $a'_i$ for $i=1,2$ be the $\sigma_i$-geodesic arc in $\hhat Y$
that represents $\pi_Y(\alpha_i)$ -- 
that is, a choice of (geodesic) lift of $\alpha_i$ to the annulus that
connects the boundaries. By definition, 
$$
\twist_\gamma(\alpha_1,\alpha_2) = \twist_{\hhat Y}(a'_1,a'_2).
$$
The right hand side decomposes as the sum
$$
\twist_{E_1}(a'_1\intersect E_1,a'_2\intersect
E_1)
+
\twist_{\bbar B}(a'_1\intersect \bbar B,a'_2\intersect
\bbar B)
+
\twist_{E_2}(a'_1\intersect E_1,a'_2\intersect
E_2),
$$
by two applications of 
the additivity property (\ref{additivity annuli})
(where $\bbar B$ is the closure of $\hhat B$, i.e. the lift of
$\bcollar(\gamma)$ to $\hhat Y$).
The absolute values of the first and third terms 
are bounded by $\shear_\gamma(\sigma_1,\sigma_2)$, by the definition and 
inequality (\ref{twist n dist}). 

The arcs $a_1$ and $a_2$ in $\bcollar(\gamma)$ lift to two arcs in
$\bbar B$, which we still call $a_1$ and $a_2$, and which are disjoint
from or equal to $a'_1\intersect \bbar B$ and $a'_2\intersect \bbar B$,
respectively. Thus $\twist_{\bbar B}(a_i,a'_i\intersect \bbar B) = 0$, 
and using the additivity inequality (\ref{additivity arcs}) twice, we
conclude
$$
|\twist_{\bcollar(\gamma)}(a_1,a_2) - \twist_{\bbar B}(a'_1\intersect
\bbar B, a'_2\intersect \bbar B)| \le 2.
$$
The estimate (\ref{twists in collar}) follows.
\end{proof}

We can bound the shear between two metrics in the following setting:

\begin{lemma}{BG shear bound}{\rm (\cite{minsky:boundgeom})}
Suppose $R$ is a subsurface of $S$ which is convex in two hyperbolic metrics
$\sigma$ and $\tau$, and that $\sigma$ and $\tau$ are locally
$K$-bilipschitz in the complement of $R$. Suppose that one component
of $R$ is an annulus $B$ which is equal to both
$\collar(\gamma,\sigma)$ and $\collar(\gamma,\tau)$ for a certain
curve $\gamma$.  Then
$$
\shear_\gamma(\sigma,\tau) \le \delta_0 K,
$$
where $\delta_0$ depends only on the topological type of $S$.
\end{lemma}
We remark that this lemma in \cite{minsky:boundgeom} is stated for a
slightly different definition of $\collar(\gamma,\sigma)$. However,
the same argument applies for our definition, resulting in a different
constant $\delta_0$.

Another property of shearing that follows immediately from the 
definition and the triangle inequality
is subadditivity: If $\gamma$ has the same collar with
respect to three metrics $\sigma_1,\sigma_2$ and $\sigma_3$, then
\begin{equation}\label{shear subadditive}
\shear_\gamma(\sigma_1,\sigma_3) \le 
\shear_\gamma(\sigma_1,\sigma_2) + 
\shear_\gamma(\sigma_2,\sigma_3).
\end{equation}

\subsection{Sweeping through the model}
\label{sweeping}

Let us elaborate on the sequence of split-level surfaces $\{F_i\}_{i\in\II'}$
constructed in the proof of Theorem \ref{embed model}. Each $F_i$
corresponds via the projection $S\times\R \to S$ to a subsurface
$Z_i\subset S$ which is the complement of a union of annuli
$\collar(\eta_i)$, where $\eta_i$ is a pants decomposition coming from
the resolution of the hierarchy. 

We can add ``middle'' surfaces $F_{i+1/2}$ as follows: Each transition
$F_i \to F_{i+1}$ corresponds to an edge $e_i$ in a 4-geodesic, and
hence a block $B_i = B(e_i)$, so that $F_i$ and $F_{i+1}$ agree except
on $F_i\intersect B_i = \boundary_-B_i$ and
$F_{i+1}\intersect B_i = \boundary_+B_i$. Define the surface
$F_{i+1/2}$ to be $F_i\intersect F_{i+1}$ union the middle surface
$D(e_i)\times\{0\}$. The projection of this to $S$, which we call
$Z_{i+1/2}$, is the complement of  $\collar(\eta_{i+1/2})$, where 
$\eta_{i+1/2}$ is $\eta_i\intersect\eta_{i+1}$. 

Now for $s$ either integer or half-integer,
Let $\sigma'_s$ be the metric on $Z_s$ inherited from $F_s$. By
construction, this metric is ``standard'' on each component (see
Section \ref{define metric}) and hence
can be extended to a hyperbolic metric $\sigma_s$ on
all of $S$, such that the components of $\collar(\eta_s)$ are 
standard collars with respect to $\sigma_s$. Note that $\sigma_s$
is not unique but we will fix some choice for each $s$.

Now consider the interval $J'(v)\subset\II'$, consisting of those $i$
for which $v$ is a component of $\eta_i$. 
Since $v$ is internal, $t=\max J'(v)< \sup\II'$ and $b=\min J'(v) > \inf
\II'$. Thus there is a transition $\eta_{b-1}\to \eta_b$ that replaces
some vertex $u$ by $v$, and 
a transition $\eta_t \to \eta_{t+1}$ replacing $v$ by some $w$. 
The block $B_{b-1} = B(e_{b-1})$ contains the bottom annulus of
$\boundary \storus v$, 
and the block $B_t = B(e_t)$ contains the top annulus of
$\boundary\storus v$.

In order to understand how $\Momega(v)$ is related to $\Homega(v)$, we
will have to understand something about the relation between the metrics
$\sigma_{b-1/2}$ and $\sigma_{t+1/2}$:

\newcommand\HALF{1/2}
\begin{lemma}{Shear bound}
$$
\shear_{v}(\sigma_{b-\HALF}, \sigma_{t+\HALF}) 
\le K \Im \Momega(v)
$$
where $K$ depends only on the topological type
of $S$.
\end{lemma}

\begin{proof}
We will apply Lemma \ref{BG shear bound} to the transitions
$\sigma_s \to \sigma_{s+\HALF}$, where $s$ (for the rest of the proof)
is an integer or half-integer in $[b-\HALF,t]$.

For every integer $s\in[b,t]$,  $Z_{s+\HALF}$ has a component $W(e_s)$ of type
$\xi=4$ which is the union of one or two components of $Z_s$ and a
collar. The metrics $\sigma_s$ and
$\sigma_{s+\HALF}$ agree pointwise on $Z_s\setminus W(e_s) = 
Z_{s+\HALF}\setminus W(e_s)$. 
On each component $Y$ of $Z_s$ which is contained in $W(e_s)$, the metrics
$\sigma_s$ and $\sigma_{s+\HALF}$ are related by a uniform bilipschitz
constant $L$, which only depends on our model construction. This is
because these two metrics on $Y$ come from subsurfaces of one of the
blocks, and the identification is inherited from one of our finite
number of standard blocks, where there is {\em some} bilischitz constant.
The corresponding statement holds for $s$ and $s-\HALF$. 

For each $s\in[b,t]$, let $V_s$ denote the union of components  (one or two)
of $Z_s$ that are adjacent to $\collar(v)$. 
Let $Q$ denote the union of intervals $[s,s+\HALF]$ in $[b,t]$ for
which $V_s = V_{s+\HALF}$. 

Let $[x,y]$ be a component of $Q$. Then $V_x = V_y$. Let $R$ be the
subsurface $S\setminus V_x$. Then $R$ contains $\collar(v)$ as a
component, and is convex in both metrics $\sigma_x$ and $\sigma_y$, 
since each boundary component $\gamma$ of $R$ is a boundary component
of the standard $\collar(\gamma)$ which is contained in $R$. 
Since every transition $s\to s+\HALF$ in $[x,y]$ does not involve
$V_s=V_x$, we conclude that $\sigma_x$ and $\sigma_y$ are pointwise
identical in the complement of $R$. Thus we may apply Lemma \ref{BG
shear bound} to conclude that
\begin{equation}
\label{Q shear bound}
\shear_v(\sigma_x,\sigma_y) \le \delta_0.
\end{equation}

Let $s$ be such that $(s,s+\HALF)$ is not in $Q$, and suppose first that
$s$ is an integer. Let $R$ be the complement of
$V_s$. This is still convex in both metrics $\sigma_s$ and $\sigma_{s+\HALF}$
just as in the previous paragraph, and has
$\collar(v)$ as a component. This time the transition does involve
$V_s$, so at least one component of $V_s$ is contained in the
component $W(e_s)$ of $Z_{s+\HALF}$, and 
the two metrics on this component are related by a uniform bilipschitz
constant $L$. 
On all other components they are identical, 
so we have a bilipschitz bound of  $L$ on the complement of $R$. Again
applying Lemma \ref{BG shear bound} we obtain
\begin{equation}
\label{transition shear bound}
\shear_v(\sigma_s,\sigma_{s+\HALF}) \le \delta_0L.
\end{equation}
The same bound holds if $s$ is a half-integer,
since then $s+\HALF$ is an integer and we can apply the same argument in
the opposite direction.

Thus, applying the subadditivity property (\ref{shear subadditive}),
we conclude that 
\begin{equation}
\shear_v(\sigma_{b-\HALF},\sigma_{t+\HALF}) \le \delta_0\#Q + \delta_0 L N
\end{equation}
where $\#Q$ is the number of components of $Q$, 
and $N$ is the number of intervals $(s,s+\HALF)$ outside of $Q$.
We note now that every such $(s,s+\HALF)$ with $s\ne b-\HALF, t$
corresponds to half of a block 
$B$ with the property that $\boundvert B$ contains an annulus on
$\boundary \storus v$ (in other words the domain subsurface of $B$ is
adjacent to $\collar(v)$). The intervals $(b-\HALF,b)$ and $(t,t+\HALF)$
correspond to the blocks meeting $\storus v$ on the bottom and top,
respectively, in annuli of height $\ep_1/2$. 
The sum of heights of all annuli of $\storus v$ is exactly $\ep_1\Im
\Momega(v)$, and hence we obtain $N \le 2\ep_1\Im \Momega(v)$. 

It is evident that $\#Q \le N-1$, and we conclude
\begin{equation}
\shear_v(\sigma_{b-\HALF},\sigma_{t+\HALF}) \le
2\delta_0(1+L)\ep_1\Im\Momega(v). 
\end{equation}
\end{proof}

\subsection{Comparing $\Homega$ and $\Momega$}
\label{compare H M}

Recall now that $u$ is the predecessor of $v$ on the 4-edge $e_{b-1}$,
and let $u^-$ denote the geodesic representative of $u$
in the metric $\sigma_{b-\HALF}$. Note that the length of $u^-$ in this
metric is uniformly bounded above and below.  Similarly $w$ is the successor of $v$ on $e_{t}$, and
we let $w^+$ be its geodesic representative in the metric
$\sigma_{t+\HALF}$, also with the same length bounds. 

Let $a^-$ be a component (there may be two) of the intersection of
$u^-$ with $\collar(v)$, and let $a^+$ be a component of the
intersection of $w^+$ with $\collar(v)$. Writing $\storus v $ 
as $\collar(v) \times [p,q]$, we can use the curve
$$
\mu = \boundary (a^+ \times [p,q])
$$
as a meridian.
The real part of $\Momega(v)$ is the amount of twisting of $\mu$ 
around the $v$ direction in the torus, and this gives rise to the
following estimate: 

\begin{lemma}{twist is real part}
$$
|\twist_{\bcollar(v)}(a^-,a^+) - \Re \Momega(v)| = O(1).
$$
\end{lemma}

\begin{proof}
Recalling the definition of the marked torus parameters in 
\S\ref{thick thin}, and the discussion in \S\ref{define metric},
we have an orientation-preserving identification of the
torus $\boundary\storus v$ as the quotient of 
$\C/(\Z+\Momega(v)\Z)$, which is an isometry with respect to $1/\ep_1$
times the model metric on $\boundary\storus v$.
After possibly
translating we may assume that the annulus $\collar(v)\times\{p\}$ 
lifts to a horizontal strip $V_1 = 
\{z:\Im z \in[0,k_1]\}$
and $\collar(v)\times\{q\}$ lifts to a horizontal strip $V_2 =
\{z:\Im z \in[k_2,k_3]\}$ with $0<k_1<k_2<k_3 < \Im \Momega(v)$
(See figure \ref{toruslift}).

\realfig{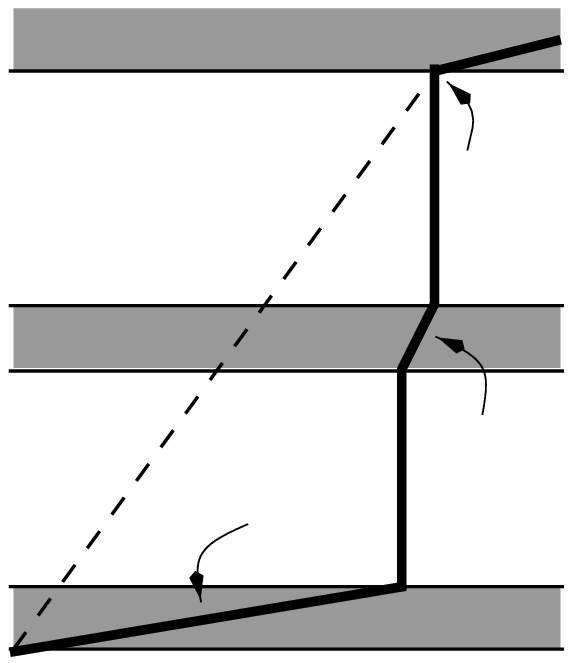}{A representative of the
meridian of $\boundary\storus v$, lifted to $\C$.}

The meridian $\mu$ lifts to four arcs in the cover (and their
translates): $a^+\times\{p\}$ lifts to an arc $a^+_1$ in $V_1$
connecting $0$ to a point in $\R+ik_1$,
$a^+\times\{q\}$ lifts to an arc $a^+_2$ in $V_2$ connecting 
$\R+ik_2$ to $\R+ik_3$,
and the arcs $\boundary a^+ \times[p,q]$ lift
to vertical arcs connecting $a^+_1$ to $a^+_2$, and $a^+_2$ to
$\Momega(v)$,  as shown in Figure \ref{toruslift}.

Thus, it is immediate that $\Re\Momega(v) = b_1 + b_2$, where $b_i$ is
the real part of the vector from the bottom to the top of $a^+_i$. 
Since the length of $a^+$ in the metric $\sigma_{t+1/2}$ is $O(1)$,
the same is true for its length in the annulus
$\collar(v)\times\{q\}$, and hence $b_2 = O(1)$. On the other hand,
$b_1$ may be large because in the lower annulus
$\collar(v)\times\{p\}$, it is $a^-$ that has length $O(1)$, 
and $a^+$ twists around $a^-$ possibly many times. The number of
twists is
determined by $|\twist_{\bcollar(v)}(a^-,a^+)|$, and the sign
convention is such (see Section \ref{complexes}) that 
$b_1 = \twist_{\bcollar(v)}(a^-,a^+) + O(1)$. 
The lemma follows.
\end{proof}

Now we note, applying (\ref{twists in collar}) and Lemma \ref{Shear
bound}, that 
\begin{align}
|\twist_v(u,w) - 
\twist_{\bcollar(v)}(a^-,a^+)| & = 
O(\shear_v(\sigma_{b-\HALF},\sigma_{t+\HALF})) \notag \\ 
& = O(\Im \Momega(v)).
\end{align}

We also have
\begin{equation}
|\twist_v(u,w) - [h_v]| = O(1)
\end{equation}
by the definition of annulus geodesics in hierarchies, and
inequality \ref{twist n dist}.
Recall that $[h_v] = \Re \Homega (v)$.

Combining these with Lemma \ref{twist is real part}, we have
\begin{equation}\label{MH real parts}
|\Re \Momega(v) - \Re \Homega(v)| = O(\Im \Momega(v)).
\end{equation}
(Here we are using $\Im\Momega(v)\ge 1$ to subsume $O(1)$ terms
into the $O(\Im\Momega(v))$ term).

We also have 
\begin{equation}\label{MH imag parts}
\Im \Homega(v) = \Im \Momega(v).
\end{equation}
This is because the sums $\sum_{h \in X_{\alpha,4}} |h|$ in the definition
(\ref{internal Homega}) for $\Homega$ count the number of $4$-edges 
associated to domains bordering $\collar(v)$, counting an edge
once for each side that its domain borders. This therefore gives
exactly the number of vertical annuli in $\boundary\storus v$, each of which
has height $\ep_1$. Counting also the annuli from the bottom and top blocks
for $\storus v$, whose heights add to $\ep_1$, and recalling that
$\ep_1\Im\Momega(v)$ is the sum of these heights,
we have the equality (\ref{MH imag parts}).

It follows immediately from (\ref{MH real parts}) and (\ref{MH imag parts})
that the distance in $\Hyp^2$ between
$\Momega(v)$ and $\Homega(v)$ is bounded.  This gives inequality
(\ref{HM bound}) of Theorem \ref{Omegas close}.

\subsection{Counting in a hierarchy}
\label{counting arguments}

In order to compare $\nomega$ and $\Homega$, 
we must consider more carefully the
structure of a hierarchy, and prove some counting lemmas
that allow us to
estimate the ``size'' of a hierarchy or certain subsets of it by
various approximations. 

Recall that if $f\fsubd g$ then there is a unique simplex $v$ of $g$
for which $D(f)$ is a component domain of $(D(g),v)$. Let us say in this
situation that $f\fsubd
g$ {\em at $v$}. Call $v$ an {\em interior simplex} of $g$ if it is
neither the first or the last. If $v$ is interior then there exactly
one $f\fsubd g$ at $v$, namely the one whose domain contains the
predecessor and successor of $v$. If $v$ is last, then it
is a single vertex (by construction) and there are at most three component
domains of $(D(g),v)$ supporting a geodesic $f\fsubd g$ (including the
annulus $\collar(v)$). If $v$ is the first vertex then there are at
most $2$ such $f$'s, since the option of $f\fsubd g$ supported on
$\collar(v)$ can only occur when $\xi(g)=4$.
This gives us: 
\begin{equation}\label{number of fsubd}
|g|-1 \le \#\{h\fsubd g\} \le |g|+4.
\end{equation}

Now consider subsets $X\subseteq H$ with the
following property, for some fixed number $M$:
\begin{enumerate}
\item[(*)] If $g\notin X$ then there are at most $M$ geodesics
$h\fsubd g$ for which there exists $h'\in X$, $h'\fsubeq h$.
\end{enumerate}

The two main examples of $X$ satisfying this property are $X=H$
(trivially), and, of interest to us:

\begin{lemma}{X alpha property}
For any vertex $v$ in $\CC(S)$, $\alpha = v\srr$ or $v\sll$, and $k\ge
  4$, the sets $X_{\alpha}$ and $X_{\alpha,k+}$
satisfy property (*) with $M$ depending only on $S$.
\end{lemma}

\begin{pf}
Let $g\notin X_\alpha$. 
If $h'\fsubeq h\fsubd g$ then $D(h')\subset D(g)$, so if $h'\in
X_\alpha$ then $\alpha$ is in $\boundary D(h')$ and hence in $D(g)$,
where it must be nonperipheral since $g \notin X_\alpha$.
Now $D(h)$ is a component domain for some
simplex $u$ of $g$, and since $D(h')\subseteq D(h)$, we have
$d_{D(g)}(u,v) \le 1$. 
It follows that $u$ is restricted to an 
interval of diameter 2 in $g$, so there are at most 3
possibilities for $u$. The discussion leading to (\ref{number of
  fsubd}) then yields at most
$M=6$ possibilities for $h$.

Now to prove the property for $X_{\alpha,k+}$:
if $g\notin X_{\alpha,k+}$ then either $g\notin X_\alpha$, in which case the
estimate follows from property (*) for $X_\alpha$
together with the fact that $X_{\alpha,k+}\subset X_\alpha$, 
or $\xi(D(g)) < k$ in which case there is no $h\fsubeq g$ in
$X_{\alpha,k+}$, and the estimate is trivial.
\end{pf}


Recalling the notation ``$x\Qeq{a,b} y$'' from
\S\ref{del estimates}, we can state:
\begin{lemma}{Approximate lengths}
If $X\subseteq H$ satisfies (*) and $\varphi:X\to \R_+$ is a function
satisfying
$$\varphi(h) \Qeq{a,b} |h|$$
for all $h\in X$, then
$$
\sum_{h\in X} |h| \Qeq{A,B} \sum_{h\in X} \varphi(h)
$$
Where $A,B$ depend only on $a,b,S$ and the constant $M$ in (*).
\end{lemma}
The point of this lemma is that, although the additive errors $b$ can
accumulate when $X$ has many members, the additive errors at one level 
are swallowed up by the multiplicative constant at a lower level.

\begin{proof}
Define
\begin{equation}\label{beta define}
\beta(g) = \sum_{\substack{h\in X\\ h\fsubeq g}} |h|, \qquad
\beta'(g) = \sum_{\substack{h\in X\\ h\fsubeq g}} \varphi(h).
\end{equation}
We will show inductively for each $m\le \xi(S)$ that
$\beta(g) \Qeq{a',b'} \beta'(g)$ for
all $g$ with $\xi(D(g))\le m$, where $a',b'$ depend on $m$. The lemma
then follows from setting $m=\xi(S)$. The base case, $m = 2$, 
is immediate from the hypothesis $\varphi(g) \qeq |g|$, with
$(a',b')=(a,b)$.

Now note that, by
Theorem \lref{Descent Sequences},
whenever $f\fsub g$ there exists a unique $h$ with 
$f\fsubeq h \fsubd g$. This allows us to inductively decompose $\beta
$ and $\beta'$: 
\begin{equation}\label{decompose beta}
\beta(g) = \begin{cases}\displaystyle
|g| + \sum_{h\fsubd g} \beta(h) & g\in X \\
\displaystyle
\sum_{h\fsubd g}\beta(h) & g\notin X 
\end{cases}
\end{equation}
and similarly with $\beta'$ replacing $\beta$ and $\varphi(g)$
replacing $|g|$.

For both $\beta$ and $\beta'$, in case $g\in X$ there are at most
$|g|+4$ terms in the summation, by (\ref{number of fsubd}).
In case $g\notin X$, there are at most $M$ non-zero terms in the
summation, by property (*).

Now we can compare $\beta$ and $\beta'$. Suppose first that 
$g\in X$. We have:
\begin{align*}
\beta(g) &= |g| + \sum_{h\fsubd g} \beta(h) \\
\intertext{and by the inductive hypothesis:}
& \le |g| + \sum_{h\fsubd g} a'\beta'(h) + b' \\
\intertext{then by (\ref{number of fsubd}):}
& \le (1+b')|g|+4b' + a'\sum_{h\fsubd g}\beta'(h)\\
\intertext{and using $\varphi(g)\qeq|g|$,}
& \le (a\varphi(g)+b)(1+b')+ 4b' + a'\sum_{h\fsubd g}\beta'(h)\\
& =  a(1+b')\varphi(g) + a'\sum_{h\fsubd g}\beta'(h) + b(1+b')+4b'\\
& \le a''\beta'(g) + b''
\end{align*}
where the last line follows from the $\beta'$ version of
(\ref{decompose beta}), using $a'' = \max(a(1+b'),a')$ and $b'' =
b(1+b')+4b'$. 

In case $g\notin X$, we have:
\begin{align*}
\beta(g) &= \sum_{h\fsubd g} \beta(h) \\
& \le \sum_{\substack{h\fsubd g\\\beta(h)>0}} a'\beta'(h) + b' \\
\intertext{and by property (*):}
& \le b'M + a'\sum_{h\fsubd g}\beta'(h)\\
& = a'\beta'(g) + b'M.
\end{align*}

The inequality $\beta'(g) \qle \beta(g)$ is obtained in the same way
(with slightly different constants).
\end{proof}

\subsection*{Counting with top-level domains}
In the following proposition we show that the size of 
$X_{\alpha,4+}$ can be estimated by the size of $X_{\alpha,4}$.
This is in keeping with the intuition that the moves in the level 4
domains are the places where ``real'' change happens, and all the rest
is a bounded amount of bookkeeping.

\begin{proposition}{Counting 4-domains}
For any $x\in \CC_0(S)$ and $\alpha = x\sll$ or $x\srr$, 
we have
$$
\sum_{h\in X_{\alpha,4}} |h|
\qeq
\sum_{h\in X_{\alpha,4+}} |h|.
$$
\end{proposition}

\begin{proof}

The proof has the same inductive structure as the proof of Lemma
\ref{Approximate lengths}.
Let
$$
\gamma(g) = \sum_{\substack{h\fsubeq g\\ h\in X_{\alpha,4}}} |h|,
\qquad
\gamma'(g) = \sum_{\substack{h\fsubeq g\\ h\in X_{\alpha,4+}}} |h|.
$$

We shall inductively prove $\gamma(g)\Qeq{a,b} \gamma'(g)$ for constants
$a,b$ depending on $\xi(D(g))$. First if $\xi(D(g))\le 4$ then
$\gamma(g)=\gamma'(g)$ is obvious. It is also obvious that $\gamma\le
\gamma'$. Now assume $\gamma'(h) \Qle{a,b} \gamma(h)$ for
$\xi(D(h))<\xi(D(g))$ and let us prove it for $g$.

We first need a few lemmas. 

\begin{lemma}{interior at least 3}
If $\xi(D(h))\ge 4$ and 
$h\fsubd g$ at an interior simplex $v$, then
$|h|\ge 3$.
\end{lemma}
\begin{proof}
By assumption, the predecessor and successor $u$ and $w$ of $v$ must
be contained in $D(h)$ and by tightness of $g$, they fill $D(h)$.
Hence their distance in $\CC_1(D(h))$ is at least 3, and the
lemma follows. 
\end{proof}

\begin{lemma}{Lots in X}
Suppose $g\in X_\alpha$ and $\xi(D(g))>4$.
If $m\fsubd g$ at $v$ and $m'\fsubd g$ at $v'$, where $v'$ is the
successor of $v$ and both are interior  in $g$,  then at least
one of $m$ and $m'$ is in $X_{\alpha,4+}$. 
\end{lemma}
\begin{proof}
Note that $D(m)$ and $D(m')$ cannot be annuli since $v$ and $v'$ are
interior vertices and $\xi(D(g))>4$. Thus it suffices to show one of
them is in $X_\alpha$.
Suppose $m'\notin X_\alpha$, i.e. $\boundary D(m')$ does not contain
$\alpha$. Since $v'$ is not last or first, $v$ must be contained in
$D(m')$. It follows that $v'$ separates $v$ from $\alpha$. This means
that the component domain of $(D(g),v)$ that meets $\alpha$ must
contain $v'$. This domain is $D(m)$, so we conclude $m\in X_\alpha$.
\end{proof}

Finally we have:
\begin{lemma}{positive gamma}
If $h\in X_{\alpha,4+}$ and $|h|\ge 3$ then $\gamma(h)
\ge 3$. 
\end{lemma}

\begin{proof}
We proceed by induction. If $\xi(D(h))=4$ the statement is
obvious. Suppose $\xi(D(h))>4$.  Since $|h|\ge 3$ it has at least 
two interior
simplices, and by Lemma \ref{Lots in X} for at least one of them
there is a $k\fsubd h$ with $k\in X_{\alpha,4+}$.
By Lemma \ref{interior at least 3}, $|k|\ge 3$. Thus by induction
$\gamma(k)\ge 3$, and since clearly $\gamma(h)\ge \gamma(k)$, we are done.
\end{proof}

We return to the proof of Proposition \ref{Counting 4-domains}.
Recall that 
we are in the case $\xi(D(g))>4$, so in particular
\begin{equation}
  \label{gamma recursive}
 \gamma(g) = \sum_{h\fsubd g} \gamma(h).
\end{equation}
Now suppose that $g\in X_\alpha$. The recursive formula for $\gamma'$ 
gives
\begin{align}
\gamma'(g) &= |g| + \sum_{h\fsubd g} \gamma'(h).\notag \\
\intertext{By the inductive hypothesis this is}
&\le |g| + \sum_{h\fsubd g}a\gamma(h) + b. \notag \\
\intertext{Using (\ref{number of fsubd}), we obtain}
&\le 4b + |g|(b+1) + \sum_{h\fsubd g}a \gamma(h) \notag \\
&= 4b+ |g|(b+1) + a\gamma(g) \label{gamma' gamma}\\
\intertext{By (\ref{gamma recursive}). If $|g|\le 2$ this becomes}
&= a\gamma(g) + 2+6b \notag
\end{align}
and we are done.
Now suppose that $|g|\ge 3$. 
By Lemma \ref{Lots in X}, for least $\bigl\lfloor\frac{|g|-1}{2}\bigr\rfloor$ of
the interior simplices of $g$,  the corresponding
$h\fsubd g$ are in $X_{\alpha,4+}$. Applying Lemmas
\ref{interior at least 3} and  \ref{positive gamma},
we have $\gamma(h)\ge 3$ for those $h$.  This tells us that
\begin{align*}
\gamma(g) = \sum_{h\fsubd g} \gamma(h) & \ge  3 \#\{h\fsubd g: \gamma(h) \ge 3\}
\\
&\ge 3\frac{|g|-3}{2}
\end{align*}
Hence $|g| \le 3 + \frac23\gamma(g)$, so that (\ref{gamma' gamma}) becomes
\begin{equation*}
  \gamma'(g)  \le 7b+3 +  \left(\frac23(b+1)+a\right) \gamma(g). 
\end{equation*}

Finally if $g \notin X_\alpha$, we use property (*) as in the previous
section to argue 
\begin{align*}
\gamma'(g) &= \sum_{h\fsubd g} \gamma'(h)\\
&\le \sum_{\substack{h\fsubd g\\\gamma'(h) > 0}} a\gamma(h) + b\\
&\le Mb + a\gamma(g).
\end{align*}
This completes the inductive step, and establishes the proposition.
\end{proof}

\subsection{Comparing $\Homega$ and $\nomega$}

A uniform bound on the difference between real parts, $|\Re \Homega  -
\Re\nomega|$, follows directly from Lemma \ref{annulus signed
estimate}, which compares $[h_v]$ to $\tw_v(\I(H),\T(H))=\tw_v(\nu_-,\nu_+)$. 
Since the imaginary parts are at least 1, this gives us
\begin{equation}\label{nH real parts}
|\Re \nomega(v) - \Re \Homega(v)| = O(\Im \Homega(v)).
\end{equation}

If we can establish a bound of the form
\begin{equation}\label{nH imag parts}
\frac1c \le \frac{\Im\Homega(v)}{\Im\nomega(v)} \le c
\end{equation}
for some uniform $c$, then a bound on
$d_{\Hyp^2}(\nomega(v),\Homega(v))$ will follow.

It will  suffice to establish
\begin{equation}
  \label{Hn compare im}
  \sum_{Y \in \YY_{\alpha,4+}}
    \Tsh K{d_Y(\nu_-,\nu_+)}
    \qeq
    \sum_{h\in X_{\alpha,4}} |h|
\end{equation}
where $\alpha = v\sll$ or $v\srr$.
For convenience let $d_Y \equiv d_Y(\nu_-,\nu_+)$ throughout this
proof. 
Choose $K > M_2$, the constant in Lemma \ref{Large Link},
so that by that lemma if $d_Y \ge K$ then $Y$ is the support of some
geodesic in $H$. With this choice, we note that
\begin{align*}
  \sum_{Y\in\YY_{\alpha,4+}}
\Tsh K{d_Y}
  &=  \sum_{h\in X_{\alpha,4+}} \Tsh K{d_{D(h)}} \\
\intertext{Now since $\Tsh K{d_{D(h)}} \qeq d_{D(h)}$ and 
$d_{D(h)} \qeq |h|$ by Lemma \ref{Large Link}, and since
$X_{\alpha,4+}$ satisfies property (*) (Lemma \ref{X alpha property}), 
Lemma \ref{Approximate lengths} gives us}
  & \qeq \sum_{h\in X_{\alpha,4+}} |h| \\
\intertext{and then Proposition \ref{Counting 4-domains} gives us}
  & \qeq \sum_{h\in X_{\alpha,4}} |h|. \\
\end{align*}
This establishes (\ref{Hn compare im}), and since both sides of this
      estimate are positive and $\Im \nomega$ and
      $\Im \Homega$ are obtained from them by adding 1,
 (\ref{nH imag parts}) follows.

\subsection{Global projection bounds}
The same type of counting arguments that led to Theorem \ref{Omegas
close} can give us the following a priori bound on $|\omega_H|$ from
bounds on projections $d_Y(\nu_+,\nu_-)$. Define 
$\YY_v = \YY_{v\sll} \union \YY_{v\srr}$.

\begin{theorem}{projections bound omega}
Given end invariants $\nu$  and an associated hierarchy $H$,
for any internal vertex $v$ in $H$ we have
$$
|\Homega(v)| \le C\biggl(\sup_{Y\in\YY_v}
  d_Y(\nu_-,\nu_+)\biggr)^a
$$
where the constants $a,C>0$ depend only on $S$.
\end{theorem}

\begin{proof}
Let 
$$
B = \sup_{Y\in\YY_v}  d_Y(\nu_-,\nu_+).
$$
It suffices to find a bound of the form
\begin{equation}\label{beta bound}
\beta(h) \le c_j B^j
\end{equation}
where $j=\xi(D(h))-1$, and $\beta$ is defined as in (\ref{beta define}) with
$X=X_\alpha$ ($\alpha=v\sll$ or $v\srr$). Since $|\Re \Homega(v)|$ is
the $\xi=2$ term of $\beta(g_H)$ and $\Im\Homega(v)$ is the sum of
terms from $\beta(g_H)$ with $\xi>2$, the bound (\ref{beta bound})
applied to $\beta(g_H)$ gives us a bound of the form 
$|\Homega(v)| = O(B^{\xi(S)-1})$, which proves Theorem \ref{projections
bound omega}.

By choice of $B$, and Lemma \ref{Large Link}, we have
$$
|g| \le B+2M_1
$$ 
for each $g\in X_\alpha$. Let $c$ be the constant in (\ref{number of
fsubd}), and choose $c_1$ such that $B+2M_1 + c < c_1 B$. 
Then for $\xi(D(h))=2$ we have
$$
\beta(h) = |h| \le B+2M_1 < c_1 B
$$
if $D(h) = \collar(v)$ and $\beta(h) = 0$ otherwise. This establishes
the base case of (\ref{beta bound}).

Now assume (\ref{beta bound}) for $\xi(D(h))<\xi(D(g))$ and let us prove it for
$\xi(D(g))$. Let $j=\xi(D(g))-1$. Suppose first $g\in X_\alpha$. 
\begin{align*}
  \beta(g)  & =     |g| + \sum_{h\fsubd g} \beta(h) \\
  & \le B+2M_1  + (|g|+c) c_{j-1} B^{j-1} \\
  & \le B+2M_1 + (B+2M_1+c)c_{j-1} B^{j-1} \\
  & \le c_1 B + c_1c_{j-1} B^{j}\\
  & \le 2c_1c_{j-1} B^{j}.
\end{align*}
For $g\notin X_\alpha$, we have, using property (*), 
\begin{align*}
\beta(g) & = \sum_{h\fsubd g} \beta(h) \\
  & \le M c_{j-1} B^{j-1}.
\end{align*}
Thus by induction we have established (\ref{beta bound}), and the
theorem.
\end{proof}

\subsection{Non-internal vertices}
\label{noninternal omega}

If $v$ is not an internal vertex of $H$,
the definitions of the coefficients must be
adjusted somewhat before proving Theorem \ref{Omegas close}.
Let us first dispense with the cases {\em not} included in that
theorem, and which correspond to cusps in the hyperbolic 3-manifold:

\begin{itemize}
\item If $v$ is {\em parabolic} in $\I(H)$ or $\T(H)$,  then
$\storus v$ 
is an unbounded solid torus, and we define
$\Momega(v)=\Homega(v)=\nomega(v)=i\infty$.

\item A component of $\boundary S$ cannot be a vertex of $H$, but it
does have an associated solid torus $U$, and it is convenient
again to define all the coefficients to be $i\infty$.
\end{itemize}

We are left with the case that $v$ is a vertex of $\base(\I(H))$ which
{\em does} have a transversal, or similarly for $\base(\T(H))$ (or
both). 

In these cases $v$ appears as a curve of length at most $L_0$ in the
top or bottom (or both) of $\boundary_\infty N$. Let $r_+(v)$ denote
$1/\ep_1$ times the height of $\collar(v,\nu_+)$, in the rescaled metric $\sigma^{m'}$
(see \S\ref{augmented core}) that makes this annulus Euclidean. If $v$
does not appear this way in the top boundary then let $r_+(v)=0$.
Define $r_-(v)$ similarly. Now we can redefine $\nomega(v)$ and
$\Homega(v)$, by adding the term 
$$ i(r_+(v) + r_-(v))$$
to the expressions in (\ref{internal Homega}) and
(\ref{internal nomega}). The definition of $\Momega$ is unchanged from
Section \ref{define metric}.

The bound (\ref{Hn bound}) on 
$ d_{\Hyp^2}(\Homega(v),\nomega(v))$ is now immediate, since we have
added the same thing to both imaginary parts. 

In order to bound $ d_{\Hyp^2}(\Homega(v),\Momega(v))$,
we have to reconsider the sweeping discussion of Section \ref{sweeping}:

Suppose that $v\in \base(\T(H))$, and let $R$ denote the component of 
$R^T_+$ (see \S \ref{cores and ends}) that contains $v$. Thus there is
a top block $B=B\sbtop(\nu_R)$ associated to $R$. The interval $J'(v)$
cannot be defined as in \S\ref{sweeping} because there may be
infinitely many $i\in\II'$ such that $v$ is a component of $\eta_i$. However, 
as in \S\ref{embedding the model}, there is a first point $i\in\II'$
such that $F_i$ contains all of $\boundary_- B$. Let
$t=\max J'(v)$ be this value of $i$. The block $B$ is
then attached to $F_t$ independently of the rest of the sweep, and so 
we define $F_{t+1/2}$ to be $F_t \setminus \boundary_-
B\sbtop(\nu_R)$, union $\boundary_o B'$ (where we recall from
\S\ref{blocks and gluing} that $B'$ is $B$ minus the cusp annuli
$\boundary R\times[0,\infty)$), and $\boundary_o B' = R\times\{0\}$).
Thus $\boundary_o B$, in the isotopy class of $R$, can play the role of the
middle surface  $D(e_t)\times\{0\}$ in \S\ref{sweeping}. We project
$F_{t+1/2}$ to $Z_{t+1/2}\subset S$,  and define $\sigma'_{t+1/2}$ on
the $R$ component of $Z_{t+1/2}$ to be the projection of
$\sigma_\infty|_{\boundary_o B}$, the Poincar\'e metric associated to
$\nu_R$, as in \S\ref{define metric}. We define $\sigma_{t+1/2}$ as
before by extending across collars, though we may have to make a small
(uniform) adjustment to $\sigma'$ first near $\boundary R$, since
$\sigma_\infty$ has cusps rather than compact collars associated with
$\boundary R$. 

If $v\in\base(\I(H))$ we may use the bottom blocks to analogously
define $F_b$, $F_{b-1/2}$, 
$\sigma_{b-1/2}$, etc. (and otherwise we use the same definition as in
\S \ref{sweeping}). 

The shear computation that yields Lemma \ref{Shear bound} now goes through as
before -- the argument, and particularly the appeal to 
Lemma \ref{BG shear bound}, is insensitive to the fact that
$\collar(v)$ may now  
have large radius in $\sigma_{b-1/2}$ and/or $\sigma_{t+1/2}$.

In order to complete the comparison of $\Momega(v)$ and $\Homega(v)$ 
we need to consider changes to \S\ref{compare H M}, and particularly
to Lemma \ref{twist is real part}. If $v\in\base(\T(H))$ then 
the curve $w$ should be chosen to be the transversal of $v$ in the
marking $\T(H)$ -- recall that this means that $w$ is a minimal-length
curve crossing $v$ with respect to the $\sigma_\infty$ metric. Letting
$w^+$ be its geodesic representative in $\sigma_{t+1/2}$ and letting
$a_+$ be a component of $w^+\intersect \collar(v)$, we note that we
no longer have an upper bound on the length of $a_+$, but that its
twist in $\collar(v)$ with respect to a geodesic arc orthogonal to the
boundaries is at most 2, by the minimal-length choice of $w$. This
means that, lifting $\boundary U(v)$ to $\C$ and using the
notation of \S\ref{compare H M}, that the real part $b_2$ of the
vector associated to $a_2^+$ is still $O(1)$ (although the imaginary
part may be large). A similar argument applies when $v\in\base(\I(H))$
and we choose $u$ to be the transversal to $v$ from $\mu_-$.
The proof of Lemma \ref{twist is real part} then proceeds as 
before. 

Thus we obtain (\ref{MH real parts}) as before, namely
$$
|\Re \Momega(v) - \Re \Homega(v)| = O(\Im \Momega(v)).
$$ 
The equality (\ref{MH imag parts}), slightly adjusted to  
$$
\Im \Homega(v) = \Im \Momega(v) \pm O(1),
$$
follows as before, using the additional information that the top (resp.
bottom) annulus of $\boundary \storus v$ has width
$r_+(v)$ (resp. $r_-(v)$), and that these are the quantities added by
definition to $\Im \Homega(v)$.  Thus the estimate on
$d_{\Hyp^2}(\Homega(v),\Momega(v))$ follows as before.

Theorem \ref{projections bound omega}, which is stated for the
internal case, is easily generalized to yield
\begin{equation}\label{projections bound omega 2}
|\Homega(v)| \le C\biggl(\biggl(\sup_{\substack{Y\subset S\\v\in[\boundary Y]}}
  d_Y(\nu_-,\nu_+)\biggr)^a + r_+(v) + r_-(v)\biggr).
\end{equation}
The proof is exactly the same counting argument as for the internal
case, with the added imaginary part in the definition of $\Homega(v)$
yielding the extra terms.

\section{The Lipschitz model map}
\label{lipschitz}
We are now ready to prove the Lipschitz Model Theorem, whose statement
appears in the introduction.
We will build the map $f$ in several stages.

\subsubsection*{Step 0:} Let $f_0:\hhat S\times\R\to N$ be any map in
the homotopy class determined by $\rho$.
Via the embedding $\modl\subset \hhat S\times\R$
we restrict $f_0$ to $\modl$.

\subsubsection*{Step 1:} On each three-holed sphere $Y$ appearing as the gluing
boundary of a block, we find a homotopy from $f_0$ to $f_1$ so that $f_1$ is 
uniformly Lipschitz. This is done as follows:
There is 
a map $p_Y:Y\to N$,  homotopic to $f_0|_Y$,
and which is pleated on $\boundary Y$ (with the
usual proviso if a 
component of $\boundary Y$ corresponds to a cusp).
Let $\sigma_Y$ be the metric induced  by $p_Y$ on $int(Y)$, and let $Y_0 =
Y\setminus\collar(\boundary Y,\sigma_Y)$. Then, 
since Theorem \lref{Upper Bounds} gives us a uniform upper bound for the
$\sigma_Y$-lengths  of $\boundary Y$, there is a uniform $K_0$ 
such that  there is a
homeomorphism $\varphi_Y:Y \to Y_0$ homotopic to the identity, which is
$K_0$-bilipschitz with respect to the model metric on $Y$ and the
induced metric $\sigma_Y$ on $Y_0$. Now define
$$
f_1|_Y = p_Y \circ \varphi_Y.
$$
This is clearly a $K_0$-Lipschitz map, and homotopic to $f_0|_Y$. 

\subsubsection*{Step 2:} Extend $f_1$ to a map $f_2$, defined
on the middle surfaces of internal blocks, and
on the outer boundaries of  boundary blocks.

For any internal block $B=B(e)$, 
consider the middle subsurface
$W = D(e) \times \{0\}$ in $B$ (in its original identification as a
subset of $D(e)\times[-1,1]$).
Extend $\boundary W$ to a curve system $v_0$ which cuts $S$ into
components that are either 3-holed spheres, or $W$ itself. Then
$v^-=v_0\union e^-$ and $v^+=v_0\union e^+$ are pants decompositions
differing by an elementary move,
and we can define a {\em  halfway surface}
$g_{v^-,v^+}\in\pleat_\rho(v_0)$, as in 
Section \ref{pleated}. Let $\sigma_e$ be the metric induced by this
map on $int(W)$. Now, 
Proposition \lref{Upper Bounds}  gives us a uniform bound
$$\ell_\rho(x)\le D$$ where $x$ is any component of $\boundary W$, or
$e^\pm$, since these are all vertices in the hierarchy.
Lemma \ref{halfway bounds} then gives us uniform upper bounds 
in the induced metric $\sigma_e$:
$$
\ell_{\sigma_e}(e^\pm) \le D+C.
$$
This immediately gives lower bounds on 
$\ell_{\sigma_e}(e^\pm)$ as well, since the curves cross each
other. Furthermore it means that $W$ under $\sigma_e$ is not far from
being our ``standard'' hyperbolic structure on $W$. More precisely, 
letting $W_0=W\setminus \collar(\boundary W,\sigma_e)$
there is an identification
$\varphi_W:W\to W_0$
which is $K_1$-bilipschitz with respect to the model metric on $W$ and
$\sigma_e$ on $W_0$. 
Hence as in the previous paragraph, 
$$
f_2|_W = g_{v^-,v^+} \circ \varphi_W
$$
is the desired definition, and is clearly homotopic to $f_1|_W$.

\medskip 

Let $B$ be a boundary block associated to a subsurface $R$, 
and let us define 
$f_2$ on its outer boundary $\boundary_o B$. 
Recall that $\boundary_o B$ is identified with a component $R_\infty$
of $\boundary_\infty N$, and its metric is equal to the rescaled
Poincar\'e metric, $\sigma^m$. In Lemma \ref{infinity to aug core}
we describe a map $\hat p^{-1}\circ p_\infty$ that takes $R_\infty$
to the corresponding component of $\boundary \hhat C_N$, and is a
uniformly bilipschitz homeomorphism.

By composing this map with the identification $\boundary_o B \to
R_\infty$, we obtain a homeomorphism from $\boundary_o B$
to the corresponding component of $\boundary\hhat C_N$, satisfying
uniform bilipschitz bounds. This is our desired map $f_2$. 

\subsubsection*{Step 3:} We can extend $f_2$ to a continuous map $f_3$
defined on all of $\modl$, and still homotopic to $f_0$. This 
is possible by an elementary argument
because each of the surfaces where $f_2$ is defined so far is collared
in $\modl$. By appropriately defining the extension on these collars
using the homotopy from $f_2$ to the restriction of $f_0$, we can then
extend to the whole manifold. 

\subsubsection*{Step 4:} We  ``straighten'' $f_3$ on $\modl[0]$ to obtain a
map $f_4$, as follows: 

Fix an  internal block $B=B(e)$. 
Let $f_4$ agree with $f_3$ on the gluing boundaries $\boundary_\pm B$,
and on the middle surface $W=D(e)\times\{0\}$.

On $D(e)\times[-1/2,1/2]$, define $f_4 = f_3 \circ q$ where 
$q(y,t) = (y,0)$. This is certainly Lipschitz with a uniform
constant. 
Now consider a component of $\boundary_+ B$, which has the form
$Y\times \{1\}$ where $Y$ is a component of $D(e)\setminus
\collar(e^+)$.
$f_4$ is already defined on $Y\times\{1/2,1\}$, so for each $x\in Y$
extend $f_4$ to $\{x\}\times[1/2,1]$  as the unique constant speed
parameterization 
of the geodesic connecting the endpoints in the homotopy class
determined by $f_3$ (negative curvature is used here for the
definition to be unique, and for $f_4$ to be continuous and homotopic
to $f_3$). This extends $f_4$ to $Y\times[1/2,1]$.

To prove that $f_4|_{Y\times[1/2,1]}$ is uniformly Lipschitz it suffices
to find an upper 
bound on the lengths of the geodesics $f_4(\{x\}\times[1/2,1])$. This is
exactly the  ``figure-8'' argument of \cite[Lem 9.3]{minsky:torus}.
That is, since $f_4$ is already $K_1$-Lipschitz on $Y\times\{1/2,1\}$, 
there is a figure-8 $X\subset Y$, or bouquet of two circles, which is
a deformation retract of $Y$ with bounded tracks(see Figure \ref{fig8}), and such
that $f_4(X\times\{1/2\})$ 
and $f_4(X\times\{1\})$ have uniformly bounded lengths (depending on
$K_1$). If $f_4|_{X\times[1/2,1]}$ has a very long trajectory then at some
middle point $t\in[1/2,1]$  the two loops of $f_4(X\times\{t\})$ are
either both very short or nearly parallel to each other. In either
case discreteness, or J\o rgensen's inequality, is violated. 

For a component $Y\times\{-1\}$ of $\boundary_- B$ the definition 
of $f_4$ on $Y\times[-1,-1/2]$ is analogous.

\realfig{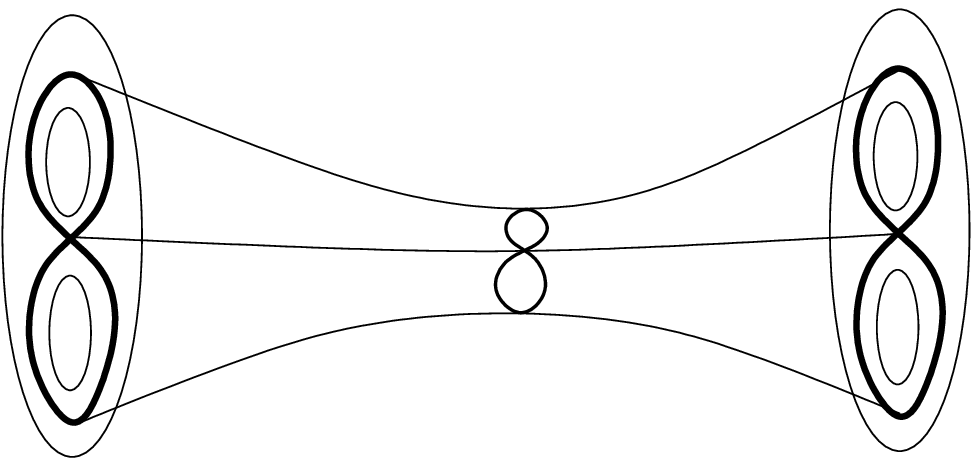}{A homotopy of a bounded figure-8 cannot
  be too long}

The very same discussion works for boundary blocks, except that the
outer boundary (minus its vertical annuli) takes the role of the
middle surface. 

The map $f_4$ is now defined on all of $\modl[0]$, 
uniformly Lipschitz (with respect to the
path metric on $\modl[0]$), and homotopic to $f_0$.

\subsubsection*{Step 5:}
We adjust $f_4$ slightly to obtain a map $f_5$ satisfying: 
\begin{enumerate}
\item $f_5(\modl[0])\subset \hhat C_N$
\item Whenever $\ell_\rho(v)\le \ep_1/2$ for a vertex $v$ of $H_\nu$, 
$f_5(\boundary U(v))  \subset \MT_{\ep_1}(v)$.
\end{enumerate}

$f_4$ restricted to the middle surfaces, gluing surfaces and
outer boundaries of boundary blocks already has image contained in
$\hhat C_N$. For internal blocks, the image is in fact already in
$C_N$ itself, and since the extension was done with geodesic arcs
$f_4$ takes each internal block into $C_N$ by convexity. The
$f_4$-image of a boundary blocks
could leave $\hhat C_N$, but only by a bounded 
amount because of the Lipschitz bound. In Section \ref{augmented core}
we show in particular that the geometry of $\boundary \hhat C_N$ is
relatively tame: it consists of parts of $\boundary C^1_N$ and
Margulis tube boundaries, meeting together with outward dihedral
angles of more than $\pi/2$. It follows that there is a uniformly
sized collar neighborhood of the boundary and a uniform Lipschitz
retraction of this neighborhood back into $\boundary C^1_N$. Composing
$f_4$ with this retraction yields a map with image in $\hhat C_N$. 

Let $v$ be a vertex with $\ell_\rho(v) < \ep_1$. 
The tube boundary $\boundary U(v)$ 
maps by a uniformly Lipschitz map and is homotopic into
$\MT_{\ep_1}(v)$. It can therefore be pushed into 
$\MT_{\ep_1}(v)$ by a uniformly Lipschitz map which comes from the
orthogonal projection to a lift of $\MT_{\ep_1}(v)$ to $\Hyp^3$.
Using this homotopy in a collar neighborhood of  $\boundary U(v)$, we
can adjust $f_4$ to get a map, still Lipschitz with a uniform
constant, that satisfies (2). This map is $f_5$.

\subsubsection*{Step 6:} We extend $f_5$ to all the non-parabolic
tubes in $\UU$. The resulting map $f_6$ is still homotopic to $f_0$.
Note that some continuous extension exists because $f_5$ is homotopic
to $f_3|_{\modl[0]}$.
To get geometric control we extend using a ``coning'' argument:
For each $U$ (which we
recall is isometric to a hyperbolic tube) we take a totally geodesic
meridian disk $\DD$. Foliating $\DD$ by geodesic segments emanating
from one boundary point,
we extend $f_6$ to the interior of $\DD$ to be totally geodesic on each of
these segments. $U\setminus \DD$
is now isometric to a
convex region in $\Hyp^3$ which we can again foliate by geodesics
emanating from a single point, and repeat. 

The extension has these properties: 
\begin{enumerate}
\item Whenever $\ell_\rho(v) < \ep_1$ for a vertex $v$ of $H$, 
$f_6(\storus v)$ 
is inside $\MT_{\ep_1}(v)$. 
\item  For any $k$ there exists $L(k)$ such that $f_6$
is $L$-lipschitz on $\modl[k]$.
\item Given $k>0$ there exists $\ep(k)\in(0,\ep_0)$ such that $f_6(\modl[k])$
  avoids the $\ep(k)$-thin part of $N$.
\end{enumerate}

For property (1), we apply property (2) of $f_5$ to see that
$f_5(\boundary U(v))\subset \MT_{\ep_1}(v)$ whenever $\ell_\rho(v) < \ep_1$.
The extension of $f_6$ to $U$ using geodesics is theorefore
also contained in $\MT_{\ep_1}(v)$.

For property (2), we use the fact that
if $|\omega_M(v)| < k$, then the hyperbolic tube $\storus v$ comes
from a compact set of possible isometry types. The boundary length of the
meridian disk, for example, is bounded by $k\ep_0$, and this together with
the negative curvature of the target and our method of extension by
coning gives us some Lipschitz bound $L$. Note that property (2) 
also implies part (5) of the statement of the Lipschitz Model Theorem.

Property (3) follows from property (2) and the following argument:
First note (following an argument of Thurston in \cite{wpt:II})
that through every point $x\in\modl[0]$ there is a pair of
loops $\alpha,\beta$ that generate a non-abelian subgroup  of
$\pi_1(S)$, and have uniformly bounded lengths. The images of these
therefore cannot be contained in a Margulis tube. This, together with
the Lipschitz bound on $f_5$, gives a uniform
$r_0>0$ so that $f_5(\modl[0])$ cannot penetrate more than $r_0$ into
the $\ep_1$-margulis tubes of $N$. Because of the $L(k)$-Lipschitz
property of $f_6$ on the tubes of $\UU\intersect \modl[k]$, and their
bounded geometric type, there is a uniform $r_1(k)$ so that
$f_6(\modl[k])$ cannot penetrate more than $r_1$ into the
$\ep_0$-Margulis tubes. This gives the desired $\ep(k)$,
by inequality (\ref{definite tube nesting}).

\subsubsection*{Step 7:} We next
adjust $f_6$ to obtain a map $f_7$
such that, for a certain
$k_1$, the image $f_7(\modl[k_1])$  avoids the interiors of the ``large'' tubes
$\MT[k_1]$ in $N$ (recall from the introduction that $\MT[k]$ is the
union of $\ep_1$-Margulis tubes in $N$, if any, corresponding to the homotopy classes
of the model tubes $\UU[k]$).
We will need the following lemma:

\begin{lemma}{Big omega short curve}
Given $\ep$, there exists $k'(\ep)$ such that if
$|\Momega(v)| \ge k'$ then $\ell_\rho(v)\le\ep$.
\end{lemma}

\begin{proof} 
Suppose first that $v$ is an internal vertex of $H$.

In \cite{minsky:kgcc} we prove that, if
$\ell_\rho(v)>\ep$, then
$$ d_Y(\nu_-,\nu_+) \le B$$
for each $Y\subset S$ such that $v\in[\boundary Y]$, where $B$ depends
only on $\ep$ and $S$. 
(This is a combination of Theorem B of \cite{minsky:kgcc} with the
proof of Theorem A.)
By theorem \ref{projections bound omega}, this implies
$$
|\Homega(v)|\le n,
$$
where $n$ depends only on the topology of $S$ and on $B$.
An upper bound on $|\Momega(v)|$ then follows from Theorem \ref{Omegas
close}, and this gives the $k'$ of the lemma. 

Suppose $v$ is non-internal. If it is parabolic then 
$|\Momega(v)|=\infty$ and $\ell_\rho(v)=0$, so the lemma continues to
hold.  If $v$ is not parabolic then we have the definition of
$\Momega(v)$ and $\Homega(v)$ from \S\ref{noninternal omega}, which
has an added term $i(r_+(v)+r_-(v))$. If $\ell_\rho(v)>\ep$ then we
obtain (via Bers' inequality \cite{bers:inequality}) a lower bound on
its length in the conformal boundary of $N_\rho$, which yields upper
bounds on $r_\pm(v)$. Combining this with the previous argument 
and the restatement (\ref{projections bound omega 2}) of
Theorem \ref{projections bound omega} for the noninternal case, 
we obtain an upper bound on $|\Momega(v)|$.
\end{proof}

Now let $k_0 = k'(\ep_1/2)$, let $\ep_4 = \ep(k_0)/2$ (where $\ep$ is
the function in property (3) of Step 5),
and let $k_1 = k'(\ep_4)$.
Consider a homotopy class of simple curves $v$ with $\ell_\rho(v) \le \ep_4$ --
In particular $v$ can be any vertex of $H$ with
$|\Momega(v)|\ge k_1$, by Lemma \ref{Big omega short curve}, but $v$
can also be a parabolic homotopy class, or a vertex of $\CC(S)$ which a-priori may
not appear in $H$. In the latter case let $\storus v = \emptyset$.

We claim that  $\modl \setminus \storus v$
must all be mapped outside of $\MT_{2\ep_4}(v)$ by $f_6$. For $\modl[k_0]$ this
follows from property (3) in Step 6 and
the choice of $\ep_4$. 

If $w\ne v$ is parabolic or a vertex of $H$ 
with $|\Momega(w)|\ge k_0$, by Lemma \ref{Big omega
short curve} we have $\ell_\rho(v)\le\ep_1/2$, and by Step 6, 
$\storus w$ is mapped into
$\MT_{\ep_1}(w)$, which is disjoint from $\MT_{\ep_1}(v)$. This
establishes the claim.

Since $\ell_\rho(v)\le \ep_4$, 
using the ``outward'' orthogonal projection 
of the tube $\MT_{\ep_1}(v)$ minus its geodesic core
to its outer boundary
$\boundary \MT_{\ep_1}(v)$, we can obtain a map
$\varphi_v:\MT_{\ep_1}(v) \to \MT_{\ep_1}(v)$ which is Lipschitz with
uniform constant, homotopic to the identity, and equal on the collar
$\MT_{\ep_1}(v)\setminus\MT_{2\ep_4}(v)$ 
to the outward orthogonal projection.
Let $\Phi:N_\rho \to N_\rho$ be equal to
$\varphi_v$ on $\MT_{\ep_1}(v)$ for each $v$ with $\ell_\rho(v)\le \ep_4$
and to the identity elsewhere. 

Let $f_7 = \Phi\circ f_6$.  This map takes each (non-parabolic) component of
$\UU[k_1]$ to its corresponding component of $\MT[k_1]$, and the
complement of all these to the complement of $\MT[k_1]$. 
It is Lipschitz with a uniform constant on $\modl[k_1]$.
$f_7$ is not
defined on the parabolic tubes, but their boundaries are mapped
(because of the retraction $\Phi$) to the boundaries of the
corresponding $\ep_1$-cusps.

\subsubsection*{Step 8:}
To obtain  $f_8$ we extend the definition of $f_7$ to the
parabolic tubes of $\UU$. Let $U$ be such a tube and $\MT$ its
corresponding $\ep_1$-cusp neighborhood. We already know that
$f_7$ maps $\boundary U$ to $\boundary\MT$. Let us identify $U$ with
$\boundary U \times [0,\infty)$ and $\MT$ with $\boundary\MT \times
[0,\infty)$. We can then define $f_8|_U$ with the rule
$f_8(x,t) = (f_7(x),t)$. This has the property that, if
$f_7|_{\boundary U}:\boundary U \to \boundary \MT$  is a proper map, 
then so is $f_8|_U:U\to\MT$.

$f_8 $ is our final map $f$. 
It takes each component of $\UU[k_1]$ to the corresponding Margulis
tube $\MT_{\ep_1}(v)$, and the complement $\modl[k_1]$ to the
complement of these tubes. On $\modl[k_1]$ it is still Lipschitz with
a uniform constant. 

\subsubsection*{Proper and degree 1}
We wish to show that $f$ is proper in both senses: that it takes 
$\boundary \modl$ to $\boundary \hhat C_N$, and that the inverse
images of compact sets are compact. 

By our construction (Step 2), if $\boundary \modl$ is nonempty
then $f$ maps it to $\boundary \hhat C_N$ by orientation-preserving
homeomorphism. Hence the first notion of properness holds. 

Let us show that $f|_{\modl[0]}$ is proper. If $\{x_i\in \modl[0]\}$
leaves every compact set, let us show the same for
$f(x_i)$.  Each internal block is compact, and the noncompact pieces
of the finitely many boundary blocks are the added vertical annuli, which by
construction are mapped properly. Hence we may assume that each $x_i$
is contained in a different internal block $B_i$. Writing $B_i = B(e_i)$ for a
4-edge $e_i$, let $v_i = e^+_i$. This is a sequence of distinct
vertices of $H$. Each $v_i$ has a representative $\gamma_i$ in $B_i$
of length bounded by a constant, and since $f$ is Lipschitz we obtain
a bound on the lengths of $f(\gamma_i)$. This means that $f(\gamma_i)$
must leave every compact set in $N$. Since the blocks have
bounded diameters, the same is true for $f(B_i)$ and hence $f(x_i)$. 

Now let $U_i$ be any infinite sequence of distinct, 
non-parabolic tubes in $\UU$. 
The generator curves of $\boundary U_i$ map to curves of bounded
length in $N$, in distinct homotopy classes, and hence must leave
every compact set. For tubes satisfying $|\Momega(U_i)| <  k_1$ there
is an upper bound on image diameter, and so their images must leave
every compact set as well. The images of  tubes with $|\Momega(U_i)|
\ge k_1$ are the corresponding Margulis tubes, which are all distinct
and hence again must leave every compact set. 

Let $U$ be a parabolic tube. The boundary $\boundary U$ is in
$\modl[0]$ and hence maps properly to the boundary of the
corresponding cusp tube in $N$. The extension that we constructed in
Step 7 is therefore automatically proper. 

We conclude that $f$ is proper. The boundary $\boundary \modl$, if
nonempty, is mapped by an orientation-preserving homeomorphism to
$\boundary \hhat C_N$, and in particular has degree 1. Hence in this
case $f$ has degree 1. 

In the remaining cases $\modl$ is all of $S\times\R$. Since $f$ is a
homotopy-equivalence, to show that it has degree 1 it suffices 
to show that it preserves the ordering of ends (determined by our
orientation convention in both domain and range).

Let $B_i$ be a sequence of blocks of $\modl[0]$ going to infinity in
$S\times\R$, in the positive direction (without loss of generality).
The associated 4-edges $e_i$ must eventually be subordinate to one of
the infinite geodesics of $H$, which terminates in the forward
direction in a lamination component $\nu_R$ of $\nu_+$. Thus the vertices
$v_i$ converge to $\nu_R$ in $\UML(S)$. By Thurston's theorem on
ending laminations (\S \ref{cores and ends}), the geodesic
representatives $f(\gamma_{v_i})^*$ in $N$ must exit the end of $N_0$
associated with $\nu_R$, which is a ``$+$'' end. Since $f(B_i)$ is a
bounded distance from either $f(\gamma_{v_i})^*$ or its Margulis tube,
the block images exit this end as well. Hence $f$ has degree 1. 

This establishes the Lipschitz Model Theorem, where $K$ is the uniform Lipschitz
constant obtained by this argument, and $k=k_1$.
\qed

\subsection*{The extended model map}

It is now a simple matter to prove the Extended Model Theorem, as
stated in the introduction. 

We note that the outer boundary of $\modl$ is naturally identified
with the 
``bottom'' boundary of $E_\nu$, namely $\boundary_\infty N
\times \{0\}$, and furthermore that the metric with which both these
boundaries were endowed is the same $\sigma^m$. Thus we can glue 
$\boundary \modl$ to the bottom boundary of $E_\nu$, obtaining the
extended model $\ME_\nu$. The model map $f$ restricted to $\boundary
\modl$ is also exactly the same as the map $\varphi$. Thus we can
combine the maps together to obtain a map $f':\ME_\nu \to N$. 
The desired properties of this map follow immediately from 
Lemma \ref{aug core exterior} and the Lipschitz Model Theorem.

\section{Length bounds}
\label{margulis}

We can now complete the proof of the Short Curve Theorem, stated in
the introduction. 

To prove part (1), let $\bar\ep = \ep_4$.  We saw in the proof of the
Lipschitz Model Theorem that $f(\modl[k_1])$ avoids the $\ep_4$-thin
part, and that each component of $\UU[k_1]$ maps into the Margulis
tube associated to its homotopy class. Since $f(\modl)$ covers all of
$\hhat C_N$ and all Margulis tubes are contained in $\hhat C_N$, 
it follows that all $\ep_4$-Margulis tubes are in fact covered by
tubes of $\UU[k_1]$.

Part (2) is exactly the statement of Lemma \ref{Big omega short
curve}. Thus it remains to prove Part (3), which is the content of
this lemma:

\begin{proposition}{length lower bound}
There exists a constant $c$ depending only on $S$, so that 
\begin{equation}\label{lambda lower bound}
|\lambda_\rho(v)| \ge \frac{c}{|\Momega(v)|}.
\end{equation}
for any vertex $v$ of $H$.
Furthermore the real part $\ell_\rho(v)$ satisfies
\begin{equation}\label{ell lower bound}
\ell_\rho(\gamma_v) \ge \frac{c}{|\Momega(v)|^2}.
\end{equation}
\end{proposition}

\begin{pf}
If $v$ is a parabolic vertex then $\lambda_\rho(v) = \ell_\rho(v) = 0$ and 
$|\Momega(v)| = \infty$, so the inequalities hold in the natural
extended sense. 

Let $v$ be a non-parabolic vertex of $H$. 
If $\ell_\rho(v) \ge \ep_4$ then, since $|\Momega(v)|$ is always at
least 1, our inequalities hold for appropriate choice of $c$.  

Assume now that $\ell_\rho(v) < \ep_4$.  Then as we saw in the proof
of the Lipschitz Model Theorem, 
the model map $f$ maps 
the complement of $\storus v$ to the complement of $\MT_{\ep_1}(v)$. 
Since $f$ has degree 1, it must take $\storus v$ to $\MT_{\ep_1}(v)$
with degree 1.  

It follows that the image of the meridian of $\storus v$ is a meridian
of $\MT_{\ep_1}(v)$. The meridian of $\storus v$  has length
$\ep_1|\Momega(v)|$ by definition of $\omega_M$, and the model map is
$K$-lipschitz. On the other hand we know that the meridian of a
hyperbolic tube of radius $r$ is at least $2\pi\sinh r$ (see Section
\ref{thick thin}). Hence
$$
2\pi \sinh r(v) \le K\ep_1|\omega_M(v)|.
$$
The Brooks-Matelski inequality (\ref{r lambda bound}) then implies
$$
2\pi \sinh \left(\log\frac{1}{|\lambda(\gamma_v)|} - c_1\right) \le K\ep_1|\omega_M(v)|
$$
and the Meyerhoff inequality (\ref{r ell bound}) gives
$$
2\pi \sinh \left(\half\log\frac{1}{\ell(\gamma_v)} - c_2\right) \le
K\ep_1|\omega_M(v)|. 
$$
A brief calculation yields (\ref{lambda lower bound}) and (\ref{ell
lower bound}).
\end{pf}


\providecommand{\bysame}{\leavevmode\hbox to3em{\hrulefill}\thinspace}

\end{document}